\documentclass[10pt]{article}
\usepackage{float}
\usepackage[english]{babel}
\usepackage{epsf,array}
\usepackage{latexsym,bbm}
\usepackage{epsfig,appendix}
\usepackage[]{algorithm2e}
\usepackage {amssymb,amsthm,amsfonts}
\usepackage {amsmath,comment,nccmath}
\usepackage{color,enumitem}
\usepackage{setspace}
\usepackage{url,caption,subcaption}
\usepackage{hyperref,appendix}
\usepackage{placeins}

\usepackage{rotating}
\usepackage[longnamesfirst]{natbib}
\bibpunct{(}{)}{;}{´a}{,}{,}
\renewcommand{\cite}{\citep}

\newtheorem{lemma}{Lemma}
\newtheorem{theorem}{Theorem}
\newtheorem{prop}{Proposition}
\newtheorem{corollary}[theorem]{Corollary} 

\allowdisplaybreaks

\makeatother
\makeatletter
\renewcommand{\@fnsymbol}[1]{\@arabic{#1}}
\makeatother

\newcommand{\BlackBox}{\ensuremath{\blacksquare}}

\title{Nonparametric spectral density estimation using interactive mechanisms under local
	differential privacy }

\author{Cristina Butucea\footnote{CREST, ENSAE, Institut Polytechnique de Paris, 5, avenue Henry Le Chatelier, 91120 Palaiseau, France} \and Karolina Klockmann \footnote{Universit\"at Kassel,
Heinrich-Plett-Str. 40, 34132 Kassel, Germany }\textsuperscript{2,3} 
\and Tatyana Krivobokova\footnote{ ISOR, Universit\"at Wien,
Oskar-Morgenstern-Platz 1, 1090 Wien, Austria } }
\date{ }

\begin{document}
\newpage
\maketitle

%%%%%%%%%%%%%%%%%%%%%%%%%%%%%%%%%%%%%%%%%%%%%%%%%%%%%%%%%%%%%%%%%%%%%
%ABSTRACT
%%%%%%%%%%%%%%%%%%%%%%%%%%%%%%%%%%%%%%%%%%%%%%%%%%%%%%%%%%%%%%%%%%%%%
\begin{abstract}
	\baselineskip=15pt \noindent 
	\\
	We study the problem of estimating the spectral density of a centered stationary Gaussian time series under local differential privacy constraints. Specifically, we propose new interactive privacy mechanisms for three tasks: recovering a single covariance coefficient, recovering the spectral density at a fixed frequency, and global recovery. Our approach achieves faster rates through a two-stage process: we first apply the Laplace mechanism to the truncated value, and then use the resulting privatized sample to learn about the dependence mechanism in the time series. 
	For spectral densities belonging to Hölder and Sobolev smoothness classes, we demonstrate that our algorithms improve upon the non-interactive mechanism of \cite{kroll2024nonparametric} for small privacy parameter $\alpha$, since the pointwise rates depend on $n\alpha^2$ instead of $n\alpha^4$. Moreover, we show that the rate $(n\alpha^4)^{-1}$ is optimal for estimating a covariance coefficient with non-interactive mechanisms. However, the $L_2$ rate of our interactive estimator is slower than the pointwise rate. We show how to use these procedures to provide a bona fide locally differentially private estimator of the entire covariance matrix. A simulation study validates our findings. 
\end{abstract}
\baselineskip=20pt

%%%%%%%%%%%%%%%%%%%%%%%%%%%%%%%%%%%%%%%%%%%%%%%%%%%%%%%%%%%%%%%%%%%%%
%INTRODUCTION
%%%%%%%%%%%%%%%%%%%%%%%%%%%%%%%%%%%%%%%%%%%%%%%%%%%%%%%%%%%%%%%%%%%%%
\section{Introduction}
In recent years, data privacy has become a fundamental problem in statistical data analysis. As the collection, storage, and analysis of personal data increase, the challenge lies in publishing accurate statistical results while protecting the privacy of individuals whose data is used. At the same time, traditional privacy techniques, such as
anonymization, are becoming increasingly ineffective due to the growing availability of publicly
accessible information, which facilitates the re-identification of individuals by linking anonymized
data with external data sets. To address these concerns, differential privacy has been proposed as a rigorous mathematical framework that offers strong privacy guarantees. First introduced by \cite{dwork2006differential}, differential privacy is typically classified into two main frameworks: global (or central) differential privacy (DP) and local differential privacy (LDP).

Global differential privacy relies on a centralized curator who collects confidential data from $n$ individuals and produces a privatized output derived from the complete data set. Only this privatized output is shared publicly, ensuring that modifying a single entry in the data set has minimal impact on the output's probability distribution. This makes it difficult for an adversary to determine whether a specific individual’s data is included in the data set.

In contrast, local differential privacy eliminates the need for a trusted curator. %, offering stronger privatization guarantees.
Under LDP, each individual independently generates a privatized version of their data on their own device before sharing it.  Consequently, only these privatized versions are available for analysis, ensuring that the original data remains private. While individuals maintain control over their original data, some level of interaction among the $n$ individuals may be permitted. LDP offers stronger privacy guarantees than DP at the expense of lower data utility for analysis.

In this paper, we focus on local differential privacy to privatize a single sequence of $n$ correlated observations $(X_1,...,X_n)\in \mathbb{R}^n.$  The ordering of the sequence may arise from the design of a study in which related individuals are surveyed sequentially, or from other factors, such as sensor networks where measurements are taken across different locations in a structured manner. 
Applications of this setting include snowball survey designs, in which individuals are interviewed sequentially, often introducing inherent correlations, whether spatial, social, or otherwise, among participants \citep{goodman1961snowball}. In such surveys, the variable of interest, $X_i$, might represent annual income or health-related attributes. 
In the following, we refer to $(X_1,...,X_n)$ as an observation of a time series and to $X_i$ as the data of the $i$-th individual, but the sequential structure may be caused by factors other than time.  
The proposed privacy mechanisms are also applicable for scenarios where the data $(X_1,...,X_n)$ belongs to a single individual, and sharing only a privatized aggregated statistic is insufficient.  For instance,  $(X_1,...,X_n)$ could represent an individual's energy consumption data. To optimize energy planning, providers may require detailed consumption patterns over time rather than just average usage. However, disclosing fine-grained consumption data could expose sensitive information, such as when a person is at home. This highlights the need to release a privatized version of the sequence $(X_1,...,X_n)$ that preserves global characteristics, such as the correlation structure, while preventing inference of individualized behavioral patterns. \cite{waheed2023privacy} proposed a generalized random response mechanism to enable the secure sharing of energy consumption data. Specifically, depending on the outcome of a coin flip, either $X_i$ is reported truthfully or replaced with a uniformly chosen value from $\{X_j\}_{j=1,...,i-1,i+1,...,n}$. This results in a new sequence that satisfies local differential privacy for an appropriately chosen coin-flip probability. \cite{shanmugarasa2024local} consider the privatization of a vector-valued time series representing the energy consumption of different household appliances. At each time instance, the energy consumption data is discretized to a set of energy levels and encoded as a binary vector, after which the randomized response mechanism of \cite{warner1965randomized} is applied to ensure local differential privacy.
While these works focus on algorithmic approaches, they lack a comprehensive analysis of how estimators of statistical quantities, such as the covariance function or the spectral density, are affected by the applied privacy mechanisms.

To formalize our setting, let $(X_t)_{t\in\mathbb{Z}}$ be a real-valued centered stationary Gaussian time series with (auto-)covariance function $\sigma: \mathbb{Z} \to \mathbb{R}, \, j \mapsto \sigma_j$, defined by $\sigma_j=\text{Cov}(X_t,X_{t+j})$ with arbitrary $t\in\mathbb{Z}$. We assume that the sequence of covariances $(\sigma_j)_{j\in\mathbb{Z}}$ is absolutely summable. This guarantees the existence of the spectral density function $f$, which is then defined by 
\begin{equation}
	\label{sdf}
	f(\omega)=\frac{1}{2\pi}\sum_{j\in\mathbb{Z}} \sigma_j\exp(-\texttt{i}j\omega), \quad \omega\in[-\pi,\pi],
\end{equation}
where $\texttt{i}$ is the imaginary unit.  In particular, $f$ is a non-negative, $2\pi$-periodic function which is symmetric at zero. Since the absolute summability of the sequence of covariances $(\sigma_j)_{j\in\mathbb{Z}}$ implies its square summability, the covariance function can be recovered from the spectral density function by the inverse Fourier transform. %Moreover, we assume that the sequence of covariances $(\sigma_k)_{k\in\mathbb{Z}}$ is square summable. Then, we can apply inverse Fourier transform and recover the covariance function from the spectral density function. 

We observe the sequence $(X_1,\ldots,X_n)$ drawn from the time series $(X_t)_{t\in\mathbb{Z}}$. In this paper, we are interested in the estimation of a single covariance coefficient $\sigma_j$ for some fixed integer $j$, of the spectral density $f(\omega)$ at a fixed frequency $\omega\in[-\pi,\pi]$, and of the whole function $f$ under local differential privacy constraints.

When designing locally differentially private mechanisms, two approaches can be distinguished: non-interactive and sequentially interactive mechanisms. In the non-interactive model, each individual independently generates a private representation $Z_i$ of their original data $X_i$. In contrast, the sequentially interactive model allows the $i$-th  individual to access previously privatized data $Z_1,...,Z_{i-1}$,  along with their own original data $X_i$ to create their private output $Z_i$.

The framework of differential privacy was introduced for i.i.d. data, and fewer methods exist for the privatization of a sequence of correlated data. 
The literature on privatizing such data can be broadly categorized into transform-based, model-based, and nonparametric approaches. In the transform-based approach, the correlated data is transformed, for example, with a Fourier or wavelet transform, so that the data become approximately uncorrelated \citep{rastogi2010differentially,lyu2017privacy,leukam2021privacy}.
Then, existing privacy mechanisms can be applied, and the data is backtransformed. However, this is only possible in the central setting, when access to the whole original data is available. 
Another line of literature employs parametric models to describe data correlation within the framework of differential privacy, and designs mechanisms specific to this type of correlation. For instance,  \cite{zhang2022differentially} consider a Gaussian $AR(1)$ model, where a privatized version $Z_i$ is generated as a convex combination of the true data $X_i$ and the linear minimum mean squared error estimate $\hat{X}_i$.  The estimate $\hat{X}_i$ is derived from $Z_{i-1}$ and estimators for the mean, variance, and autocorrelation based on previously released data $Z_1,...,Z_{i-1}$. By adding appropriately scaled, independent Gaussian noise, the mechanism achieves approximate local differential privacy. 
\cite{cao2018quantifying}  and \cite{xiao2015protecting} model temporal correlation with a Markov chain. % and study a scheme for privacy budget allocation over time. 
Specifically, \cite{xiao2015protecting} construct a Hidden Markov Model from the noise-perturbed data sequence and release the sequence inferred from it.  %HMM, known correlation structure
%\cite{zhu2014correlated}, %correlation used to define correlated query sensitivty to increase the noise
%\cite{wang2017cts}, %correlated Laplace noise with exactly the same correltaion structure
%\cite{zhang2022differentially}. %weighted AR(1) LDP, only eps-delt-LDP
Finally, there is a nonparametric approach that imposes fewer restrictions on the type of correlation.  To achieve global differential privacy, \cite{wang2017cts} proposed adding a sequence of correlated Laplace noise, with the correlation based on an estimate of the original sequence's autocorrelation. \cite{kroll2024nonparametric} considers linear stationary processes with sub-Gaussian marginal distributions and assumes that the spectral density belongs to a certain smoothness class. Local differential privacy is ensured by truncating each data point $X_i$ and adding independent Laplace noise. To our knowledge, \cite{kroll2024nonparametric} is the only work that considers minimax optimal estimation of the spectral density under local differential privacy.
For additional references on the privatization of sequential data, see the survey by \cite{zhang2022correlated}. Privatization of correlated, but non-sequential data, with multiple observations has been extensively studied in numerous works; see, for example, \cite{zhu2014correlated} and the references therein.
In the context of covariance estimation,   \cite{amin2019differentially} studied private estimation of the $d$-dimensional covariance matrix from a sample of $n$ i.i.d. vectors $Y_1,...,Y_n\in\mathbb{R}^d$, deriving error bounds under global differential privacy. Assuming a sub-Gaussian distribution for $Y_i$ and entry-wise sparsity of the covariance matrix, \cite{wang2021differentially} proposed an estimator for the covariance matrix from locally privatized data. \cite{wang2020tight} established the matching lower bounds, proving minimax optimality of this estimator. However, our setup corresponds to $n=1$, in which case these approaches and their corresponding convergence rates differ from ours. Indeed, the cited papers privatize the entire vector $Y_i \in\mathbb{R}^d$ rather than its components. Our results show that coordinate-wise local differential privacy converges faster as the vector dimension increases under stationarity.
Finally, alternative formulations of differential privacy that allow for correlated data, such as Bayesian differential privacy and the Pufferfish framework, have been introduced by \cite{kifer2014pufferfish} and \cite{yang2015bayesian}, respectively. \cite{mcelroy2023flip} introduce the concept of Linear Incremental Privacy and propose an all-pass filtering method to protect a univariate stationary time series. This approach preserves the second-order structure of the time series while also accounting for potential auxiliary information the attacker may possess.

\subsection{Our Contributions}
We propose new sequentially interactive privacy mechanisms and estimators of a single covariance coefficient $\sigma_j$, of the spectral density at a fixed frequency $\omega$: $f(\omega)$, and of the whole spectral density function $f$ from a sample of size $n$ generated from a centered Gaussian stationary time series.  We assume that the spectral density function belongs either to the class of Sobolev smooth spectral densities  $f\in W^{s,2}(L)$, $s>1/2,\, L>0$, or to the class of Hölder smooth spectral densities $f\in W^{s,\infty}(L_0,L)$, $s>1/2,\, L,\,L_0>0$, both explicitly defined in Section~\ref{sec:2}.
Previously, a non-interactive privacy mechanism for spectral density estimation under local differential privacy constraints was proposed by \cite{kroll2024nonparametric}.  Their estimator $\hat{f}^{\text{NI}}$ is shown to attain over the class $W^{s,2}(L)$ a convergence rate  (up to log-factors) of
$$\mathbb E \|\hat{f}^{\text{NI}}-{f}\|_2^2 \leq c(n \alpha^4)^{-\frac{2s}{2s+1}},$$
where $\alpha \in(0,1)$ is the privacy parameter which gets smaller for higher privacy level and $c>0$ is some constant depending on $s,\,L$. 
Their derived lower bounds are of the order $(n \alpha^2)^{-\frac{2s}{2s+1}}$, %\textcolor{red}{I don't think he had any $\alpha$ in the lower bounds?!} \textcolr{blue}{KK: Kroll has $ (e^\alpha-1)^2\approx \alpha^2$ in the lower bound}
leaving open the minimax optimal convergence rate.  For better comparability to our results, we first derive the convergence rate for $\sigma_j$ and $f(\omega)$ with the non-interactive privacy mechanism of \cite{kroll2024nonparametric}, which, as summarized in Table~\ref{tab:rates}, depend on $n \alpha^4$ as well. 

Next, we propose sequentially interactive mechanisms, also called interactive mechanisms here, and estimators  for a single Fourier coefficient $\sigma_j$ and the spectral density $f(\omega)$ at a fixed frequency.
For both, we show that our interactive mechanisms improve upon the currently known rates as $\alpha^4$ is replaced by $\alpha^2$. Specifically, we show the following bounds on the mean-squared error  for $\alpha\in(0,1)$ and sufficiently large $n$ (up to log-factors),
\begin{align*}
	\mathbb E|\hat{\sigma}_j-\sigma_j|^2\leq c \left (n \alpha^2 \right)^{-1}, &\quad \text{ if } \sum_{k\in\mathbb{Z}} \sigma_k^2 \leq M_1,\\
	\mathbb E|\hat{f}(\omega)-{f}(\omega)|^2 \leq c_1 (n \alpha^2)^{-\frac{2s}{2s+1}}, & \quad \text{ if }f\in W^{s,\infty}(L_0,L), \\
	\mathbb E|\hat{f}(\omega)-{f}(\omega)|^2 \leq c_2(n \alpha^2)^{-\frac{2s-1}{2s}}, & \quad \text{ if }f\in W^{s,2}(L),
\end{align*} where $c,c_1,c_2>0$ are constants depending only on $M_1>0$ and $s,\, L,\, L_0$, respectively.

It is important to note that these privacy mechanisms depend on the index $j$ and, respectively, on the frequency $\omega$, and thus cannot be directly generalized to estimate $f$ globally. 
Employing the mechanism of \cite{duchi2018minimax} to privatize vectors in an $\ell_\infty$ ball, we propose a sequentially interactive privacy mechanism and an estimator for the spectral density function $f$ and show the following bound on the mean integrated squared error, for sufficiently large $n$ (up to log-factors),
$$ 
\mathbb{E}\|\hat{f}-f\|_2^2\leq c_3 (n \alpha^2)^{-\frac{2s}{2s+2}}, \quad \text{if } f\in W^{s,2}(L) \text{ or } f\in W^{s,\infty}(L_0,L),
$$
where $c_3>0$ is a constant depending only on $s,\, L,\, L_0$.

This is the first setup in the literature where nonparametric rates differ for the pointwise mean-squared error and the mean integrated squared error. It is still an open question whether these rates are minimax optimal over all privacy mechanisms. Table~\ref{tab:rates} summarizes the rates attained using the non-interactive mechanism of \cite{kroll2024nonparametric} and the rates attained using our sequentially interactive privacy mechanisms.

{%
	\renewcommand{\arraystretch}{1.5} % Only affects this table
	\begin{table}[h] 
		\centering
		\begin{tabular}{|c|c|p{6cm}|c|}
			\hline
			$\alpha-$LDP mechanism& $\mathbb{E}| \hat{\sigma}_j- \sigma_j|^2$ & {\hspace{1cm} $\mathbb{E}|\hat{f}(\omega)-{f}(\omega)|^2$} & $\mathbb{E}\|\hat{f}-f\|^2_2$\\
			\hline
			Non-interactive    & $(n\alpha^4)^{-1}$  & $(n \alpha^4)^{-\frac{2s}{2s+1}}$ if $f\in  W^{s,\infty}(L_0,L)$ \newline
			$(n \alpha^4)^{-\frac{2s-1}{2s}}$ if $f\in  W^{s,2}(L)$ 
			& $(n \alpha^4)^{-\frac{2s}{2s+1}}$\\
			\hline
			Sequentially  &  $(n\alpha^2)^{-1}$ & $(n \alpha^2)^{-\frac{2s}{2s+1}}$ if $f\in  W^{s,\infty}(L_0,L)$  & $(n \alpha^2)^{-\frac{2s}{2s+2}}$ \\ 
			interactive     & & 
			$(n \alpha^2)^{-\frac{2s-1}{2s}}$ if $f \in  W^{s,2}(L)$  & $\wedge (n \alpha^4)^{-\frac{2s}{2s+1}}$\\
			\hline 
		\end{tabular}
		\caption{Convergence rates up to log-factors for recovering a single covariance coefficient $\sigma_j$, the spectral density $f$ at a fixed point $\omega$ and globally. Here, $a\wedge b:=\min\{a,b\}$.}
		
		\label{tab:rates}
	\end{table}
}
Using the close connection between the spectral density $f$ and the covariance matrix $\Sigma=(\sigma_{|i-j|})_{i,j=1,...,n}$ for stationary processes, we also provide an estimator $\check \Sigma$ and convergence rates with respect to the Hilbert-Schmidt norm under $\alpha$-LDP constraints.

These results confirm the importance of using sequentially interactive privacy mechanisms for achieving local differential privacy. Faster convergence rates with interactive mechanisms in time-series settings are intuitive, as past observations provide more valuable information than in i.i.d. settings.  A key reason for the improvement in our setting is that the quantities of interest are second‑order functionals. Non‑interactive mechanisms privatize each observation $X_i$ separately. As a result, products of privatized data lead to error terms involving products of independent privacy noises, whose variance scales with the fourth moment of the noise and yields an $\alpha^4$ rate. In contrast, interactive mechanisms allow privatizing statistics closer to the
target second-order functional, for instance the product of $X_i$ and previously privatized data. In this case, the dominant error depends only on the second moment of the privacy noise, resulting in an $\alpha^2$ rate.

Finally, we show that for estimating a single covariance coefficient $\sigma_j$, the rate $(n\alpha^{4})^{-1}$ is minimax optimal over all non-interactive privacy mechanisms and estimators. In order to do that, we prove a general information-theoretic bound for the transcripts $Z_{1:n}$ of the sensitive data through an arbitrary privacy mechanism. The starting point is that the Fisher information of the parameter (here it is $\sigma_j$) contained in a  sample $Z_{1:n}=(Z_1,...,Z_n)$ also satisfies a chain rule:
$$
I(Z_{1:n}) = I(Z_1) + I(Z_2|Z_1) + \ldots +I(Z_n|Z_{1:(n-1)}),
$$
where, for all $i$ from 1 to $n$, $
I(Z_i|Z_{1:(i-1)}) $ denotes the Fisher information of the conditional distribution, averaged over the conditioning variables.
Lemma \ref{FisherInfo} shows that for all $i=1,\ldots,n$ it holds that 
\begin{align*}
	I(Z_i|Z_{1:(i-1)})  \lesssim (e^\alpha -1)^2 I(X_i|Z_{1:(i-1)}) ,
\end{align*}
where, by convention, there is no conditioning when $i=1$. Moreover, we have that 
\begin{align*} 
	I(X_i|Z_{1:(i-1)})
	&\lesssim  (e^\alpha - 1)^2 I(X_i,X_{i-1}|Z_{1:(i-2)}) + I(X_i|Z_{1:(i-2)}) .
\end{align*}
This lemma holds for any privacy mechanism, whether interactive or not. It is, however, applied to non-interactive mechanisms in conjunction with a particular choice of distribution for the sensitive variables, in order to show the minimax lower bound of order $1/(n\alpha^4)$ for privately estimating a single covariance coefficient $\sigma_j$. This clarifies that non-interactive mechanisms cannot attain faster minimax rates. As we further explain in Section~\ref{sec:thmLB}, the above inequality shows that only a small amount of information propagates through a non-interactive mechanism, and this step does not extend to interactive mechanisms.

\section{Building Private Samples and Estimators} \label{sec:2}
Let $(X_1,...,X_n)$ be a sample from a real-valued, centered stationary Gaussian time series with sequence of covariances  $(\sigma_j)_{j\in \mathbb{Z}}$ and spectral density function $f:[-\pi,\pi]\to \mathbb{R}_{\geq 0}$ as defined in (\ref{sdf}).  %In particular, the spectral density is a $2\pi$-periodic and symmetric function. 
We assume that the spectral density function is smooth, that is, $f$ either belongs to the periodic Sobolev class of functions 
$$W^{s,2}(L){:=}\left \{ f(\cdot){=}\sum_{j\in\mathbb{Z}} \sigma_j \exp(-\texttt{i}j \cdot) :\, \sum_{j\in\mathbb{Z}} |\sigma_j|^2 |j|^{2s}\leq L\right \},$$
%$$W^{s,2}(L){:=}\left \{ f(\cdot){=}\sum_{j\in\mathbb{Z}} \sigma_j \exp(-\texttt{i}j \cdot) :\, f {\geq} 0 \text{ and } \sum_{j\in\mathbb{Z}} |\sigma_j|^2 |j|^{2s}\leq L\right \}$$
where $s>1/2,$ or to the periodic Hölder class of functions
%$$W^{s,\infty}(L_0,L){:=}\left \{ f(\cdot){=}\sum_{j\in\mathbb{Z}} \sigma_j \exp(-\texttt{i}j \cdot) :\, f {\geq} 0 \text{ and } |f^{(s-1)}(x){-}f^{(s-1)}(y)|{\leq} L \\, \forall\, x,y{\in} [-\pi,\pi] \right \} ,$$ 
%$$W^{s,\infty}(L_0,L){:=}\left \{ f\in W^{s,2}(L):\, |f^{(s-1)}(x){-}f^{(s-1)}(y)|{\leq} L_0|x-y| \\, \forall\, x,y{\in} [-\pi,\pi] \right \},$$ 
\begin{align*} 
	W^{s,\infty}(L_0,L){:=} \left \{ f(\cdot){=}\sum_{j\in\mathbb{Z}} \sigma_j \exp(-\texttt{i}j \cdot) :\, \right. & \|f\|_\infty<L, \, f^{(0)},...,f^{(\ell-1)} \text{ continuous and } \\
	& |f^{(\ell)}(x){-}f^{(\ell)}(y)|{\leq} L_0|x-y|^{s-\ell} \, \forall \, x,y{\in} [-\pi,\pi] \Bigg\},
\end{align*}
where $s>1/2,$  $\ell=\lfloor s\rfloor$  is the largest integer strictly less than $s$ and $L_0,L>0$ are some constants. The corresponding norms are denoted as $\|\cdot\|_2$ for the $L_2$ norm and $\|\cdot\|_\infty$ for the $L_\infty$ norm. Furthermore, we use the notation $a\lor b=\max\{a,b\}$ and $a\land b=\min\{a,b\}$ for $a,b\in\mathbb{R}$. To simplify, we sometimes omit constant factors and use  $\lesssim$ for less than or equal up to constants. The notation $a \asymp b$ is used  to indicate that $a$ and $b$ grow at the same asymptotic rate, i.e., there are constants $c_1,c_2>0$  such that $c_1b \leq a\leq c_2b$.

Note that, for $s>1/2$, $f\in W^{s,2}(L)$ implies that 
\begin{equation*} 
	\sum_{j \in \mathbb{Z}} \sigma_j^2 \leq L
\end{equation*} 
and for  $f\in W^{s,2}(L)$ or $ f\in  W^{s,\infty}(L_0,L)$, there exists a constant $M>0$ such that 
\begin{equation} \label{constants} \sum_{j\in\mathbb{Z}} |\sigma_j| <M,
\end{equation} 
%\begin{equation} \label{constants}\text{ } \, M>0 \text{  such %that } \sum_{j\in\mathbb{Z}} |\sigma_j| <M,
%\end{equation} 
see \cite{katznelson2004introduction}.

Next, we  introduce the definition of $\alpha-$local differential privacy (see, e.g., \citealp{duchi2018minimax}).
%{\bf Local differential privacy}
In this framework, private samples $Z_1,...,Z_n$ are successively obtained by  drawing from conditional distributions \( Q_i(\cdot \mid X_i, Z_1, \ldots, Z_{i-1}),\, i=1,...,n \). 
Let $\alpha \in (0,\infty)$. We say that
a sequence of conditional distributions $(Q_i)_{i=1,...,n}$ provides $\alpha-$local differential privacy   ($\alpha-$LDP)  or
that $Z_1, ..., Z_n$ are $\alpha-$local differentially private views of $X_1, ...,X_n$ if
\begin{equation} \label{def:aLDP}
	\sup_{Z_1, \ldots, Z_{i-1}, X_i, X_i'} \frac{Q_i(\cdot \mid X_i, Z_1, \ldots, Z_{i-1})}{Q_i(\cdot \mid X_i', Z_1, \ldots, Z_{i-1})} \leq e^{\alpha}, \quad \text{ for all } i = 1, \ldots, n,
\end{equation}
where we use the convention that \( \{Z_1, \ldots, Z_{i-1}\} \) is the empty set for \( i = 1 \). In particular, the ratio in Equation~(\ref{def:aLDP}) quantifies the sensitivity to changes in the input. By bounding this ratio, one ensures that the probability of observing the same private output does not differ significantly given the previously privatized data and the two different inputs $X_i$  and  $X_i'$.  The parameter $\alpha$ quantifies the privacy of the individual: the smaller the value of $\alpha$, the stronger the privacy guarantee. Note that $Z_i$ may also be a vector. As each \( Q_i \) is allowed to use previously released samples \( Z_1, \ldots, Z_{i-1} \), we call this privacy mechanism \textit{sequentially interactive}.

A subclass of LDP mechanisms is \textit{non-interactive}, which are not allowed to use previously released data. A non-interactive $\alpha-$LDP privacy mechanism generates each \( Z_i \) solely based on \( X_i \), via a conditional distribution of the form \( Q_i(\cdot \mid X_i) \) satisfying Equation~\ref{def:aLDP}. 

Local differential privacy is preserved under deterministic post-processing. Therefore, deterministic functions of the privatized samples, such as the estimators considered in this paper, satisfy the same privacy guarantee. Furthermore, estimators of nonnegative quantities may be replaced by their positive-part versions without affecting the privacy guarantees. For the estimation problems considered in this paper, this modification also does not affect the established risk bounds. %Hence, we omit the positive-part truncation in the subsequent analysis.

\subsection{Non-interactive Privacy Mechanism } \label{sec:1}
%We briefly summarize the non-interactive privacy mechanism and the spectral density estimator of \cite{kroll2024nonparametric}.
In this section, we give new results on the convergence rate of $\hat{\sigma}^{\text{NI}}_j$ and $\hat{f}^{\text{NI}}(\omega)$, which are estimators of $f(\omega)$ and $\sigma_j$ obtained with the non-interactive $\alpha-$LDP mechanism of \cite{kroll2024nonparametric}.
Let $(X_1,...,X_n)$ be a sample from a  Gaussian stationary time series $(X_t)_{t\in\mathbb{Z}}$ with  mean zero and a spectral density such that $f\in W^{s,\infty}(L_0,L)$ or $f\in W^{s,2}(L)$. Let $\tau_n>0$ be a threshold depending on the sample size $n$.
The original data $X_i$ is first truncated by setting
$$ \widetilde{X_i}=(X_i\land \tau_n)\lor(-\tau_n),\quad i=1,...,n.$$
Subsequently, it is perturbed by adding a Laplace distributed random variable $\xi_i$ with  $\xi_1,...,\xi_n\overset{\mathrm{i.i.d.}}{\sim} \mathrm{Lap}(2\tau_n/\alpha)$, yielding 
$$Z_i= \widetilde{X_i}+\xi_i,\quad i=1,...,n.$$ This mechanism is non-interactive and satisfies the definition of $\alpha-$LDP in Equation \eqref{def:aLDP}. Define  the bias corrected periodogram by $$\hat{I}_n(\omega)=\frac{1}{2\pi n}\left | \sum_{t=1}^n Z_t\exp(-\texttt{i}t\omega) \right|^2- \frac{4\tau_n^2}{\pi \alpha^2} .$$ Then, a spectral density estimator is defined by the partial Fourier sum  of order $m\in\mathbb{N}$
$$
\hat{f}^{\text{NI}}_m(\cdot)=\frac{1}{2\pi}\sum_{j=-m}^m \hat{\sigma}^{\text{NI}}_j \exp(-\texttt{i}j \cdot), \, \text{ where } \hat{\sigma}^{\text{NI}}_j =\int_{-\pi}^\pi \hat{I}_n(\omega) \exp(\texttt{i}j \omega)\, \mathrm{d}\omega.
$$
\cite{kroll2024nonparametric} showed that if $m\asymp \{1/n \lor \tau_n^4/(n\alpha^4)\}^{-\frac{1}{2s+1}}$ and $\tau_n^2=56 \log^{1+\delta}(n)$ for some constant $\delta >0$, then, for some constant $c>0$ and sufficiently large $n$

\begin{equation*}
	\mathbb E\|\hat{f}^{\text{NI}}_m-f\|^2_2 %\asymp \left(\frac{\tau_n^4}{n \alpha^4} \right)^{\frac{2s}{2s+1}}
	\leq c \left [\left\{\frac{1}{n} \lor \frac{\log^{2+2\delta}(n)}{n \alpha^4} \right\}^{\frac{2s}{2s+1}} \lor \frac{\log^{2+2\delta}(n)}{n} \right ].
\end{equation*} %see Theorem 3.2 in \cite{kroll2024nonparametric}.

We complement this result by proving the convergence rates for the estimation of the spectral density $\hat{f}^{\text{NI}}_m(\omega)$ at a fixed frequency $\omega$ and of a single covariance coefficient $\hat{\sigma}^{\text{NI}}_j$.
\begin{prop} \label{prop:ratekroll} Let $(X_t)_{t\in\mathbb{Z}}$ be a centered stationary Gaussian time series such that $\|\sigma\|_2^2=\sum_{k\in\mathbb{Z}}\sigma_{k}^2\leq M_1$ for some constant $M_1>0$. Define $\tau_n^2=56 \log^{1+\delta}(n)$ for some $\delta >0$. Then,  for sufficiently large $n$, the following holds 
	\begin{enumerate}
		\item[i)]  Let $R\in\mathbb{N}$  such that $R\leq\sqrt{n}$. For every $j\in\{-R,...,R\}$, it holds $$\sup_{\sigma:\, \|\sigma\|_2^2\leq M_1}\mathbb E|\hat{\sigma}^{\text{NI}}_j-\sigma_j|^2\leq c_1 \log^{2+2\delta}(n) \left \{\frac{1}{n} \lor \frac{1}{n\alpha^4} \right \},$$ where the constant $c_1>0$ depends on $M_1$.
		\item[ii)] If   $m\asymp ( \frac{1}{n} \lor \frac{\tau_n^4}{n\alpha^4})^{-\frac{1}{2s+1}}$, then for  every $\omega \in [-\pi,\pi]$, it holds that
		$$
		\sup_{f\in W^{s,\infty}(L_0,L)} \mathbb E|\hat{f}^{\text{NI}}_m(\omega)-{f}(\omega)|^2\leq c_2 \left [ \left\{\frac{\log(n)}{n} \lor \frac{\log^{3+2\delta}(n)}{n \alpha^4} \right\}^{\frac{2s}{2s+1}} \lor \frac{\log^{2+2\delta}(n)}{n} \right ],$$ where the constant $c_2>0$ depends on $s,\,L_0,L$ and $M$, with $M$ given in (\ref{constants}).
		\item[iii)] If $m\asymp ( \frac{1}{n} \lor \frac{\tau_n^4}{n\alpha^4})^{-\frac{1}{2s}}$, then for  every $\omega \in [-\pi,\pi]$, it holds that  $$\sup_{f\in W^{s,2}(L)} \mathbb E|\hat{f}^{\text{NI}}_m(\omega)-{f}(\omega)|^2 \leq c_3 \left [\left\{\frac{1}{n} \lor \frac{\log^{2+2\delta}(n)}{n \alpha^4} \right\}^{\frac{2s-1}{2s}} \lor \frac{\log^{2+2\delta}(n)}{n} \right ],$$
		where the constant $c_3>0$ depends on $s,\,L$ and $M$, with $M$ given in (\ref{constants}).
	\end{enumerate}
\end{prop}
In particular,  to estimate $\sigma_j$, only square-summability of the Fourier coefficients, but no further smoothness assumption on $f$ is required. The result for $\hat{\sigma}^{\text{NI}}_j$ holds for $j\in\{-(n-1),...,(n-1)\}$  fixed or slowly increasing with $n$ as long as $j=\mathcal{O}(n^{1/2}).$  If the normalization factor is replaced by $n/(n-j)$, that is, if $\hat{\sigma}^{\text{NI}}_jn/(n-j)$ is considered, then the result holds for $j=0,1,..,R$ with $0<R<n$ and $R=o(n)$. Note that the exponent $(2s-1)/(2s)$ in the pointwise convergence rate for a Sobolev smooth spectral density function is optimal in the classical nonparametric (non-private) setting.  

%%%%%%%%%%%%%%%%%%%%%%%%%%%%%%%%%%%%%%%%%%%%%%%%%%%%%%%%%%%%%%%%%%
%SUBSECTION: Interactive privacy mechanism for $\sigma_j$
%%%%%%%%%%%%%%%%%%%%%%%%%%%%%%%%%%%%%%%%%%%%%%%%%%%%%%%%%%%%%%%%%%
\subsection{Interactive Privacy Mechanism for \texorpdfstring{$\sigma_j$}{sigma	\textunderscore j}} \label{sec:sigmaj}
%Let  $K\in \mathbb{N}$ such that $K \ll n$. 
%\textcolor{red}{KK: previous extra parameter $K$ in this section is not needed and set to $j$ for estimating $\sigma_j$}
In this section, a new interactive privacy mechanism and an estimator for the covariance coefficient $\sigma_j$ are introduced. % with $j\in\{-R,...,R\}$, where $R\in\mathbb{N}$ such that $n>R>0$. 
Since  $\sigma_j=\sigma_{-j}$, it is sufficient to consider $j\in\{0,1,...,n-1\}$. Let $\xi_1,...,\xi_n\overset{\mathrm{i.i.d.}}{\sim}\mathrm{Lap}(4\tau_n /\alpha)$ be a set of Laplace random variables. If $j\neq 0$, then for the first $j$ iterations, $i=1,...,j$, 
the $i$-th data holder truncates its observation $X_i$  at a threshold $\tau_n>0$ and adds Laplace noise $\xi_i$, yielding
\begin{equation}
	\label{eq:Zi}
	Z_i=\widetilde{X}_i+\xi_i,  \, \, \text{ where }\,   \widetilde{X}_i= (X_i \land \tau_n ) \lor (-\tau_n ).    
\end{equation}
Consider a second set of Laplace  random variables: $\tilde{\xi}_{j+1},...,\tilde{\xi}_n\overset{\mathrm{i.i.d.}}{\sim}\mathrm{Lap}(4\tilde{\tau}_n /\alpha)$  for some  $\tilde{\tau}_n>0$, which is independent of $\xi_1,...,\xi_n$ and $X_1,...,X_n$. Starting from  iteration $i \geq j+1$, we build both $Z_i$ as before and a second random variable $\bar{Z}_{i,j}$ by %. For this, we first define   $$W_{i,j}=X_i \cdot Z_{i-j}.$$ Next, $W_{i,j}$ is trimmed at the level $\tilde \tau_n$ and we build  $\bar{Z}_{i,j}$ by adding the Laplace noise $\tilde{\xi}_i$, i.e.,
\begin{equation}
	\label{eq:barZij}
	\bar{Z}_{i,j} = \widetilde  W_{i,j} + \tilde \xi_i,  \, \, \text{ where }\,  \widetilde W_{i,j} = \left( W_{i,j} \wedge \tilde{\tau}_n \right) \vee (-\tilde{\tau}_n) \text{ and } W_{i,j}=X_i \cdot Z_{i-j}.
\end{equation}
In the case of $j=0$, we build directly, for all $i=1,\ldots, n$, the random variables
\begin{equation*}
	\label{eq:barZi0}
	\bar{Z}_{i,0} = \widetilde  W_{i,0} +  \xi_i, \, \, \text{ where } \,    \widetilde W_{i,0}=(X_i^2\wedge\tau_n ) \vee (-\tau_n)
\end{equation*}
and $\xi_1,...,\xi_n\overset{\mathrm{i.i.d.}}{\sim} \mathrm{Lap}(2{\tau}_n /\alpha)$.
The  proposed mechanism is sequentially interactive (for $j\neq 0$) and satisfies the definition of  $\alpha-$local differential privacy as shown in the following lemma.
\begin{lemma} \label{lemma:aLDP_sigma}
	Let $j\in\{0,1,...,n\}$. Then, for $i=1,...,j$,  $Z_i$ is an $\alpha-$local differentially private view of $X_i$, and for $i=j+1,...,n$,  $(Z_{i},\bar Z_{i,j})$ is an   $\alpha-$local differentially private view of $X_i$.	
\end{lemma}
The covariance coefficient $\sigma_j$ is then estimated by
\begin{equation} \label{est:sigmaj}
	\hat \sigma_j = \frac 1{n-j} \sum_{i=j+1}^n\bar{Z}_{i,j}.
\end{equation}
We show that the mean-squared error of $\hat{\sigma}_j$ is of the order  $1/(n\alpha^2).$ This is an improvement compared to the convergence rate in Proposition~\ref{prop:ratekroll} for the non-interactive privacy mechanism, which depends on $\alpha^4$. As in Proposition~\ref{prop:ratekroll},  to estimate $\sigma_j$, no smoothness assumption on $f$ is required. 
\begin{theorem} \label{theo:ratesigma} 
	Let $(X_t)_{t\in\mathbb{Z}}$ be a centered stationary Gaussian time series  such that $\|\sigma\|_2^2=\sum_{k\in\mathbb{Z}} \sigma_k^2\leq M_1$ for some constant $M_1>0$. Define $ \tau_n^2=8\log^{1+\delta}(n)$ and $\tilde \tau_n=16\log^{1+\delta}(n)\tau_n^2$  for some $\delta>0$. Let $R\in\mathbb{N}$  such that $0<R<n$ and $R=o(n)$. Then, for  $j\in\{0,1,\ldots,R\}$ and some constant $c>0$, the estimator $\hat{\sigma}_j$ of $\sigma_j$,  as defined in (\ref{est:sigmaj}), satisfies for sufficiently large $n$ that
	$$\sup_{\sigma:\, \|\sigma\|_2^2\leq M_1} \mathbb E|\hat{\sigma}_j-\sigma_j|^2\leq c \left \{ \frac{1}{n} \lor \frac{\log^{4+4\delta}(n)}{n \alpha^2} \right\},$$
	where the constant $c>0$ depends on $M_1$.
\end{theorem}
In particular, the result for $\hat{\sigma}_j$ holds for any $j\in\{-R,...,R\}$ fixed as long as $R = o(n)$.
%%%%%%%%%%%%%%%%%%%%%%%%%%%%%%%%%%%%%%%%%%%%%%%%%%%%%%%%%%%%%%%%%%
%SUBSECTION: Interactive privacy mechanism for $f(\omega)$
%%%%%%%%%%%%%%%%%%%%%%%%%%%%%%%%%%%%%%%%%%%%%%%%%%%%%%%%%%%%%%%%%%
\subsection{Interactive Privacy Mechanism for Pointwise Estimation} \label{sec:f0}
In this section, an  interactive privacy mechanism and an estimator for the spectral density function at a fixed $\omega\in[-\pi,\pi]$ are introduced. Let  $K\in \mathbb{N}$ be such that $K \ll n$. For the first $K$ iterations, we proceed as in Section~\ref{sec:sigmaj}, i.e., truncating $X_i$ with some threshold $\tau_n>0$ and adding an independent Laplace noise $\xi_i\sim \mathrm{Lap}(4\tau_n /\alpha)$ to obtain $Z_i$ in \eqref{eq:Zi}. Starting from iteration $i\geq K+1$, we build both $Z_i$ in \eqref{eq:Zi} and a second random variable $\widetilde Z_i$ as follows.   Let 
\begin{equation} \label{Vi}
	V_i=  X_{i}^2+\sum_{1\leq |k|\leq K} a_kX_{i}Z_{i-|k|}\exp(-\texttt{i}\omega k),
\end{equation}
where $a_k=1$ for $0\leq |k|\leq K/2$ and $a_k=2(1-|k|/K)$ for $K/2< |k|\leq K$.
Next, $V_i$ is trimmed at the level $\tilde \tau_n$, for some threshold $\tilde \tau_n>0$, and we build $\widetilde Z_i$ by adding a Laplace distributed random variable  $\tilde \xi_i$  with $\tilde\xi_{K+1},...,\tilde\xi_n\overset{\mathrm{i.i.d.}}{\sim} \mathrm{Lap}(4\tilde{\tau}_n/\alpha)$, yielding
\begin{align} \label{eq:tildeZi}
	\widetilde{Z}_{i}= \widetilde  V_{i}+ \tilde{\xi}_{i}, \, \, \text{ where }\,\widetilde  V_{i} = ( V_{i}\land \tilde{\tau}_n ) \lor (-\tilde{\tau}_n ).
\end{align}
%\begin{lemma} \label{lemma:aLDP_sdf0}
%	The privacy mechanism generating the sample $ Z_1,...,Z_K$ in \eqref{eq:Zi}, and $(Z_{K+1},\tilde Z_{K+1}),$ $\ldots, (Z_n,\tilde{Z}_n)$ in \eqref{eq:Zi} and \eqref{eq:tildeZi} is $\alpha-$local differentially private.	
%\end{lemma}
The  proposed mechanism is sequentially interactive and satisfies the definition of  $\alpha-$local differential privacy as shown in the following lemma.
\begin{lemma} \label{lemma:aLDP_sdf0}
	For $i=1,...,K$, $Z_i$ is an $\alpha-$local differentially private view of $X_i$. For $i=K+1,...,n$,  $(Z_{i},\widetilde Z_{i})$ is an $\alpha-$local differentially private view of $X_i$.	
\end{lemma}
The estimator of $f(w)$ is then defined by
\begin{equation*}
	%\label{est_pointwise}
	\hat f_K(\omega) = \frac 1{2\pi(n-K)} \sum_{i=K+1}^n \widetilde Z_i.
\end{equation*}
Note that the random variables $V_i = V_i(\omega)$ and $\widetilde Z_i = \widetilde Z_i(\omega)$ depend on the point of interest $\omega$. The next theorem shows that the estimator $\hat{f}_K(\omega)$ based on the sequentially interactive privatization mechanism attains the  usual nonparametric convergence rates in the non-private setting up to log-factors and a loss of $1/\alpha^2$. Compared to Proposition~\ref{prop:ratekroll}, the loss of $1/\alpha^4$ for non-interactive privacy mechanisms is improved to the smaller loss of $1/\alpha^2$ as $\alpha$ tends to 0.
\begin{theorem} \label{theo:ratef0} Let  $\omega \in [-\pi,\pi]$.
	Define $ \tau_n^2=8\log^{1+\delta}(n)$ and $\tilde{\tau}_n^2=1024\tau_n^{6}(K+1)$ for some $\delta>0$.  
	Then,  it holds for sufficiently large $n$ that
	\begin{enumerate}
		\item[$i)$]  If $K\asymp \{1/n \lor \tau_n^{6}/(n\alpha^2)\}^{-\frac{1}{2s+1}}$, then \begin{align*}
			&\sup_{f\in W^{s,\infty}(L_0,L)} \mathbb E|\hat{f}_K(\omega)-{f}(\omega)|^2 \leq c_1\left\{\frac{1}{n} \lor\frac{\log^{3+3\delta}(n)}{n \alpha^2} \right\}^{\frac{2s}{2s+1}},
		\end{align*}
		where the constant $c_1>0$ depends on $s,\,L_0,\,L$ and $M$, with $M$ given in (\ref{constants}).
		\item[$ii)$] If $K\asymp \{1/n \lor\tau_n^6/(n\alpha^2)\}^{-\frac{1}{2s}}$, then
		\begin{align*}
			&\sup_{f\in W^{s,2}(L)} \mathbb E|\hat{f}_K(\omega)-{f}(\omega)|^2 \leq c_2 \left\{\frac{1}{n}\lor\frac{\log^{3+3\delta}(n)}{n \alpha^2} \right\}^{\frac{2s-1}{2s}},
		\end{align*}
		where the constant $c_2>0$ depends on $s,L$ and $M$, with $M$ given in (\ref{constants}).
	\end{enumerate}       
\end{theorem}
Note that the weights $a_k$ in the definition of $V_i$ in Equation~\eqref{Vi} result in a reduced bias for Hölder-smooth functions. Consequently, the first factor in the rate in Theorem~\ref{theo:ratef0}$\,i)$ is $n^{-1}$ rather than $\log(n)n^{-1}$, as stated in Proposition~\ref{prop:ratekroll}$\,ii)$.

%%%%%%%%%%%%%%%%%%%%%%%%%%%%%%%%%%%%%%%%%%%%%%%%%%%%%%%%%%%%%%%%%%
%SUBSECTION: Interactive privacy mechanism for $f$
%%%%%%%%%%%%%%%%%%%%%%%%%%%%%%%%%%%%%%%%%%%%%%%%%%%%%%%%%%%%%%%%%%
\subsection{Interactive Privacy Mechanism for Global Estimation}
%Conjecture: $(n\alpha^2)^{-2s/(2s+2)}$.\\
There are two possible extensions of the methods in Section~\ref{sec:sigmaj} and \ref{sec:f0} to design an $\alpha-$LDP mechanism and an estimator for the whole spectral density $f$: either by building a $(K+1)$-dimensional vector of  $W_{i,j}, \, j=0,...,K,$ to estimate $K$ covariance coefficients $\sigma_0,...,\sigma_K$ in the $i$-th iteration. Or, a vector of $V_i(\omega_k), k=1,...,N$, is built in the $i$-th iteration to estimate the spectral density on  an $N$-dimensional grid $(\omega_k)_{k=1}^N$ of $[0,\pi]$ and then to smooth it. The choice of privatization of the vector is crucial. Using coordinate-wise Laplace noise, the variance scales with $K$, and a suboptimal exponent of $2s/(2s+3)$ in the first case and of $2s/(2s+4)$ in the second case can be derived. Therefore, we employ the mechanism of \cite{duchi2018minimax} to privatize vectors in an $\ell_\infty$ ball.
%In this section, we describe a privacy mechanism to obtain private data that is used to estimate the whole spectral density function $f$. 
For the first $K$ iterations, we proceed similarly as in Section~\ref{sec:sigmaj}, i.e., we truncate at threshold $\tau_n>0$ and add independent Laplace noise $\xi_i\sim \mathrm{Lap}(4\tau_n /\alpha)$ such that
$$Z_i=\widetilde{X}_i+\xi_i,  \, \, \text{ where }\,   \widetilde{X}_i= (X_i \land \tau_n ) \lor (-\tau_n ).$$ 
Starting from iteration $i\geq K+1$, we build both $Z_i$ as above and the vector \begin{align*}
	\widetilde W_i = (\widetilde W_{i,0},\ldots , \widetilde  W_{i,K})^\top \in [-\tilde \tau_n , \tilde \tau_n]^{K+1}, 
\end{align*} where 
\begin{align*}
	&\widetilde W_{i,j}= (W_{i,j} \land \tilde  \tau_n) \lor (-\tilde  \tau_n) \text{ and } \begin{cases}
		W_{i,j}=X_i\cdot Z_{i-j}, &j=1,...,K,\\
		W_{i,0}=X_i^2, &j=0,
	\end{cases}
\end{align*}
with some threshold $\tilde \tau_n >0$.
In particular, $\widetilde W_{i}\in B_\infty(\tilde{\tau}_n)$%\subset \mathbb{R}^{K+1}$
, where $B_\infty(r)=\{ x\in\mathbb{R}^{K+1} \mid \|x\|_\infty\leq r\}$ is the $\ell_\infty$-ball of radius $r$. % since $\|\tilde{W}_i\|_\infty \leq \tilde \tau_n$.
We apply the privacy mechanism for $\ell_\infty$ balls of \cite{duchi2018minimax}, corrected by \cite{butucea2023phase}. This mechanism applies to multivariate vectors whose sup-norm is bounded by $r$. The outcomes are the vertices of a hypercube with coordinates in $\{\pm B\}^{K+1}$ which are split in two sets according to the value of $\alpha$. \cite{duchi2018minimax} showed that a Laplace mechanism applied to each entry of the vector will introduce too much variance in the transcripts, a factor $K$ that makes all estimators based on these randomizations suboptimal. The hypercube mechanism we briefly describe here reduces the variance to its minimum and yields optimal estimators.

To obtain a private view $\check Z_i \in \mathbb{R}^{K+1}$ of $\widetilde W_i$ in the $i$-th iteration, we perform the following steps:
\begin{itemize}
	\item  If $\widetilde W_{i}=  (w_0,...,w_K)^\top$, define  a random vector $\widetilde Y_i=(\widetilde Y_{i,0},...,\widetilde Y_{i,K})^\top$ with
	$$ \widetilde{Y}_{i,j} := \begin{cases}
		+\tilde \tau_n, &\text{ with probability } \frac{1}{2}+ \frac{w_j}{2\tilde \tau_n},\\
		-\tilde \tau_n, &\text{otherwise.} \\%\text{ with probability } \frac{1}{2}- \frac{w_j}{2\tilde \tau_n}\\
	\end{cases} $$
	\item Sample $T_{i}\sim Ber(\pi_\alpha)$, where $\pi_\alpha=e^\alpha/(e^\alpha+1)$, and set 
	$$ \widetilde {Z}_i {:=} \begin{cases}  \text{Uniform} \left \{{z}{\in}\{\pm B\}^{K+1} \, \mid \, \langle  {z}, \widetilde{Y}_{i} \rangle{>} 0 \text{ or } \left(\langle {z}, \widetilde{Y}_{i} \rangle {=}0 \text{ and } \tilde z_1{=}  \frac{B}{\tilde \tau_n}\widetilde Y_{i,0}\right) \right \},  \quad \text{if } T_i{=}1,\\
		\text{Uniform} \left \{{z}{\in} \{\pm B\}^{K+1}\, \mid \, \langle {z}, \widetilde{Y}_{i} \rangle {<} 0 \text{ or } \left(\langle {z}, \widetilde{Y}_{i} \rangle {=}0 \text{ and }  z_1{=}\frac{-B}{\tilde \tau_n}\widetilde Y_{i,0} \right) \right \}, \text{ if } T_i{=}0,
	\end{cases}$$ where $B$ is chosen as
	\begin{equation*}
		\label{Kd}
		B= \tilde \tau_n \frac{e^\alpha+1}{e^\alpha-1}C_K, \text{ where } \frac{1}{C_K}=
		\begin{cases}
			\frac{1}{2^{K}}\binom{K}{\frac{K}{2}}, & \text{ if $K+1$ is odd},\\
			%\frac{1}{2^{K}+\frac{1}{2} \binom{K+1}{(K+1)/2}} \binom{K}{(K+1)/2}& \text{ if $K+1$ is even}.
			\frac{(K-1)!(K-1)}{2^{K}(\frac{K+1}{2}-1)!\frac{K+1}{2}!}, & \text{ if $K+1$ is even}.
		\end{cases}
	\end{equation*}
	%\begin{cases}
	%\frac{1}{2^{d-1}}\binom{d-1}{\frac{d-1}{2}} & \text{ if $d$ is odd}\\
	%\frac{(d-2)!(d-2)}{2^{d-1}(\frac{d}{2}-1)!\frac{d}{2}!} & \text{ if $d$ is even}.
	%\end{cases}
	\item Define the vector $\check{Z}_{i}$ by $\check Z_{i}=\widetilde{Z}_{i}$ if $K+1$ is odd, and by its components
	$$
	\check Z_{i,j} =
	\begin{cases}
		\frac{K-1}{2K}\widetilde{Z}_{i,0}, & \text{if } j=0,\\
		\widetilde{Z}_{i,j}, & \text{if } j = 1,\ldots ,K ,
	\end{cases}
	$$
	if $K+1$ is even.
	%\end{equation}
	
\end{itemize}
We show that this mechanism satisfies the definition of $\alpha-$LDP.

\begin{lemma} \label{lemma:EaLDP_f}
	The described mechanism has the following properties.
	\begin{itemize}
		\item[i)] For all $i=1,...,K$,  $Z_i$ is an $\alpha-$local differentially private view of $X_i$, and for all $i=K+1,...,n$, the sample $(Z_{i},\check Z_{i})$ is an $\alpha-$local differentially private view $(X_i,W_{i})$.
		\item[ii)]
		For all $i = K+1, \ldots, n$ and $j=0,...,K$, it holds
		\begin{align*}
			&\mathbb E[\check Z_{i}\mid X_i, Z_{1:(i-1)}, \check Z_{1:(i-1)}] =\mathbb E[\check Z_{i}\mid X_i, Z_{1:(i-1)}] = \widetilde W_i\\
			\text{and } \quad &\mathrm{Var}(  \check Z_{i,j}\mid X_i, Z_{1:(i-1)}, \check Z_{1:(i-1)} ) = \mathrm{Var}(  \check Z_{i,j}\mid X_i, Z_{1:(i-1)} ) \leq  B^2,
		\end{align*}
		where $B\lesssim \tilde \tau_n \sqrt{K+1}/\alpha$, if $\alpha\in(0,1)$.
	\end{itemize}
\end{lemma}

The estimator of $f$ is then defined by
\begin{equation*}
	\label{est_pointwise}
	\check f_K(\omega) = \frac 1{2\pi(n-K)} \sum_{i=K+1}^n  \sum_{0\leq |k|\leq K} \check Z_{i,|k|} \exp(-\texttt{i}k \omega), \quad \omega\in[-\pi,\pi],
\end{equation*} where $\check Z_{i,k}$ is the $k$-th entry of $\check Z_{i}$, for $k=0,...,K.$

Finally, we show that the convergence rate of this estimator depends on $\alpha^2$, but with a loss in the exponent. 
\begin{theorem} \label{theo:ratef}
	Define  $ \tau_n^2=8\log^{1+\delta}(n)$ and $\tilde \tau_n=16\log^{1+\delta}(n)\tau_n^2$ for some $\delta>0$.
	If $K\asymp\{\tilde \tau_n^{2}/(n\alpha^2)\}^{-1/(2s+2)}\land n^{-1/(2s+1)}$, then for sufficiently large $n$
	\begin{align*}
		&\sup_{f\in W^{s,\infty}(L_0,L)}\mathbb{E}\|\check{f}_K-f\|_2^2\leq c_1 \left( \frac{\log^{4+4\delta}(n)}{n \alpha^2} \right)^{\frac{2s}{2s+2}} \lor \left( \frac{1}{n}\right)^{\frac{2s}{2s+1}} \\
		\text{and} \quad &\sup_{f\in W^{s,2}(L)}\mathbb{E}\|\check{f}_K-f\|_2^2 \leq c_2 \left( \frac{\log^{4+4\delta}(n)}{n \alpha^2} \right)^{\frac{2s}{2s+2}}\lor \left( \frac{1}{n}\right)^{\frac{2s}{2s+1}},
	\end{align*} where the constant $c_1>0$ depends on $s,\,L$, and the constant $c_2>0$ additionally on $L_0$.
\end{theorem}

Combining this with the result of \cite{kroll2024nonparametric} gives the following upper bound on the $L_2$ risk over the class of Sobolev smooth spectral density function $W^{s,2}(L)$ for the low-privacy regime $\alpha\in (0,1)$:
$$
\left (\frac{\log^{4+4\delta}(n)}{n \alpha^2}\right)^{\frac{2s}{2s+2}} \land \left (\frac{\log^{2+2\delta}(n)}{n\alpha^4}\right)^{\frac{2s}{2s+1}},
$$
where the minimum depends on the magnitude of $\alpha.$ For small enough values of $\alpha=\alpha_n$, such that $n \alpha_n^{2(2s+3)} \log^{4s + 4s \delta}(n) $ tends to 0, then our interactive procedure achieves faster rates than the non-interactive procedure.

%\textcolor{red}{f uses squaKK: Can this be improved with a smaller dimension of the vectors to b privatized?  Idea: instead of $n-K$ estimates $\tilde W_{K+1,j},...,\tilde W_{n,j}$ for $\sigma_j$, try only $\lfloor (n-K)/\sqrt{K} \rfloor$ estimates for $\sigma_j$ by building only a vector $\tilde W_i$ of size $\sqrt{K}$ in each iteration. - Did not work, see comment in proof.} %Then, I expect instead $\frac{2M_1\sqrt{K}}{(n-K)}  + \frac{\tilde{B}^2}{n-K}$ with $\tilde{B}=\tilde \tau_n K^{1/4}$, which are balanced and I think yield the rate $(1/(n\alpha^2))^{2s/(2s+1.5)}$ as the bias of $f$ with $K$ Fourier coefficients is still $K^{-2s}$.}
The pointwise risk of the estimator $\check{f}_K$ is slower than the one of $\hat{f}_K(\omega)$ as shown in Theorem~\ref{theo:ratef0}.
\begin{corollary} \label{cor:ratef0}
%Define $ \tau_n^2=8\log^{1+\delta}(n)$ and $\tilde{\tau}_n^2=1024\tau_n^{6}(K+1)$ for some $\delta>0$ 
Define $\tau_n$ and $\tilde\tau_n$ as in Theorem~\ref{theo:ratef}. Then it holds for sufficiently large $n$ that
\begin{enumerate}
	\item[$i)$]  %If $K\asymp \{1/n \lor \tau_n^{6}/(n\alpha^2)\}^{-\frac{1}{2s+2}}$, then
	If $K\asymp\{\tilde \tau_n^{2}/(n\alpha^2)\}^{-1/(2s+2)}\land n^{-1/(2s+1)}$, then
	\begin{align*}
		&\sup_{f\in W^{s,\infty}(L_0,L)} \mathbb E|\check{f}_K(\omega)-{f}(\omega)|^2 \leq %c_1\left\{\frac{\log(n)}{n} \lor\frac{\log^{5+4\delta}(n)}{n \alpha^2} \right\}^{\frac{2s}{2s+2}},
		c_1\left(\frac{\log^{5+4\delta}(n)}{n \alpha^2} \right)^{\frac{2s}{2s+2}} \lor \left(\frac{\log(n)}{n}\right)^{\frac{2s}{2s+1}},
	\end{align*}
	where the constant $c_1>0$ depends on $s,\,L_0,\,L$ and $M$, with $M$ given in (\ref{constants}).
	\item[$ii)$] %If $K\asymp \{1/n \lor\tau_n^6/(n\alpha^2)\}^{-\frac{1}{2s+1}}$, then
	If $K\asymp\{\tilde \tau_n^{2}/(n\alpha^2)\}^{-1/(2s+1)} \land n^{-1/(2s)}$, then
	\begin{align*}
		&\sup_{f\in W^{s,2}(L)} \mathbb E|\check{f}_K(\omega)-{f}(\omega)|^2 %\leq c_2 \left\{\frac{1}{n}\lor\frac{\log^{4+4\delta}(n)}{n \alpha^2} \right\}^{\frac{2s-1}{2s+1}},
		\leq c_2 \left(\frac{\log^{4+4\delta}(n)}{n \alpha^2} \right)^{\frac{2s-1}{2s+1}} \lor \left(\frac{\log(n)}{n}\right)^{\frac{2s-1}{2s}} ,
	\end{align*}
	where the constant $c_2>0$ depends on $s,L$ and $M$, with $M$ given in (\ref{constants}).
\end{enumerate}     
\end{corollary}
Let us remark that an estimator for the Toeplitz covariance matrix $\Sigma=(\sigma_{|i-j|})_{i,j=1,...,n}$ can be easily derived from our procedure. The benefit of building an estimator of the spectral density function is that it becomes easier to produce a bona fide covariance estimator. For this, one must truncate the negative values of the estimator: $\check f_K^+:=\check f_K \lor 0$, and then apply the Fourier transform to obtain the estimators $\sigma^\dag_j,\,j=0,...,n-1$. Since $\check f_K^+$ is non-negative and symmetric, the resulting plug-in estimator $\check \Sigma=( \sigma^\dag_{|i-j|})_{i,j=1,...,n}$ is a proper positive semi-definite Toeplitz matrix.
By the relation of the Hilbert-Schmidt norm, denoted by $\|\cdot\|_{HS}$, and the $L_2$ norm, it follows
$$
\mathbb{E}\|\check \Sigma -\Sigma\|^2_{HS}\leq c  \left( \frac{\log^{4+4\delta}(n)}{n \alpha^2} \right)^{\frac{2s}{2s+2}}.
$$ 
%It is worth noting that estimating $\check{\sigma}_j= \frac{1}{n-K} \sum_{i=K+1}^n   \check Z_{i,|j|}$ directly from the $ \check Z_{i,|j|}$ may result in a non-negative definite plug-in estimator. 

\section{Non-interactive Mechanisms are Necessarily Slower}

The information theoretic inequalities by \cite{duchi2018minimax} were developed for independent initial random variables $X_1,...,X_n$. In our case, the dependence is at the core of the problem. Therefore, we develop in the next section a new information-theoretic inequality to hold for arbitrary multivariate likelihoods of the joint vector $X_{1:n}=(X_1,...,X_n)$ and any sequentially interactive (or non-interactive) mechanisms. The second section is an application of our inequality to show that our faster rates $1/(n\alpha^2)$ for estimating the covariance value $\sigma_j$ cannot be attained by any non-interactive mechanisms since the best they can do is $1/(n\alpha^4)$. We use the new inequality of the first section in conjunction with a particular construction of a Gaussian time series $X_1,\ldots,X_n$ in order to prove that only slower rates can be attained using non-interactive mechanisms.

\subsection{Information Theoretic Inequality for Dependent Data}

The Fisher information of the parameter $\theta$ contained in a  sample $Z_{1:n}=(Z_1,...,Z_n)$ also satisfies a chain rule
$$
I(Z_{1:n}) = I(Z_1) + I(Z_2|Z_1) + \ldots +I(Z_n|Z_{1:(n-1)}),
$$
with, for all $i$ from 1 to $n$,
$$
I(Z_i|Z_{1:(i-1)}) := \mathbb E_{Z_{1:(i-1)}} \left[ \int \frac{(\dot{\mathcal{L}}_\theta^{Z_i|Z_{1:(i-1)}} (z_i))^2}{\mathcal{L}_\theta^{Z_i|Z_{1:(i-1)}}(z_i) } \mathrm{d}z_i \right] .
$$

\begin{lemma} \label{FisherInfo}
For all $i=1,\ldots,n$, we have
\begin{align}\label{term1}
	I(Z_i|Z_{1:(i-1)})  \lesssim (e^\alpha -1)^2 I(X_i|Z_{1:(i-1)}) ,
\end{align}
where, by convention, there is no conditioning when $i=1$. Moreover, for an arbitrary $\epsilon>0$, 
\begin{align} \label{term2}
	I(X_i|Z_{1:(i-1)})
	&\leq (1+ \frac 1\epsilon) (e^\alpha - 1)^2 I(X_i,X_{i-1}|Z_{1:(i-2)}) + (1+\epsilon) e^{2\alpha} I(X_i|Z_{1:(i-2)}) .
\end{align}
and
\begin{align} \label{term3}
	I(X_i|Z_{1:(i-1)})  \leq (1 + \frac 94 (e^\alpha-1)^2) \left( I(X_i,X_{i-1}|Z_{1:(i-2)}) 
	+  I(X_i|Z_{1:(i-2)})
	\right).
\end{align}    
\end{lemma}
Note that the coefficients of the terms on the right-hand side in Equation \eqref{term2} are much smaller for the first term and decrease to 1 for the second as $\alpha >0$ tends to 0, whereas in \eqref{term3} the coefficients are balanced and decrease to 1 as $\alpha$ tends to 0. We think both versions may be useful in different contexts.

\subsection{Non-interactive Mechanisms are Slower for Estimating a Covariance Coefficient} \label{sec:thmLB}
%\textcolor{blue}{KK: Title: Minimax lower bound for estimating a covariance coefficient under non-interactive mechanisms ? - I think a more informative title is preferable, as suggested by the referee}
By the results of \cite{kroll2024nonparametric}, $1/(n\alpha^2)$ is a lower bound for estimating $\sigma_0$. Here, we show that for estimating $\sigma_j$ for some fixed $j\geq 1$ using non-interactive privacy mechanisms, the lower bound  increases to $1/(n\alpha^4)$, showing optimality of the plug-in estimator of \cite{kroll2024nonparametric} based on a non-interactive privacy mechanism. However, we recall that our interactive procedure and the associated estimator attain the faster estimation rate $1/(n \alpha^2)$.

Let us denote by $\mathcal{Q}^{\mathrm{NI}}_\alpha$ the class of non-interactive privacy mechanisms which are $\alpha-$LDP. We assume here that $0< \alpha \leq \alpha_0 <1$ for some fixed $\alpha_0$. 
\begin{theorem}\label{thm:LBcoeff} Assume $X$ is an $n-$dimensional vector generated from a centered stationary Gaussian time series with $\|\sigma\|_2^2=\sum_{j \in \mathbb{Z}} \sigma_j^2 \leq M_1$. Then, for any arbitrary, fixed $K\geq 1$, we have
$$
\inf_{Q \in \mathcal{Q}^{\mathrm{NI}}_\alpha} \inf_{\hat \sigma_K} \sup_{\sigma:\|\sigma\|_2 ^2\leq M_1} \mathbb{E}\left[ |\hat \sigma_K - \sigma_K|^2\right] \geq \frac{C}{n \alpha^4},
$$
for sufficiently large $n$. 
\end{theorem} 
%\textcolor{blue}{ $\inf_{Q \in \mathcal{Q}^{\mathrm{NI}}_\alpha} \inf_{\hat \sigma_j} \sup_{\sigma:\|\sigma\|_2 ^2\leq M_1}$ }
Thus, the faster rate $1/(n\alpha^2)$ that our interactive method attains cannot be obtained by any non-interactive method.

\emph{Sketch of proof of Theorem~\ref{thm:LBcoeff}}: We proceed by reducing the set of parameters to a simple $K-$dependent centered and stationary Gaussian vector, that is, only $\sigma_0$ and $\sigma_K$ are non-zero. We let $\sigma_K$ vary in a small interval and apply the Van Trees inequality to obtain the first bound from below. The bound involves the Fisher information that the privatized vector $Z_{1:n}$ brings on the parameter of interest. After the telescopic expansion of this Fisher information for conditional likelihoods, Lemma~\ref{FisherInfo} provides a one-step contraction relating $I(Z_i|Z_{1:(i-1)})$ to $I(X_i|Z_{1:(i-1)})$ under local differential privacy. In order to obtain an $\alpha^{-4}$ lower bound for non-interactive mechanisms, we additionally control the latter information by using the fact that each private sample $Z_i$ only depends on the sensitive $X_i$ at each time step $i$ and not on the previous transcripts $Z_1,...,Z_{i-1}$. This restricts how much information about the dependence parameter can propagate through the transcript and yields an additional $\alpha^2-$order contraction. This additional step does not apply to interactive mechanisms because $Z_i$ may adapt to past releases, thereby making the transcript more informative about the correlation parameters. 

\section{Simulation Study} 
\label{sec:simulation}
We compare the finite-sample performance of the proposed sequentially interactive (SI) private estimators with the non-interactive (NI) private estimators of \cite{kroll2024nonparametric} for estimating a covariance coefficient and the spectral density at a fixed frequency.
In particular, we focus on the dependence of the mean squared error (MSE) on the privacy level $\alpha$ for a fixed, sufficiently large sample size $n$, to validate the privacy–accuracy trade-off predicted by our theoretical results in the previous sections.
We consider three mean-zero Gaussian stationary processes $(X_1,\dots,X_n)$ with $n=1000$ and variance $1.44$:
\begin{enumerate}
\item[(1)]  an autoregressive process AR($0.8$) with Gaussian innovations and variance $1.44$, 
\item[(2)] a Gaussian process with spectral density $f(x)= 0.2\{|\cos(x)|^{0.8}+0.45\}$ for $x\in\mathbb{R}$,
\item[(3)] a Gaussian process with covariance coefficients $\sigma_k= 1.44(1+|k|)^{-5.1}$ for $k=0,...,n-1$.
\end{enumerate}
The three examples are chosen to represent Gaussian stationary processes with substantially different dependence structures. 
Example~(1) has exponentially decaying covariance coefficients and an analytic spectral density, Example~(2) has slowly decaying covariances $\sigma_k=\mathcal{O}(k^{-0.8})$ and a $0.8$-Hölder‑continuous spectral density, %$f\in W^{0.8,\infty}(2,2)$, 
and Example~(3) exhibits faster polynomially decaying covariances and a moderately smooth spectral density $f\in W^{4,\infty}$.

We estimate the variance $\sigma_0$, the lag-$2$ covariance coefficient $\sigma_2$, and the spectral density at frequency $\omega=\pi/5$. Performance is measured by the Monte Carlo estimate of the mean squared error, computed over $300$ independent replications. For spectral density estimation, we use
$
K =\lceil 3 \{ n \,\lor\, (n\alpha^2)/\tau_n^{6}\}^{1/(2s+1)} \rceil
$ in the sequentially interactive case and  $
m = \lceil 3 \{ n \,\lor\, (n\alpha^4)/\tau_n^{4}\}^{1/(2s+1)}\rceil$ in the non-interactive case,
where $\lceil \cdot \rceil$ denotes the ceiling function. We set $s=3$ and $\delta=0.001$ for all examples, as the influence is negligible, provided $\delta$ is chosen small. Negative variance estimates are replaced by $10^{-3}$, and negative spectral density estimates are truncated to zero.
For comparison, we furthermore estimate the  MSE of the corresponding standard estimators on the original data, which are $\hat{\sigma}_j=\frac{1}{n}\sum_{k=1}^{n-j}X_kX_{k+j}$ and $f(x)=\frac{1}{2\pi}\sum_{j=-K}^K$ $ \hat{\sigma}_j \exp(-\texttt{i}j x) $ with $K=  \lceil 3 n^{1/(2s+1)}\rceil$. %For  Example (1), the true quantities are $\sigma_0=1.44$ and $\sigma_2=0.352$ and $f(\pi/5)= 1.51$.\\

We present the results for Example (1) in detail. The remaining examples are reported in the Appendix and lead to similar conclusions.
Figure~\ref{EX1_loglog_plot_MSE_varyalpha}  displays the logarithm of the estimated MSE as a function of $\log(\alpha)$ for NI-based estimators (left) and SI-based estimators (right), for $\sigma_0,\, \sigma_2$ and $f(\pi/5)$. The logarithm of the theoretical MSE in the privacy‑dominated regime is overlaid as dashed gray lines with slopes $-4$ for NI-based estimation and $-2$ for SI-based estimation.

Using the theoretical truncation levels $\tau_{\mathrm{theory}}$ and $\tilde{\tau}_{\mathrm{theory}}$ as defined in  Proposition~\ref{prop:ratekroll}, Theorem~\ref{theo:ratesigma}, and Theorem~\ref{theo:ratef0}, the logarithm of the empirical MSE closely follows the predicted slopes for small and moderate values of $\alpha$. As $\alpha$ increases further, the MSE  coincides with the MSE of the corresponding non-private estimator. 
Overall, the simulations confirm the theoretical scaling of the MSE with respect to $\alpha$.

For both NI-based and SI-based estimators, the theoretical truncation levels tend to be too large in finite samples, leading to very large MSE for strong privacy (small $\alpha$).
In practice, moderate reductions of $\tau_{\mathrm{theory}}$ and $\tilde{\tau}_{\mathrm{theory}}$ can substantially reduce the MSE for strong privacy levels without affecting the $\alpha$-rate. For  NI-based variance estimation, reducing $\tau_{\mathrm{theory}}^{\mathrm{NI}} \approx 20$ to $\tau \approx 3$ uniformly improves or preserves performance across $\alpha$. In contrast, more aggressive truncation, that is, choosing $\tau$ close to or smaller than the true parameter value, introduces a truncation bias and leads to an increase in the MSE for larger values of $\alpha$. For estimation of $\sigma_2$ and $f(\pi/5)$ from NI‑privatized data, similar qualitative behavior is observed. For SI-based variance estimation, the theoretical truncation level is close to the boundary at which truncation bias becomes visible. Slightly smaller truncation values reduce the MSE for $\alpha\in(0,1)$ but lead to a bias for larger $\alpha$. %This indicates that the theoretical truncation level is close to sharp for SI‑based variance estimation.
For estimation of $\sigma_2$ and $f(\pi/5)$ from SI‑privatized data, a moderate simultaneous reduction of the truncation parameters $\tau$ and $\tilde{\tau}$ leads to a noticeable decrease of the MSE over a wide range of $\alpha$, while overly strong truncation again introduces a bias for larger $\alpha$. No substantial difference in sensitivity between $\tau$ and $\tilde{\tau}$ is observed in this setting.

Across all three estimation problems, the SI-based estimators achieve substantially lower MSE than their NI counterparts in the strong privacy regime, consistent with the improved scaling in $\alpha$. For variance estimation, the SI-based estimator outperforms the NI-based estimator uniformly over $\alpha$, whereas for covariance and pointwise spectral density estimation, the advantage is most pronounced for $\alpha\in(0,1)$, with the precise crossover depending on the exact choice of the truncation parameters.

Increasing the sample size $n$ yields similar observations regarding the alignment of estimated MSE with the theoretical results and the influence of the truncation parameter.  If the variance of the underlying process increases, the true values of $\sigma_0,\sigma_2$ and $f(\pi/5)$ increase, and truncation bias becomes relevant earlier.

In  Example~(2) and (3), where $\sigma_2$ is close to zero,  
aggressive truncation reduces the variance of the estimator while introducing only negligible bias. As a result, the MSE is smaller than that of the non-private estimator for moderate and large $\alpha$. Otherwise, the results for Examples~(2) and (3) are qualitatively similar to those of Example~(1), despite their different dependence structures. The simulations suggest that, in finite samples,  truncation and privacy noise dominate the magnitude of the MSE under local differential privacy and not the dependence structure.
The code for the simulation study is available on \href{https://github.com/kklockmann/Nonparametric-Spectral-Density-Estimation-LDP}{GitHub}.

\begin{figure}[tbp]
% \ContinuedFloat
\begin{subfigure}{0.48\textwidth} 
	\includegraphics[width=0.85\textwidth]{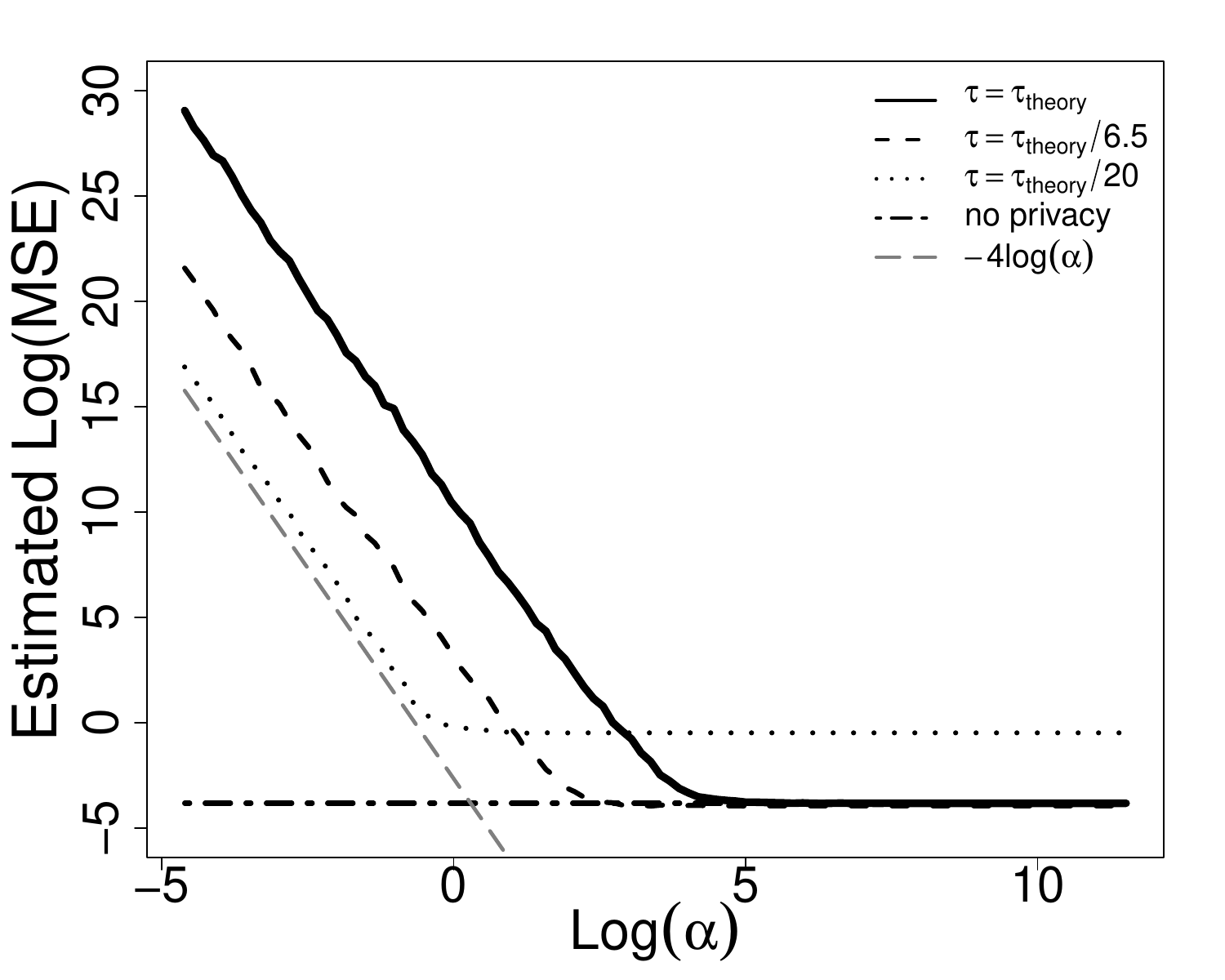}
	\caption{Estimation of $\sigma_0{=}1.44$ from NI-privatized data.}
	\label{EX1_loglog_plot_MSE_varyalpha_parta}
\end{subfigure}
\begin{subfigure}{0.48\textwidth}
	\includegraphics[width=0.85\textwidth]{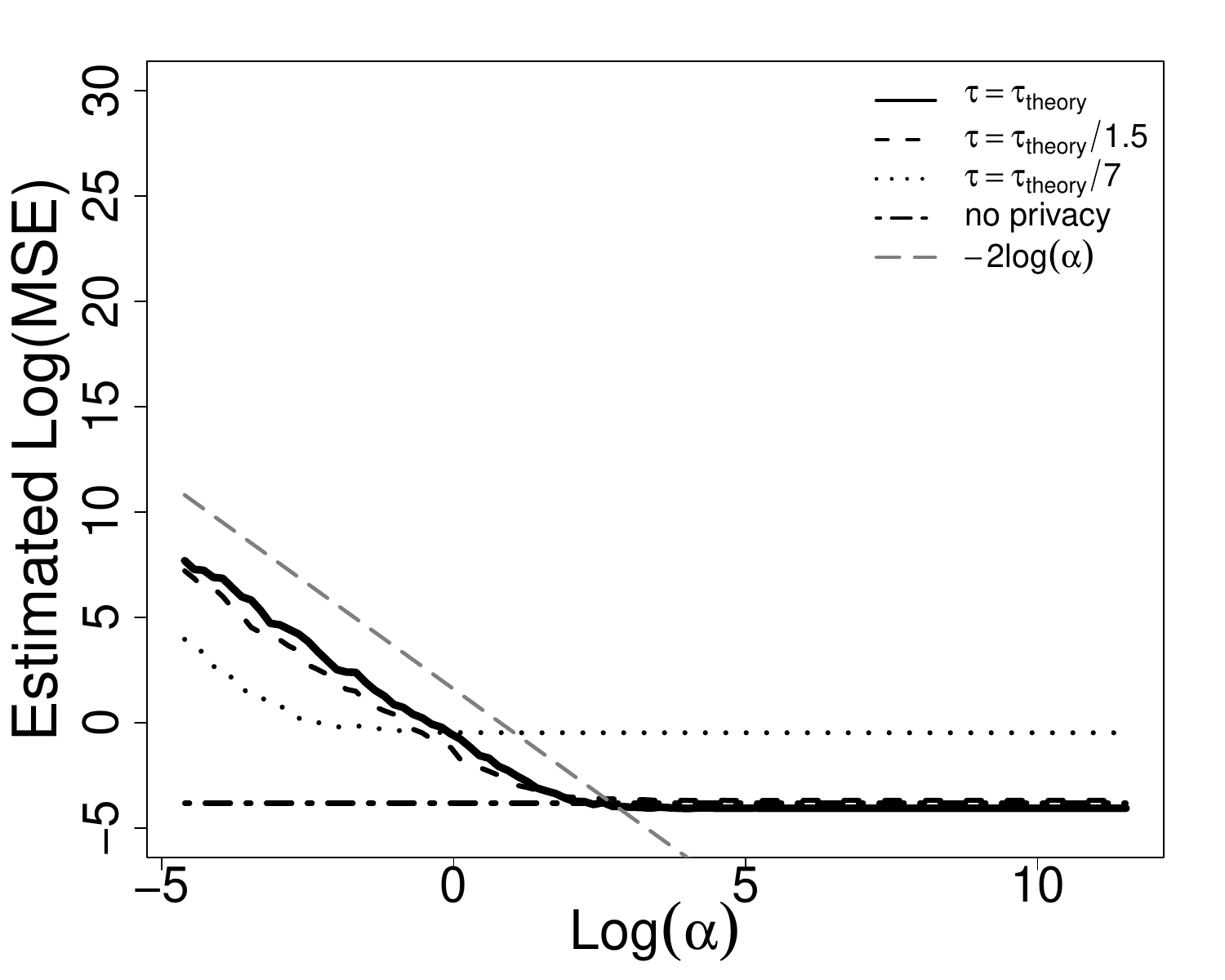}
	\caption{Estimation of $\sigma_0{=}1.44$ from SI-privatized data.}
\end{subfigure}
\begin{subfigure}{0.48\textwidth}
	\includegraphics[width=0.85\textwidth]{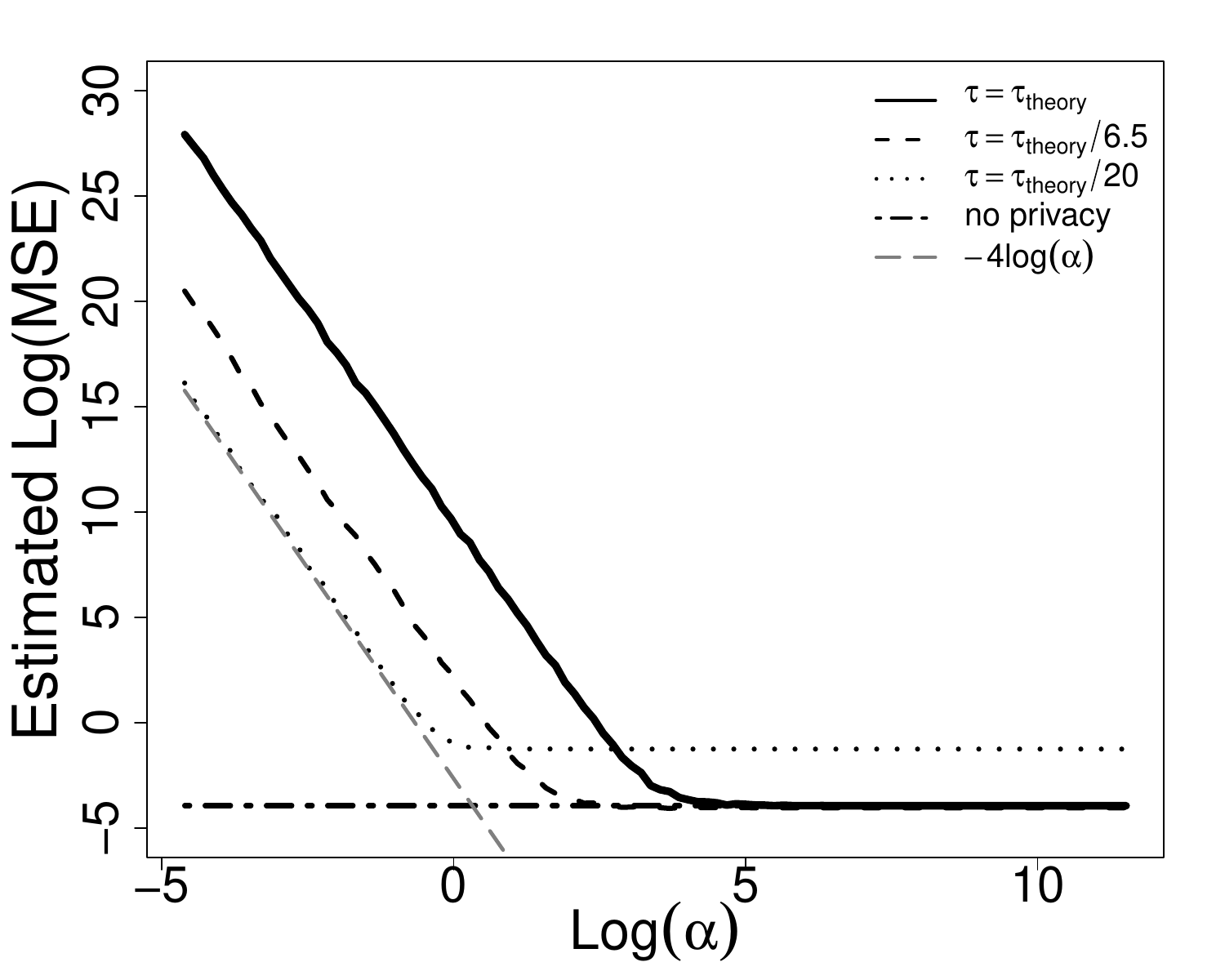} 
	\caption{Estimation of $\sigma_2{=}0.92$ from NI-privatized data.}
\end{subfigure}
\begin{subfigure}{0.48\textwidth}
	\includegraphics[width=0.85\textwidth]{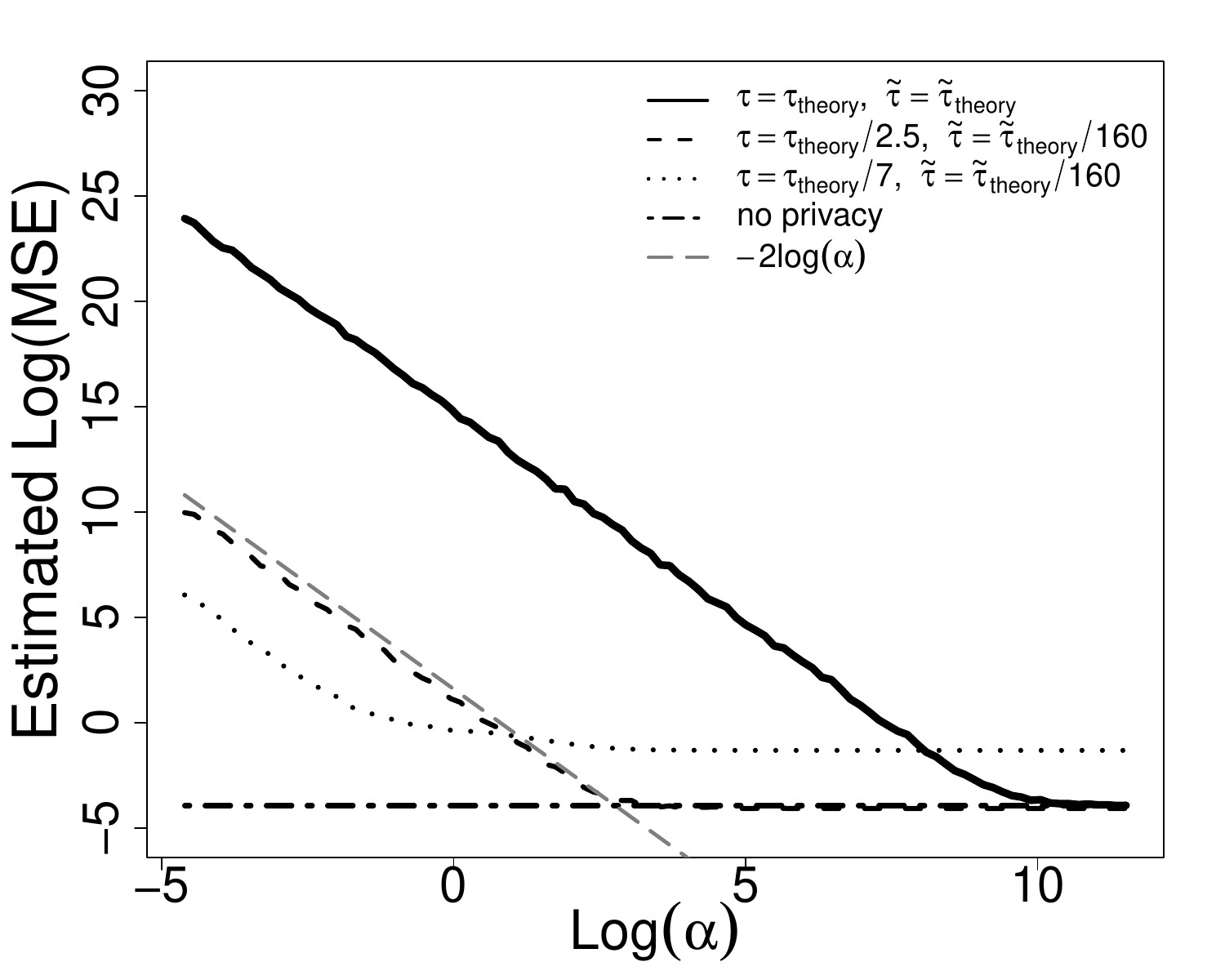} 
	\caption{Estimation of $\sigma_2{=}0.92$ from SI-privatized data.}
\end{subfigure}	
\begin{subfigure}{0.48\textwidth}
	\includegraphics[width=0.85\textwidth]{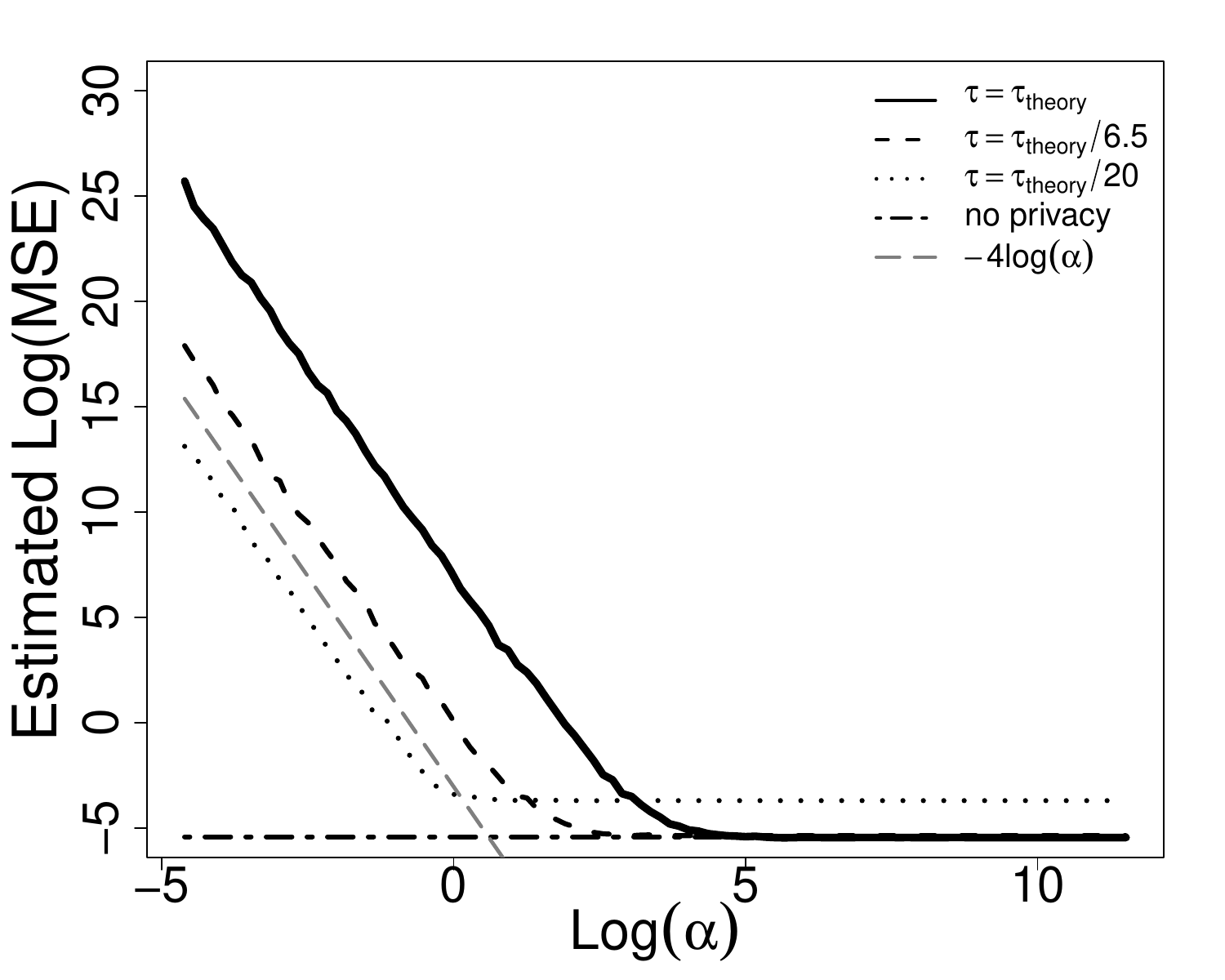} 
	\caption{Estimation of $f(\pi/5){=} 0.24$ from NI-privatized data.}
\end{subfigure} \hspace*{0.31cm}
\begin{subfigure}{0.48\textwidth}
	\includegraphics[width=0.85\textwidth]{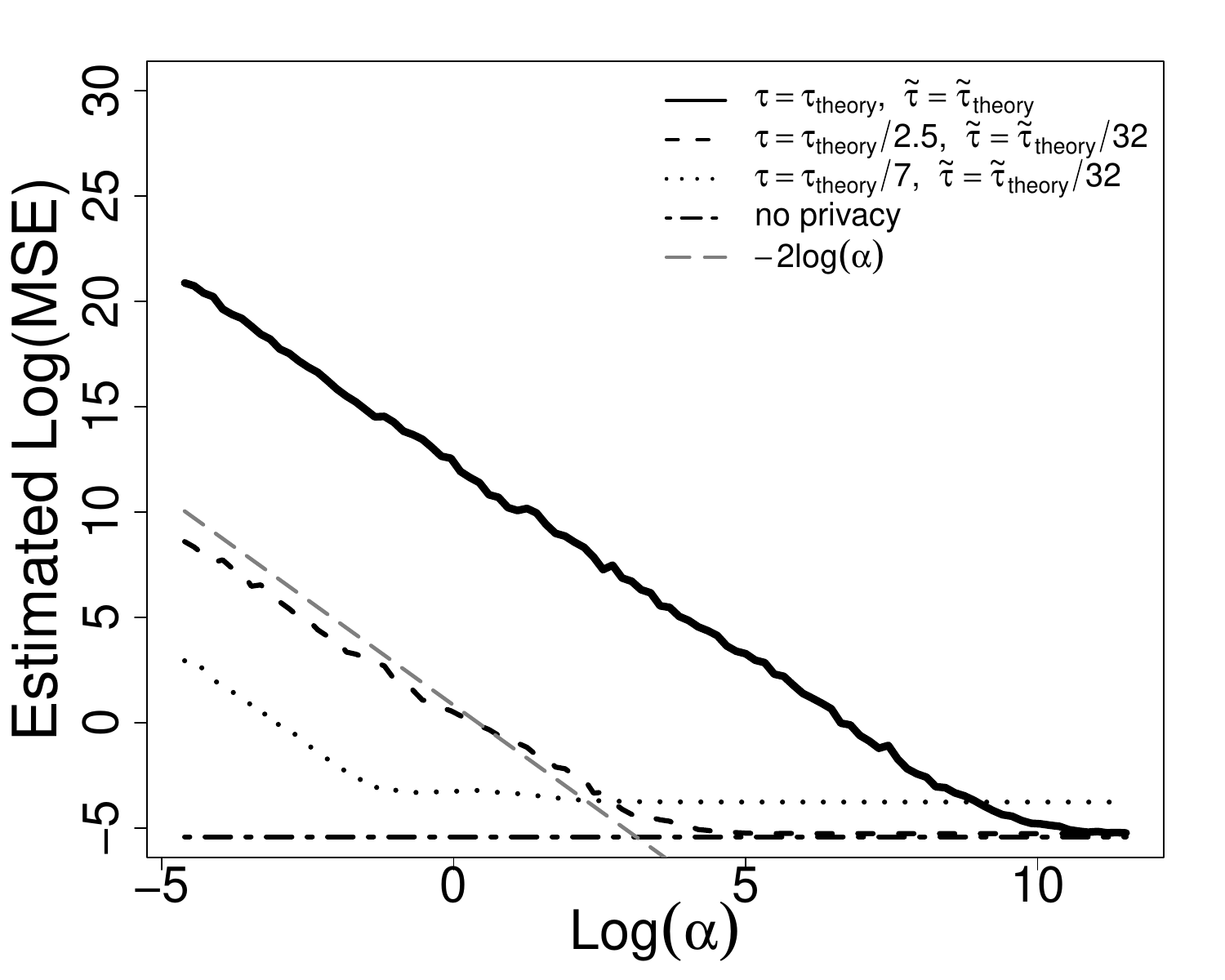} 
	\caption{Estimation of $f(\pi/5){=} 0.24$  from SI-privatized data.}
\end{subfigure}
\caption{
	Example~(1): Log-log plot of the estimated mean squared error (MSE) versus the privacy level $\alpha$ for estimation of  variance  (top), covariance coefficient (center) and  pointwise spectral density (bottom) using NI‑privatized data (left) and  SI‑privatized data (right).
	Gray dashed lines indicate the theoretical $\alpha$‑scaling of the MSE for NI‑ and SI‑based estimators, including $\log(n)$-factors.
}
\label{EX1_loglog_plot_MSE_varyalpha}
\end{figure}

%%%%%%%%%%%%%%%%%%%%%%%%%%%%%%%%%%%%%%%%%%%%%%%%%%%%%%%%%%%%%%%%%%
%APPENDIX: Proofs
%%%%%%%%%%%%%%%%%%%%%%%%%%%%%%%%%%%%%%%%%%%%%%%%%%%%%%%%%%%%%%%%%%
\appendix
\noindent
\section{Proofs}
Throughout this section, we use $c, c_1, C, C_1, \dots $ to denote generic constants that are independent of  $K, n,\alpha$. The indicator function is written as  $\mathbf{1}(\text{condition})$, which equals $1$ if the condition holds  and 0 otherwise. By $\mathbb{I}_A$, we denote a random variable that takes the value 1 if the event $A$ occurs and 0 otherwise. %To simplify the notation, the constants are sometimes skipped, and we write $\lesssim$ for less than or equal to up to constants.
\subsection{Proof of Proposition~\ref{prop:ratekroll} }
Let 
$\hat{I}_n(\omega)=\hat{I}_n^Z(\omega)- \frac{4\tau_n^2}{\pi \alpha^2}  $
with $\hat{I}_n^Z(\omega)= \frac{1}{2\pi n}\left | \sum_{t=1}^n Z_t\exp(-\texttt{i}t\omega) \right|^2$.
Then,  %we can rewrite $\hat{I}_n^Z(\omega)$ as 
\begin{align*}
\hat{I}_n^Z(\omega){=}\frac{1}{2\pi n}\left | \sum_{t=1}^n Z_t\exp(-\texttt{i}t\omega) \right|^2{=}\frac{1}{2\pi n} \sum_{s=1}^n \sum_{t=1}^n Z_s Z_t\exp(\texttt{i}\omega(s{-}t)){=}\frac{1}{2\pi} \sum_{|h|<n} \hat{\gamma}(h)\exp({-}\texttt{i}h\omega),
\end{align*}
where $\hat{\gamma}(h)= \frac{1}{n} \sum_{t=1}^{n-|h|} Z_tZ_{t+|h|}$ is the sample covariance estimator.
Thus, for $j{\in}\{{-}R,...,R\}$,
\begin{align}
\hat{\sigma}^{\mathrm{NI}}_j &=\int_{-\pi}^\pi \hat{I}_n(\omega) \exp(\texttt{i}j \omega)\, \mathrm{d}\omega=\int_{-\pi}^\pi \hat{I}^Z_n(\omega) \exp(\texttt{i}j \omega)\, \mathrm{d}\omega - \frac{8\tau_n^2}{ \alpha^2} \delta_{j,0}\nonumber\\
&=\int_{-\pi}^\pi \frac{1}{2\pi} \sum_{|h|<n} \hat{\gamma}(h)\exp(-\texttt{i}h\omega)  \exp(\texttt{i}j \omega)\, \mathrm{d}\omega - \frac{8\tau_n^2}{ \alpha^2} \delta_{j,0}\nonumber\\
&= \frac{1}{2\pi} \sum_{|h|<n} \hat{\gamma}(h)\int_{-\pi}^\pi \exp\{\texttt{i}\omega(j-h)\}\, \mathrm{d}\omega - \frac{8\tau_n^2}{ \alpha^2} \delta_{j,0}\nonumber\\
&=\sum_{|h|<n} \hat{\gamma}(h)\delta_{h,j} - \frac{8\tau_n^2}{ \alpha^2} \delta_{j,0}
%&=  \hat{\gamma}(j)- \frac{4\tau_n^2}{\pi \alpha^2} \delta_{j,0}
=\begin{cases}
	\frac{1}{ n}  \sum_{t=1}^{n} Z_t^2 - \frac{8\tau_n^2}{ \alpha^2}, & j=0,\\
	\frac{1}{ n}   \sum_{t=1}^{n-|j|} Z_tZ_{t+|j|}, & j\neq 0.
\end{cases} \label{sampleacf}
\end{align}
This implies
\begin{align*}
\hat{f}^{\mathrm{NI}}_m(\omega)=\frac{1}{2\pi}\sum_{j=-m}^m \hat{\sigma}^{\mathrm{NI}}_j \exp(-\texttt{i}j \omega)=\frac{1}{2\pi}\sum_{j=-m}^m \hat{\gamma}(j) \exp(-\texttt{i}j \omega)  - \frac{4\tau_n^2}{\pi \alpha^2}.
\end{align*}

\emph{Proof of $i)$}:
As $\hat{\sigma}^{\text{NI}}_j=\hat{\sigma}^{\text{NI}}_{-j}$, it is sufficient to consider $j=0,1,...,R$. We start with $j\neq 0$. Consider the event $A_{\text{NI}}$ and its complement defined by
\begin{align*}
A_{\text{NI}} &{=} \bigcap_{i{=}1}^n \{ \widetilde{X}_i{=} X_i  \}, 
\quad A_{\text{NI}}^c {=}\bigcup_{i{=}1}^n \{ \widetilde{X}_i{\neq} X_i \}. 
\end{align*}
\emph{$L_2$ risk on the event $A_{\text{NI}}$.} 
If the event $A_{\text{NI}}$ occurs, then $\hat{\sigma}^{\text{NI}}_j$ coincides with the untruncated estimator $\tilde{\sigma}^{\text{NI}}_j$, which is defined by 
\begin{align} %\label{stildeNI}
\tilde{\sigma}^{\text{NI}}_j & =\frac{1}{ n} \sum_{i=1}^{n-j}  (X_i +\xi_i) \cdot \left ( X_{i+j}+\xi_{i+j} \right). \label{hatsimgaj}
\end{align}
In particular,
\begin{align*}
\mathbb{E}[(\hat{\sigma}^{\text{NI}}_j - \sigma_j)^2] &= \mathbb{E}[\mathbb{I}_{A_{\text{NI}}}(\hat{\sigma}^{\text{NI}}_j - \sigma_j)^2] + \mathbb{E}[\mathbb{I}_{A_{\text{NI}}^c}(\hat{\sigma}^{\text{NI}}_j - \sigma_j)^2] %\label{decompose}
\end{align*}
\begin{align*}
\text{and } \quad \mathbb{E}[\mathbb{I}_{A_{\text{NI}}}(\hat{\sigma}^{\text{NI}}_j - \sigma_j)^2] &\leq \mathbb{E}[(\tilde{\sigma}^{\text{NI}}_j - \sigma_j)^2] %=\mathbb{E}[\{\tilde{\sigma}_j - \mathbb{E}[\tilde{\sigma}_j] + \{\mathbb{E}[\tilde{\sigma}_j] - {\sigma}_j\}^2] \\
= \mathbb{E}[(\tilde{\sigma}^{\text{NI}}_j - \mathbb{E}[\tilde{\sigma}^{\text{NI}}_j ])^2]  + (\mathbb{E}[\tilde{\sigma}^{\text{NI}}_j ] - {\sigma}_j)^2.
\end{align*}
It is straightforward to see that 
\begin{align*}
\mathbb{E}[ \tilde{\sigma}^{\text{NI}}_j]& 
=  \frac{1}{ n} \sum_{i=1}^{n-j} \mathbb{E} \left[ X_i X_{i+j}\right ] = \frac{1}{ n}\sum_{i=1}^{n-j}  \sigma_j = \frac{n-j}{n}\sigma_j.
\end{align*} As $0\leq j\leq R<n$ and $R=\mathcal{O}(n^{1/2})$, it follows $ (\mathbb{E}[\tilde{\sigma}^{\text{NI}}_j ] - {\sigma}_j)^2=\frac{j^2}{n^2}\sigma_j^2=\mathcal{O}(n^{-1})$. The hidden constant does not depend on $j$ since $\sigma_j^2\leq M_1$ by assumption.
Let $\bar{\sigma}_j = \frac{1}{ n} \sum_{i=1}^{n-j}   X_{i}X_{i+j}$ be the biased estimator built from the original data. Then,
\begin{align*}
\tilde{\sigma}^{\text{NI}}_j&= \frac{1}{ n} \sum_{i=1}^{n-j}  (X_i +\xi_i) \cdot \left ( X_{i+j}+\xi_{i+j}\right) \\
&= \bar{\sigma}_j + \frac{1}{ n} \sum_{i=1}^{n-j}  X_i\xi_{i+j} + \frac{1}{ n}  \sum_{i=1}^{n-j}   X_{i+j}\xi_{i} + \frac{1}{ n}  \sum_{i=1}^{n-j}   \xi_{i+j}\xi_{i}\\
& =:  \bar{\sigma}_j+B+C+D. \nonumber
\end{align*} 
Since all terms of $B,C$ and  $D$ are centered, uncorrelated and uncorrelated to the terms in $\bar \sigma_j$, we have
\begin{equation}
\mathbb{E}\left[(\tilde{\sigma}^{\text{NI}}_j - \mathbb{E}[\tilde{\sigma}^{\text{NI}}_j ])^2 \right ]=\mathrm{Var}[\bar{\sigma}_j]+\mathbb{E}[B^2]+\mathbb{E}[C^2]+\mathbb{E}[D^2]. \label{sigmaABCkroll}
\end{equation}
For the first term on the right-hand side of the above equality, we obtain
\begin{align*}
\mathrm{Var}[\bar{\sigma}_j]&= \frac{1}{n^2}  \sum_{i=1}^{n-j} \sum_{k=1}^{n-j}   \mbox{Cov} \left (   X_{i}X_{i+j} , X_{k}X_{k+j} \right)=\frac{1}{n^2}  \sum_{i=1}^{n-j}  \sum_{k=1}^{n-j}  \mathbb{E}\left [X_{i}X_{i+j} X_{k}X_{k+j} \right  ] -   \sigma_j^2.
\end{align*}
By Wick's formula, we have \begin{align}
\mathbb{E}\left [X_{i}X_{i+j} X_{k}X_{k+j} \right  ]&=\mathbb{E}\left [X_{i}X_{i+j}] \mathbb{E}[ X_{k}X_{k+j} \right  ] + \mathbb{E}\left [X_{i}  X_{k}] \mathbb{E}[X_{i+j}X_{k+j} \right  ] + \mathbb{E}\left [X_{i} X_{k+j} ] \mathbb{E}[ X_{k} X_{i+j}\right  ] \nonumber\\
&=\sigma_j^2 + \sigma_{|i-k|}^2 + \sigma_{|i-k-j|} \sigma_{|i+j-k|}. \label{Wick}
\end{align} 
Using the Cauchy-Schwarz inequality and $\sum_{j\in\mathbb{Z}}\sigma_j^2\leq M_1$, we get
\begin{align*}
\mathrm{Var}[\bar{\sigma}_j]&{=}\frac{1}{n^2} \sum_{i=1}^{n-j}\sum_{k=1}^{n-j} \sigma_{|i-k|}^2 + \sigma_{|i-k+j|} \sigma_{|i-k-j|} {=}\frac{1}{n^2} \sum_{h=-(n-j-1)}^{n-j-1}  (n-j-|h|)(\sigma_{|h|}^2 + \sigma_{|h+j|} \sigma_{|h-j|})\\
%&=\frac{1}{n} \sum_{h=-(n-K-1)}^{n-K-1}  \left (1-\frac{|h|}{n-K}\right)(\sigma_{|h|}^2 + \sigma_{|h+j|} \sigma_{|h-j|})\\
&\leq\frac{1}{n} \left ( \sum_{h=-n}^{n}  \sigma_{|h|}^2  +  |\sigma_{|h+j|}|| \sigma_{|h-j|}|\right )\leq\frac{2M_1}{n}.
\end{align*} %where we used the Cauchy-Schwarz inequality and that $\sum_{l\in\mathbb{Z}}\sigma_{l}^2\leq M_1$.
Next, since the $\xi_i$ are independent of the $X_i$ and centered,
\begin{align*}
\mathbb{E}[B^2]&=   \frac{1}{ n^2} \sum_{i=1}^{n-j}  \sum_{k=1}^{n-j}  \mathbb{E}[X_iX_k]\mathbb{E}[\xi_{i+j}\xi_{k+j}] =  \frac{8\sigma_0 \tau_n^2(n-j)}{\alpha^2n^2}\leq c\frac{8\sigma_0 \tau_n^2}{\alpha^2n}.
\end{align*}
Similar computations give the same bound for $\mathbb{E}[C^2]$.
For the next term in (\ref{sigmaABCkroll}) we have
$$
\mathbb{E}[D^2]=\frac{1}{n^2}\sum_{i=1}^{n-j}  \mathbb{E}[\xi_{i+j}^2\xi_{i}^2]=\frac{1}{n^2}\sum_{i=1}^{n-j}  \mathbb{E}[\xi_{i+j}^2]\mathbb{E}[\xi_{i}^2]= \frac{64  \tau_n^4(n-j)}{\alpha^4n^2} \leq c \frac{64  \tau_n^4}{\alpha^4n}.
$$
Inserting the previous bounds into (\ref{sigmaABCkroll}), we get, for $j\in\{1,2,...,R\}$, 
\begin{align*}\mathbb{E}[(\hat{\sigma}^{\text{NI}}_j - \sigma_j)^2] &\lesssim \frac{2M_1}{n} + \frac{8\sigma_0 \tau_n^2}{\alpha^2n}+\frac{64  \tau_n^4}{\alpha^4n} + \mathbb{E}[\mathbb{I}_{A_{\text{NI}}^c}(\hat{\sigma}_j - \sigma_j)^2].%\\
%& \leq C\max \left \{\frac{(\tau_n^2,\tilde \tau_n^2)}{\alpha^2(n-K)}  + \mathbb{E}[\mathbb{I}_{A^c}(\hat{\sigma}_j - \sigma_j)^2].
\end{align*}
If $j=0$, then define $\tilde \sigma_0^{\text{NI}}$ by \begin{equation} \label{tildesigma0}
\tilde{\sigma}^{\text{NI}}_0=\frac{1}{ n} \sum_{i=1}^{n}  (X_i +\xi_i)^2 - \frac{8\tau_n^2}{\alpha^2}.
\end{equation}
As $\mathbb{E}[ \tilde{\sigma}^{\text{NI}}_0]
=  \frac{1}{ n} \sum_{i=1}^{n} \mathbb{E} \left[ X_i^2\right ] +\mathbb E[\xi_i^2] - \frac{8\tau_n^2}{\alpha^2} =\sigma_0$, it follows
\begin{align*}
\mathbb{E}[\mathbb{I}_{A_{\text{NI}}}(\hat{\sigma}^{\text{NI}}_0 - \sigma_0)^2] \leq  \mathbb{E}[(\tilde{\sigma}^{\text{NI}}_0 - \sigma_0)^2] = \mathrm{Var}(\bar{\sigma}_0 )+2\mathbb{E}[C^2]+\mathbb{E}[D^2]
\end{align*} with  $C=\frac{2}{n}\sum_{i=1}^n X_i\xi_i$ and $D=\frac{1}{n}\sum_{i=1}^n \xi_i^2$. Similar calculations as before yield
\begin{align*}
\mathbb{E}[(\hat{\sigma}^{\text{NI}}_0 - \sigma_0 )^2] &\lesssim \frac{2M_1}{n} + \frac{8\sigma_0 \tau_n^2}{\alpha^2n}+\frac{96 \tau_n^4}{\alpha^4n} + \mathbb{E}[\mathbb{I}_{A_{\text{NI}}^c}(\hat{\sigma}^{\text{NI}}_0 - \sigma_0)^2].
\end{align*} 
\emph{$L_2$ risk on the event $A_{\text{NI}}^c$.}
By Parseval's equality and the results of \cite{kroll2024nonparametric} (see proof of Theorem 3.2 on page 19), it follows with $\tau_n^2=56 \log^{1+\delta}(n)$ for some $c,\delta>0$, that
\begin{align*}
\mathbb E[|\hat{\sigma}^{\text{NI}}_j-\sigma_j|^2 \mathbb{I}_{A_{\text{NI}}^c}]&\leq\mathbb E[\|\hat{f}^{\text{NI}}_m-f_m\|^2_2\mathbb{I}_{A_{\text{NI}}^c}]\leq c \left ( \frac{\tau_n^4}{n} \lor \frac{\tau_n^4}{n\alpha^4} \right), \quad j\in\{0,1,...,R\},
\end{align*} where $f_m(\cdot)=\frac{1}{2\pi}\sum_{k=-m}^m \sigma_k \exp(-\texttt{i}k\cdot)$  and $R\leq m\leq n$.
Then, for sufficiently large $n$,
$$ \mathbb E|\hat{\sigma}^{\text{NI}}_j-\sigma_j|^2 \leq c \left [  \frac{\log^{2+2\delta}(n)}{n} \lor \frac{\log^{2+2\delta}(n)}{n\alpha^4} \right].$$ %\frac{1}{n} smaller than log^{2+2\delta}(n)}{n} for n>2

\emph{Proof of $ii)$ and $iii)$}:
Let $\omega\in[-\pi,\pi]$ be arbitrary but fixed. Consider the event $A_{\text{NI}}$ and its complement defined as before.\\
\emph{$L_2$ risk on the event $A_{\text{NI}}$.} 
If the event $A_{\text{NI}}$ occurs, then $\hat{f}^{\text{NI}}_m(\omega)$ coincides with the untruncated estimator $\tilde{f}^{\text{NI}}_m(\omega)$, which is defined by 
\begin{align*}
\tilde{f}^{\text{NI}}_m(\omega)&=\frac{1}{2 \pi} \sum_{j=-m}^{m}\tilde{\sigma}^{\text{NI}}_j \exp(-\texttt{i}j\omega)\\
&= \frac{1}{2 \pi n} \sum_{j=-m}^{m} \sum_{i=1}^{n-|j|}  \left[(X_i +\xi_i) \cdot  ( X_{i+|j|}+\xi_{i+|j|})  \exp(-\texttt{i} j\omega)-\delta_{j,0}   \frac{8\tau_n^2}{\alpha^2}\right ],
\end{align*} where we used identity (\ref{sampleacf}), (\ref{hatsimgaj}) and (\ref{tildesigma0}).
In particular,
\begin{align*}\quad \quad \,\mathbb{E}[\{\hat{f}^{\text{NI}}_m(\omega) - f(\omega)\}^2]&= \mathbb{E}[\mathbb{I}_{A_{\text{NI}}}\{\hat{f}^{\text{NI}}_m(\omega) - f(\omega)\}^2] + \mathbb{E}[\mathbb{I}_{A_{\text{NI}}^c}\{\hat{f}^{\text{NI}}_m(\omega) - f(\omega)\}^2] %\label{f0decomposekroll}
\end{align*} and
\begin{align*}
\mathbb{E}[\mathbb{I}_{A_{\text{NI}}}\{\hat{f}^{\text{NI}}_m(\omega) - f(\omega)\}^2]&\leq \mathbb{E}[\{\tilde{f}^{\text{NI}}_m(\omega) - f(\omega)\}^2] \\ %= \mathbb{E}[\{\tilde{f}_m(\omega) - \mathbb{E}[\tilde{f}_m(\omega)] + \{\mathbb{E}[\tilde{f}_m(\omega)] - f(\omega)\}^2] \\
&= \mathbb{E}[\{\tilde{f}^{\text{NI}}_m(\omega) - \mathbb{E}[\tilde{f}^{\text{NI}}_m(\omega)]\}^2]  + \{\mathbb{E}[\tilde{f}^{\text{NI}}_m(\omega)] - f(\omega)\}^2.
\end{align*}
%\subsubsection{The Bias $ \{\mathbb{E}[\tilde{f}_m(\omega)] - f(\omega)\}^2$ }
We start with calculating the expectation of $\tilde{f}^{\text{NI}}_m(\omega)$.
\begin{align*}
\mathbb{E}[ \tilde{f}^{\text{NI}}_m(\omega)]& %\frac{1}{\pi(n-k)} \sum_{j=K+1}^n \mathbb{E} \left [ \sum_{|k|\leq K} X_j Z_{j-k} \exp(\texttt{i}\omega k) + \tilde{\xi}_j \right ] 
=  \frac{1}{2\pi n} \sum_{j=-m}^{m} \sum_{i=1}^{n-|j|} \left [ \left (\mathbb{E} [ X_i X_{i+|j|} ]+\mathbb{E} [ \xi_i \xi_{i+|j|} ]\right) \exp(-\texttt{i}j\omega)-\delta_{j,0}   \frac{8\tau_n^2}{\alpha^2}\right ] \\
%=\mathbb{E}[ \bar{f}_K(\omega)]
&=   \sum_{j=-m}^{m} \frac{n-|j|}{2\pi n}\sigma_j\exp(-\texttt{i}j\omega ) = \frac{1}{2\pi} \sum_{|j|\leq m}  a_j\sigma_j \exp(-\texttt{i}j\omega )
\end{align*} with $a_j=1-|j|/n$ for $j=-m,...,m$. Then,  using (\ref{constants}),
\begin{align*}
&\left|\frac{1}{2\pi} \sum_{|j|\leq m}  a_j\sigma_j \exp(-\texttt{i}j\omega )-f(\omega)\right|\\
&= \frac{1}{2\pi} \left| \sum_{|j|\leq m}  -\frac{|j|}{n}\sigma_j \exp(-\texttt{i}j\omega )-\sum_{|j|> m}  \sigma_j\exp(-\texttt{i}j\omega )\right|\\
&\leq \frac{Mm}{2\pi n}  + \frac{1}{2\pi}\left |\sum_{|j|> m}  \sigma_j\exp(-\texttt{i}j\omega )\right|.
\end{align*} 
If $f\in W^{s,2}(L),\, s>1/2,$ then $\sum_{j\in \mathbb{Z}}  |\sigma_j|^2|j|^{2s}<L$. Thus, by the Cauchy-Schwarz inequality we obtain % $|\sigma_k|\leq o(k^{-s})$ holds uniformly in $ W^{s,2}(L)$. Thus,
\begin{align} 
\left |\sum_{|j|> m}  \sigma_j\exp(-\texttt{i}j\omega )\right|&\leq  c\sum_{|j|> m}  |\sigma_j|\frac{|j|^s}{|j|^s} \leq  c\sum_{|j|> m}  |\sigma_j|^2|j|^{2s}\sum_{|j|> m}  \frac{1}{|j|^{2s}}= \mathcal{O}(m^{-s+1/2}). \label{FSapprox}
\end{align} % where we used (\ref{constants}) and that for $f\in W^{s,2}(L)$.
If $f\in W^{s,\infty}(L_0,L)$,
then by the approximation results of \citet[page 120]{zygmund2002trigonometric}
\begin{align}
\left|\frac{1}{2\pi} \sum_{|j|> m}  \sigma_j \exp(-\texttt{i}j\omega )\right|= \mathcal{O}\left ( \frac{\log(m)}{m^s} \right).  \label{FSapprox2}
\end{align}
Thus, we obtain for the squared bias 
\begin{equation*}
\{\mathbb{E}[ \tilde{f}^{\text{NI}}_m(\omega)]-f(\omega)\}^2\leq C\frac{m^2}{n^2}+ C \cdot
\begin{cases} 
\log(m)^2m^{-2s}, & f\in W^{s,\infty}(L_0,L),\\
m^{-2s+1}, &f\in W^{s,2}(L),
\end{cases}
\end{equation*} where the constant $C>0$ depends on $s,\,L$, if $f\in W^{s,2}(L)$, and additionally on $L_0$, if $ f\in W^{s,\infty}(L_0,L).$
%REMARK: our estimator constructs VP-sum, but here only Fejer sum, thus an additional log-factor,
%recall l_sup norm +Hölder space, then log-factor cannot be avoided,
%but here pointwise rate+Hölder space, optimal rate (non-private rate) has no log-factor.

To analyze the variance, define the estimator without truncation or added Laplace noise by
$$
\bar{f}_m(\omega) = \frac{1}{2 \pi n} \sum_{j=-m}^{m} \sum_{i=1}^{n-|j|} X_{i}X_{i+|j|}\exp(-\texttt{i}j\omega ).
$$
Then,
\begin{align*}
\tilde{f}^{\text{NI}}_m(\omega)=& \frac{1}{2 \pi n} \sum_{j=-m}^{m} \sum_{i=1}^{n-|j|}  \left[(X_i +\xi_i) \cdot  ( X_{i+|j|}+\xi_{i+|j|})  \exp(-\texttt{i} j\omega)-\delta_{j,0}   \frac{8\tau_n^2}{\alpha^2}\right ]\\
=& \, \frac{1}{2 \pi n} \sum_{j=-m}^{m} \sum_{i=1}^{n-|j|} X_iX_{i+|j|} \exp(-\texttt{i} j\omega) +  \frac{1}{2 \pi n} \sum_{j=-m}^{m} \sum_{i=1}^{n-|j|} X_i\xi_{i+|j|}  \exp(-\texttt{i} j\omega)\\
&\,+  \frac{1}{2 \pi n} \sum_{j=-m}^{m} \sum_{i=1}^{n-|j|} \xi_iX_{i+|j|}  \exp(-\texttt{i} j\omega) +  \frac{1}{2 \pi n} \sum_{j=-m}^{m} \sum_{i=1}^{n-|j|} \xi_i\xi_{i+|j|}  \exp(-\texttt{i} j\omega) -\delta_{j,0}   \frac{8\tau_n^2}{\alpha^2}\\
=:& \,\bar{f}_m(\omega)+B+C+D.
\end{align*} 
In particular, $\mathbb{E}[\tilde{f}^{\text{NI}}_m(\omega)]=\mathbb{E}[\bar{f}_m(\omega)]$ and 
\begin{equation}
\mathbb{E}[\{\tilde{f}^{\text{NI}}_m(\omega) - \mathbb{E}[\tilde{f}^{\text{NI}}_m(\omega)]\}^2] =\mathrm{Var}[\bar{f}_m(\omega)]+\mathbb{E}[B^2]+\mathbb{E}[C^2]+\mathbb{E}[D^2]. \label{f0ABCkroll}
\end{equation}
Using $ \bar{\sigma}_j =\frac{1}{ n} \sum_{i=1}^{n-|j|}  X_i X_{i+|j|}$%and the second-order stationarity of $(X_1,...,X_n)$
, we have
\begin{align}
&(2\pi n)^2   \mathrm{Var}[\bar{f}_m(\omega)] \nonumber\\
&=  n^2\sum_{l=-m}^m\sum_{k=-m}^m \exp\{-\texttt{i}\omega (l-k)\}\text{Cov}(\bar{\sigma}_l, \bar{\sigma}_k)\nonumber\\
&= \sum_{l=-m}^m\sum_{k=-m}^m \exp\{-\texttt{i}\omega (l-k)\}\sum_{r=1}^{n-|l|} \sum_{s=1}^{n-|k|} \text{Cov}(X_rX_{r+|l|},X_{s}X_{s+|k|}). \label{Varf0calckroll}
\end{align}
Now,
\begin{align*}
&\text{Cov}(X_rX_{r+|l|},X_{s}X_{s+|k|})\\
&=\mathbb E[X_rX_{r+|l|}]\mathbb E[X_{s}X_{s+|k|}]+ \mathbb E[X_rX_{s}] \mathbb E[X_{s+|k|}X_{r+|l|}]+\mathbb E[X_rX_{s+|k|}] \mathbb E[X_{s}X_{r+|l|}]-\sigma_{|l|}\sigma_{|k|}\\
&=\sigma_{|r-s|}\sigma_{|r-s-|k|+|l||}+\sigma_{|r-s-|k||} \sigma_{|r-s+|l||}.
\end{align*}
Substituting this into \eqref{Varf0calckroll}, we get 
\begin{align*}
&(2\pi n)^2    \mathrm{Var}[\bar{f}_m(\omega)]\\
&\leq \sum_{l=-m}^m\sum_{k=-m}^m \sum_{h=-(n-1)}^{n-1}  (n-1-|h|)\left|\sigma_{|h|}\sigma_{|h-|k|+|l||}+\sigma_{|h-|k||} \sigma_{|h+|l||}\right|\\
&\leq \sum_{l=-m}^m\sum_{h=-(n-1)}^{n-1}  (n-1-|h|) \left[\left |\sigma_{|h|}\right|\sum_{k=-m}^m\left |\sigma_{|h-|k|+|l||}\right|+ \left|\sigma_{|h+|l||}\right| \sum_{k=-m}^m\left| \sigma_{|h-|k||}\right| \right ]\\
&\leq \sum_{l=-m}^m\sum_{h=-(n-1)}^{n-1}  (n-1-|h|)(|\sigma_{|h|}|+| \sigma_{|h+|l||}|)M. 
% \leq 2M^2(2K+1).
\end{align*}
Since $\sum_{h=-(n-1)}^{n-1}  \left (1-\frac{1+|h|}{n} \right)(|\sigma_{|h|}|+| \sigma_{|h+|l||}|) \leq 2M$, it follows $\mathrm{Var}[\bar{f}_m(\omega)]\leq Cmn^{-1}.$
Similar computations as in the proof of Theorem~\ref{theo:ratef0} show that $\mathbb{E}[B^2]$ and $\mathbb E[C^2]$ can be bounded by $\frac{\tau_n^2m}{n\alpha^2}$ up to constants.
For the last term in Equation~(\ref{f0ABCkroll}) we have 
\begin{align*}
\mathbb{E}[D^2]=\mathrm{Var}[D]&= \frac{1}{4 \pi^2 n^2} \sum_{j=-m}^{m} \sum_{i=1}^{n-|j|} \mathrm{Var}[\xi_i\xi_{i+|j|}]  \leq c\frac{m\tau_n^4}{\alpha^4n}.
\end{align*}
Summing up, we get
\begin{align}
\mathbb{E}[\{\tilde{f}^{\text{NI}}_m(\omega) - f(\omega)]\}^2]=\, &  \mathbb{E}[\mathbb{I}_{A_{\text{NI}}^c}\{\hat{f}^{\text{NI}}_m(\omega) - f(\omega)\}^2]+\mathcal{O}\left ( \frac{m}{n} \lor\frac{m\tau_n^{4} }{n\alpha^ {4}}\right)\nonumber\\
&  + \begin{cases}
\mathcal{O}(\log(m)^2m^{-2s}), &f\in W^{s,\infty}(L_0,L),\\
\mathcal{O}(m^{-2s+1}), &f\in  W^{s,2}(L).
\end{cases}  \label{boundeventAkroll}
\end{align}
\emph{$L_2$ risk on the event $A_{\text{NI}}$.}
Application of the Cauchy-Schwarz inequality gives
\begin{align*} %\label{f0CSbound}
\mathbb{E}[\mathbbm{I}_{A_{\text{NI}}^c}\{\hat{f}^{\text{NI}}_m(\omega) - f(\omega)\}^2] \leq \left(\mathbb E[\{\hat{f}^{\text{NI}}_m(\omega) - f(\omega)\}^4] \cdot \Pr[ A_{\text{NI}}^c ]\right)^{1/2} .
\end{align*}
%Let $f_m(\omega)=\sum_{|l|\leq m}\sigma_l\exp(-\texttt{i}l\omega)$. defined before
From the proof of Theorem 3.2. of \cite{kroll2024nonparametric} follows
\begin{align*}
\mathbb{E}[\{\hat{f}^{\text{NI}}_m(\omega) - f(\omega)\}^4] &\leq \mathbb{E}[\|\hat{f}^{\text{NI}}_m - f\|_\infty^4]\\
&\leq \mathbb{E}[\|\hat{f}^{\text{NI}}_m - f_m\|_\infty^4] +\|{f}_m - f\|_\infty^4\\
& \leq \left(M+\frac{\tau_n^2}{\alpha^2}\right)^4+\frac{n^4\tau_n^8}{1 \land \alpha^8} + \left (\sum_{|l|>m} |\sigma_l|\right)^4.
\end{align*} 
Furthermore, it is shown by \cite{kroll2024nonparametric} that  if $\tau_n^2=56 \log^{1+\delta}(n)$, then $\Pr(A_{\text{NI}}^c)\leq n^{-6}$.
Then, 
\begin{equation}
\label{boundeventAcomplkroll}
\mathbb{E}[\mathbbm{I}_{A_{\text{NI}}^c}\{\hat{f}^{\text{NI}}_m(\omega) - f(\omega)\}^2]\leq  \left ( \frac{\tau_n^4}{n} \lor \frac{\tau_n^4}{n\alpha^4} \right).
\end{equation}
Combining  (\ref{boundeventAkroll}) and (\ref{boundeventAcomplkroll}), setting  $ \tau_n^2=56 \log^{1+\delta}(n)$ gives
\begin{align*}\mathbb{E}|\hat{f}^{\text{NI}}_m(\omega)-f(\omega)|^2
&=	\mathcal{O}\left (\frac{\tau_n^4}{n} \right) +\mathcal{O}\left(\frac{ m\lor m\tau_n^{4}\alpha^{-4} }{n}\right)  + \begin{cases}
\mathcal{O}(\log(m)m^{-2s}), &f{\in} W^{s,\infty}(L_0,L),\\
\mathcal{O}(m^{-2s+1}), &f{\in}  W^{s,2}(L).
\end{cases}
%&=&\mathcal{O}\{\log(K)^2K^{-2s}\}+\mathcal{O}\left\{ \frac{\tau_n^6K}{\alpha^ {2}(n-K)}\right\}
\end{align*}
Optimizing  $m^{-2s}+m \cdot (1\lor \tau_n^{4}\alpha^{-4})/n$ with respect to $m$ gives
$$m=\left(\frac{1}{n} \lor\frac{\tau_n^4}{n\alpha^4} \right)^{-\frac{1}{2s+1}}.$$
Analogously, optimizing  $m^{-2s+1}+m \cdot (1\lor \tau_n^{4}\alpha^{-4}) /n$ with respect to $m$ gives
$$m=\left( \frac{1}{n}\lor\frac{\tau_n^4}{n\alpha^4} \right)^{-\frac{1}{2s}}.$$
As $f\in W^{s,\infty}(L_0,L)$ resp. $f\in W^{s,2}(L)$  was chosen arbitrary, we obtain 
\begin{align*}
\sup_{f\in W^{s,\infty}(L_0,L)}\mathbb{E}|\hat{f}^{\text{NI}}_m(\omega)-f(\omega)|^2&\leq c_1 \left\{\frac{\log(n)}{n} \lor \frac{\log^{3+2\delta}(n)}{n \alpha^4} \right\}^{\frac{2s}{2s+1}} \lor \frac{\log^{2+2\delta}(n)}{n}
\end{align*} 
and 
\begin{align*}
\sup_{f\in W^{s,2}(L)}\mathbb{E}|\hat{f}^{\text{NI}}_m(\omega)-f(\omega)|^2&\leq c_2\left\{\frac{1}{n} \lor \frac{\log^{2+2\delta}(n)}{n \alpha^4} \right\}^{\frac{2s-1}{2s}} \lor \frac{\log^{2+2\delta}(n)}{n}.
\end{align*} \hfill \BlackBox
\subsection{Proof of Lemma~\ref{lemma:aLDP_sigma}}
First consider $j$ in $\{1,...,n-1\}$. It is sufficient to check that for all $i=1,...,n$, for all $z$ and $\bar z$ and for all $x$ different from $x'$ the ratio of conditional density functions $q$ of $Q$ is bounded by $e^\alpha$, that is,
\begin{align*}
\frac{q^{Z_i | X_i=x, Z_1, \ldots, Z_{i-1}}(z)}{q^{ Z_i | X_i=x', Z_1, \ldots, Z_{i-1}}(z)} &\leq e^\alpha , \quad \mbox{for } i\leq j , \text{ and }\\
\frac{q^{(Z_i,\bar Z_{i,j})|X_i=x, (Z_1, \bar Z_{1,j}),...,(Z_{i-1},\bar Z_{i-1,j})}(z,\bar z)}{q^{( Z_i,\bar Z_{i,j})|X_i=x', (Z_1, \bar Z_{1,j}),...,(Z_{i-1},\bar Z_{i-1,j})}(z,\bar z)} &\leq e^\alpha, \quad \text{for } i\geq j+1.
\end{align*} 
%where  $ \bar Z_{i,j}=0$ for $i=1,..,K$. Cristina: the problem with putting these variables to 0 is that we should include some dirac measures at 0. It is simpler not to have them for $i\leq K$.
By definition of the privacy mechanism described in \eqref{eq:Zi} and \eqref{eq:barZij}, the conditional density of $(Z_i,\bar Z_{i,j})$ given $X_i=x$ and all the previously built variables  is
\begin{align*}
&q^{Z_i | X_i=x, Z_1, \ldots, Z_{i-1}}(z)  = \frac{\alpha}{8 \tau_n} \exp\left(-\frac{\alpha}{4 \tau_n } |z- \tilde x| \right), \quad \mbox{for } i\leq j , \text{ and }\\
&q^{(Z_i,\bar Z_{i,j})|X_i=x, (Z_1, \bar Z_{1,j}),...,(Z_{i-1},\bar Z_{i-1,j})}(z,\bar z) \\
& =  \frac{\alpha}{8 \tau_n} \exp\left(-\frac{\alpha}{4 \tau_n} |z- \tilde x|\right) \cdot \frac{\alpha}{8 \tilde\tau_n} \exp\left (-\frac{\alpha}{4 \tilde\tau_n} |\bar z- \tilde U_i(x)| \right )  , \quad \text{for } i\geq j+1,
\end{align*}
where $\widetilde{X}$ denotes $x$ trimmed at $\pm \tau_n$ and $\tilde U_i(x)$ denotes $xZ_{i-j}$ trimmed at $\pm \tilde \tau_n$. Thus, the likelihood ratio above can be treated simultaneously by
\begin{align*}
&= \exp \left\{-\frac{\alpha}{4 \tau_n} |z- \tilde x|+\frac{\alpha}{4 \tau_n} |z-\tilde x'| -\left(\frac{\alpha}{4  \tilde \tau_n} |\bar z-\tilde U_i(x)| - \frac{\alpha}{4  \tilde \tau_n} |\bar z-\tilde U_i(x')|\right) \cdot \mathbf{1}(i\geq j+1)  \right\}\\
& \leq  \exp \left(\frac{\alpha}{4 \tau_n} |\tilde x-\tilde x'| + \frac{\alpha}{4  \tilde \tau_n} |\tilde U_i(x)-\tilde U_i(x')| \right).
\end{align*} We conclude after noting that $| x-\tilde x'| \leq 2\tau_n$ and that $|\tilde U_i(x)-\tilde U_i(x')|\leq 2\tilde \tau_n$ by construction. 

If $j=0$, then for all $i=1,...,n$, for all $\bar z$ and for all $x$ different from $x'$, we have
\begin{align*}
&\frac{q^{(\bar Z_{i,0})|X_i=x}(\bar z)}{q^{(\bar Z_{i,0})|X_i=x'}(\bar z)}=\exp \left\{ -\frac{\alpha}{2   \tau_n} |\bar z - \tilde U_i(x)| + \frac{\alpha}{2  \tau_n} |\bar z-\tilde U_i(x')| \right\}\leq  \exp \left(\frac{\alpha}{2 \tau_n} |\tilde U_i(x)-\tilde U_i(x')| \right),
\end{align*} where $\tilde U_i(x)$ denotes $x^2$ trimmed at $\pm \tilde \tau_n$. We conclude after noting that  $|\tilde U_i(x)-\tilde U_i(x')|\leq 2 \tau_n$ by construction. \hfill \BlackBox

%%%%%%%%%%%%%%%%%%%%%%%%%%%%%%%%%%%%%%%%%%%%%%%%%%%%%%%%%%%%%%%%%%
%SECTION: Interactive privacy mechanism for $\sigma_j$
%%%%%%%%%%%%%%%%%%%%%%%%%%%%%%%%%%%%%%%%%%%%%%%%%%%%%%%%%%%%%%%%%%

%%%%%%%%%%%%%%%%%%%%%%%%%%%
\subsection{Proof of Theorem~\ref{theo:ratesigma}} 
Let $R\in\mathbb{N}$  such that $0<R<n$ and $R=o(n)$. We start with $j\in\{1,...,R\}$. Consider the event $A$ and its complement defined by
\begin{align*}
A&{=} \bigcap_{i{=}1}^n \{ \widetilde{X}_i{=} X_i \text{ and } (\widetilde W_{i,j} {=} W_{i,j})   \mathbf{1}(i>j) \}, &  A^c&{=}\bigcup_{i{=}1}^n \{ \widetilde{X}_i{\neq} X_i \text{ or } (\widetilde W_{i,j} {\neq} W_{i,j})  \mathbf{1}(i>j)\}. 
\end{align*}
\emph{$L_2$ risk on the event $A$.} If the event $A$ occurs, then $\hat{\sigma}_j$ coincides with the untruncated estimator $\tilde{\sigma}_j$, which is defined by
\begin{align*}
\tilde{\sigma}_j&=\frac{1}{n-j} \sum_{i=j+1}^n  X_i \cdot \left ( X_{i-j}+\xi_{i-j}\right) +\tilde{\xi}_i.
\end{align*}
In particular,
\begin{align*}\mathbb{E}[(\hat{\sigma}_j - \sigma_j)^2] = \mathbb{E}[\mathbb{I}_A(\hat{\sigma}_j - \sigma_j)^2] + \mathbb{E}[\mathbb{I}_{A^c}(\hat{\sigma}_j - \sigma_j)^2] % \label{decompose2}
\end{align*}
\begin{align*}
\text{and } \quad \mathbb{E}[\mathbb{I}_A(\hat{\sigma}_j - \sigma_j)^2] &\leq \mathbb{E}[(\tilde{\sigma}_j - \sigma_j)^2] %=\mathbb{E}[\{\tilde{\sigma}_j - \mathbb{E}[\tilde{\sigma}_j] + \{\mathbb{E}[\tilde{\sigma}_j] - {\sigma}_j\}^2] \\
= \mathbb{E}[(\tilde{\sigma}_j - \mathbb{E}[\tilde{\sigma}_j ])^2]  + (\mathbb{E}[\tilde{\sigma}_j ] - {\sigma}_j)^2.
\end{align*}
It is straightforward to see that 
\begin{align*}
\mathbb{E}[ \tilde{\sigma}_j]& 
=  \frac{1}{n-j} \sum_{i=j+1}^n \mathbb{E} \left[ X_i X_{i-j}\right ] =\frac{1}{n-j} \sum_{i=j+1}^n \sigma_j = \sigma_j.
\end{align*} 
Let $\bar{\sigma}_j = \frac{1}{n-j} \sum_{i=j+1}^n   X_{i}X_{i-j}$ be the unbiased estimator built using the data. Then,
\begin{align}
\tilde{\sigma}_j&=\frac{1}{n-j} \sum_{i=j+1}^n  X_i \cdot \left ( X_{i-j}+\xi_{i-j}\right) +\tilde{\xi}_i= \bar{\sigma}_j + \frac{1}{n-j}\sum_{i=j+1}^n  X_i\xi_{i-j} + \frac{1}{n-j}\sum_{i=j+1}^n  \tilde\xi_{i} \nonumber\\
&=:  \bar{\sigma}_j+B+C. \nonumber
\end{align} 
Since all terms of $B$ and $C$ are centered, uncorrelated and uncorrelated to the terms in $\bar \sigma_j$, we have
\begin{equation}
\mathbb{E}\left[(\tilde{\sigma}_j - \mathbb{E}[\tilde{\sigma}_j ])^2 \right ]=\mathrm{Var}[\bar{\sigma}_j]+\mathbb{E}[B^2]+\mathbb{E}[C^2]. \label{sigmaABC}
\end{equation}
By Wick's formula, see Equation~(\ref{Wick}), and similar computations as in the proof of Proposition~\ref{prop:ratekroll}$i)$, we have
\begin{align*}
\mathrm{Var}[\bar{\sigma}_j]&= \frac{1}{(n-j)^2} \sum_{i=j+1}^n \sum_{k=j+1}^n  \mbox{Cov} \left (   X_{i}X_{i-j} , X_{k}X_{k-j} \right)
%&=\frac{1}{(n-j)^2} \sum_{i=j+1}^n \sum_{k=j+1}^n \sigma_{|i-k|}^2 + \sigma_{|i-k+j|} \sigma_{|i-k-j|} \\
%&=\frac{1}{(n-j)^2} \sum_{h=-(n-j-1)}^{n-j-1}  (n-j-|h|)(\sigma_{|h|}^2 + \sigma_{|h+j|} \sigma_{|h-j|})\\
%&=\frac{1}{(n-j)} \sum_{h=-(n-j-1)}^{n-j-1}  \left (1-\frac{|h|}{n-j}\right)(\sigma_{|h|}^2 + \sigma_{|h+j|} \sigma_{|h-j|})\\
%&\leq\frac{1}{(n-j)} \left ( \sum_{h=-(n-j-1)}^{n-j-1} \sigma_{|h|}^2  + \sum_{h=-n}^{n}\sigma_{|h|}^2 \right )\\
\leq\frac{2}{(n-j)}    \sum_{h=-n}^{n}\sigma_{|h|}^2  \leq\frac{2M_1}{(n-j)}.
\end{align*}
Next, since the $ \xi_i$ are independent  of the $X_i$ and centered,
\begin{align*}
\mathbb{E}[B^2]&=  \frac{1}{(n-j)^2}\sum_{i=j+1}^n \sum_{k=j+1}^n  \mathbb{E}[X_iX_k] \mathbb{E}[\xi_{i-j}\xi_{k-j}]=  \frac{8\sigma_0 \tau_n^2}{\alpha^2(n-j)}.
\end{align*}
For the last term in (\ref{sigmaABC}) we have
$$
\mathbb{E}[C^2]=\frac{1}{(n-j)^2}\sum_{i=j+1}^n  \mathbb{E}[\tilde\xi_{i}^2]= \frac{8 \tilde \tau_n^2}{\alpha^2(n-j)}.
$$
Then,
\begin{align*}\mathbb{E}[(\hat{\sigma}_j - \sigma_j)^2] &\leq \frac{2M_1}{(n-j)} + \frac{8M \tau_n^2}{\alpha^2(n-j)}+\frac{8 \tilde \tau_n^2}{\alpha^2(n-j)} + \mathbb{E}[\mathbb{I}_{A^c}(\hat{\sigma}_j - \sigma_j)^2].%\\
%& \leq C\max \left \{\frac{(\tau_n^2,\tilde \tau_n^2)}{\alpha^2(n-K)}  + \mathbb{E}[\mathbb{I}_{A^c}(\hat{\sigma}_j - \sigma_j)^2].
\end{align*}
To conclude the proof, we show that $\mathbb{E}[\mathbb{I}_{A^c}(\hat{\sigma}_j - \sigma_j)^2]$ is not larger than the other terms.\\  %$= \mathcal{O}(\frac{\max(\tau_n^2,\tilde \tau_n^2)}{\alpha^2(n-j)}).$ \\
\emph{$L_2$ risk on the event $A^c$.} 
Application of the Cauchy-Schwarz inequality gives
\begin{align*} %\label{CSbound}
\mathbb{E}[\mathbb{I}_{A^c}(\hat{\sigma}_j - \sigma_j)^2] \leq \left(\mathbb E[(\hat{\sigma}_j - \sigma_j)^4] \cdot \Pr[ A^c ]\right)^{1/2} .
\end{align*}
Then, using $|\widetilde W_{i,j}|\leq  \tilde \tau_n$, it holds for some constants $c_1,...,c_5>0$, that
\begin{align*}
\mathbb{E}[(\hat{\sigma}_j - \sigma_j)^4] &\leq c_1\mathbb{E}[\hat{\sigma}_j^4] +c_2 \leq \frac{c_3}{(n-j)^4} \left\{\mathbb{E}\left[ \left(\sum_{i=j+1}^n\widetilde W_{i,j}\right)^4\right ]+  \mathbb{E}\left[\left(\sum_{i=j+1}^n \tilde{\xi}_{i} \right)^4 \right]\right\} +c_2\\
&\leq  c_4\tilde \tau_n^{4} +  \frac{c_5\tilde\tau_n^{4}}{\alpha^4} +c_2=\mathcal{O}(\tilde \tau_n^{4}+\tilde\tau_n^{4}\alpha^{-4}).
\end{align*} 
%It remains to show that $$\Pr[A^c]=\mathcal{O}\{(n-K)^{-2}\}.$$
Moreover,
\begin{align*}
\Pr[A^c] &\leq \sum_{i=1}^n \Pr[\{ \widetilde{X}_i\neq X_i \text{ or }\{ \widetilde  W_{i,j} \not = W_{i,j}\}\cdot \mathbf{1}(i>j)\}] \\ 
&\leq \sum_{i=1}^n \Pr[ \widetilde{X}_i\neq X_i ] + \Pr[\{ \widetilde{X}_i = X_i \text{ and }  \{\widetilde  W_{i,j} \not = W_{i,j}\}\cdot \mathbf{1}(i>j)\}].
\end{align*}
We start with the first term above and use that $\{ \widetilde{X}_i\neq X_i \} = \{ |X_i| > \tau_n \} $. As $X_i\sim\mathcal{N}(0,\sigma_0)$ we have
\begin{equation} \label{Gausstail}
\sum_{i=1}^n \Pr[ \widetilde{X}_i\neq X_i ] \leq 2n\exp\left ( -\frac{\tau_n^{2}}{2\sigma_0} \right),
\end{equation} which is of the order $\mathcal{O}(n^{-3})$ for $\tau_n^2\geq 8\log^{1+\delta}(n)$ for some $\delta>0$. 
For $i>j$, we have
\begin{align*}
& \sum_{i=j+1}^n \Pr[\{ \widetilde{X}_i = X_i \text{ and }  \{\widetilde  W_{i,j} \not = W_{i,j}\}\cdot \mathbf{1}(i>j)\}]\\
&\leq \sum_{i=j+1}^n\Pr \left[\left |\widetilde X_i \cdot Z_{i-j}\right |\geq \tilde\tau_n \right ] \leq n\Pr \left[\left |\widetilde{X}_{i-j}+\xi_{i-j}\right |\geq \tilde\tau_n/\tau_n \right ]\\ \nonumber 
&\leq C n\Pr \left[\left | \xi_{i-j}\right |\geq \tilde\tau_n/\tau_n \right ] \leq  \frac{Cn\alpha}{4\tau_n}\exp\left ( - \frac{\tilde \tau_n \alpha}{4\tau_n^2} \right)
\end{align*}
as the subexponential tail is dominating and is of the order $\mathcal{O}(n^{-3})$ if $\tilde \tau_n\tau_n^{-2}\geq 16\log^{1+\delta}(n).$  
Setting $ \tau_n^2=8\log^{1+\delta}(n)$ and $\tilde \tau_n=16\log^{1+\delta}(n) \tau_n^2$ for some $\delta>0$ finally gives for large $n$ that
$$\mathbb{E}[(\hat{\sigma}_j - \sigma_j)^2]\leq c \left \{ \frac{1}{n} \lor\frac{\log^{4+4\delta}(n)}{\alpha^2n} \right\}  .$$
In particular, we used that for $j\leq R$ and $R=o(n)$ holds $1/(n-j)=\mathcal{O}(n)$. Similar computations for $j=0$. As the constants do not depend on $(\sigma_k)_{k\in\mathbb{Z}}$ but only on the global bounds $\sum_{j\in\mathbb{Z}}|\sigma_j|<M$ and $\sum_{j\in\mathbb{Z}}\sigma_j^2\leq M_1$, we obtain $$\sup_{\sigma:\,\|\sigma\|_2^2\leq M_1}\mathbb{E}[(\hat{\sigma}_j - \sigma_j)^2]\leq c \left \{ \frac{1}{n} \lor\frac{\log^{4+4\delta}(n)}{\alpha^2n} \right\} .$$ \hfill \BlackBox

%%%%%%%%%%%%%%%%%%%%%%%%%%%%%%%%%%%%%%%%%%%%%%%%%%%%%%%%%%%%%%%%%%
%SECTION: Interactive privacy mechanism for $f(\omega)$
%%%%%%%%%%%%%%%%%%%%%%%%%%%%%%%%%%%%%%%%%%%%%%%%%%%%%%%%%%%%%%%%%%

%%%%%%%%%%%%%%%%%%%%%%%%%%%%%%%%%
\subsection{Proof of Lemma~\ref{lemma:aLDP_sdf0}} In view of Lemma~\ref{lemma:aLDP_sigma} the case $i\leq K$ is proven. It is therefore sufficient to check that for all $i=K+1,\ldots,n$, for all $z$ and $\tilde z$ and for all $x$ different from $x'$ the ratio of conditional density functions $q$ of $Q$ is bounded, that is,
\begin{align*}
\frac{q^{(Z_i,\widetilde Z_i)|X_i=x, (Z_{1},\widetilde Z_{1}),\ldots , (Z_{i-1},\widetilde Z_{i-1})}(z,\tilde z)}{q^{(Z_i,\widetilde Z_i)|X_i=x', (Z_{1},\widetilde Z_{1}),\ldots , (Z_{i-1},\widetilde Z_{i-1})}(z,\tilde z)} \leq e^\alpha.
\end{align*}
By definition of the privacy mechanism described in \eqref{eq:Zi} and \eqref{eq:tildeZi}, the conditional density of $(Z_i,\widetilde Z_i)$ given $X_i=x$ and all the previously built variables  is
\begin{align*}
& q^{(Z_i,\widetilde Z_i)|X_i=x, (Z_1,\widetilde Z_1),\ldots , (Z_{i-1},\widetilde Z_{i-1})}(z,\tilde z)\\
& = %\frac{\alpha}{8 \tau_n} \exp\left(-\frac{\alpha}{4 \tau_n } |z- \tilde x| \right) \cdot I(i\leq K)  \\
%&\quad + 
\frac{\alpha}{8 \tau_n} \exp\left(-\frac{\alpha}{4 \tau_n} |z- \tilde x|\right) \cdot \frac{\alpha}{8 \tilde\tau_n} \exp(-\frac{\alpha}{4 \tilde\tau_n} |\tilde z-\tilde V_i(x,\omega)|),  \quad \text{for } i\geq K+1 ,
\end{align*}
where $\widetilde{X}$ denotes $x$ trimmed at $\pm \tau_n$ and $\tilde V_i(x)$ denotes $x^2+\sum_{1\leq |k|\leq K} a_kxZ_{i-|k|}\exp(\texttt{i}\omega k)$ trimmed at $\pm \tilde \tau_n$. Thus,
\begin{align*}
& \frac{q^{(Z_i,\widetilde Z_i)|X_i=x, (Z_1,\widetilde Z_1),\ldots , (Z_{i-1},\widetilde Z_{i-1})}(z,\tilde z)}{q^{(Z_i,\widetilde Z_i)|X_i=x', (Z_1,\widetilde Z_1),\ldots , (Z_{i-1},\widetilde Z_{i-1})}(z,\tilde z)} \\
&= \exp \left\{-\frac{\alpha}{4 \tau_n} |z- \tilde x|+\frac{\alpha}{4 \tau_n} |z-\tilde x'| -\left(\frac{\alpha}{4  \tilde \tau_n} |\tilde z-\tilde V_i(x)| - \frac{\alpha}{4  \tilde \tau_n} |\tilde z-\tilde V_i(x')|\right)  \right\}\\
& \leq  \exp \left(\frac{\alpha}{4 \tau_n} |\tilde x-\tilde x'| + \frac{\alpha}{4  \tilde \tau_n} |\tilde V_i(x)-\tilde V_i(x')| \right).
\end{align*}
We conclude by noting that $| x-\tilde x'| \leq 2\tau_n$ and that $|\tilde V_i(x)-\tilde V_i(x')|\leq 2\tilde \tau_n$.  \hfill \BlackBox

%%%%%%%%%%%%%%%%%%%%%%%%%%%%%%%%%%%%%%%%%%%%%%%%%%%%%%%%%%%%%%%%%%
%SECTION: Interactive privacy mechanism for $f(\omega)$
%%%%%%%%%%%%%%%%%%%%%%%%%%%%%%%%%%%%%%%%%%%%%%%%%%%%%%%%%%%%%%%%%%

%%%%%%%%%%%%%%%%%%%%%%%%%%%
\subsection{Proof of Theorem~\ref{theo:ratef0}} 
%Note that $f\in W^{s,\infty}(L_0,L)$ implies that there exists a constant $\tilde L$ such that $f\in W^{s,2}(\tilde L)$.
Since the proofs for the cases $f\in W^{s,\infty}(L_0,L)$ and for $f\in W^{s,2}(L)$ are identical apart from the treatment of the bias term, we handle both cases simultaneously. 
Let $\omega\in[-\pi,\pi]$ be arbitrary but fixed. Consider the event $\tilde A$ and its complement defined by
\begin{align*}
\tilde A&{=} \bigcap_{i=1}^n \{ \widetilde{X}_i{=} X_i \text{ and } (\widetilde V_i {=} V_i)  \mathbf{1}(i>K)\}, & \tilde A^c&{=}\bigcup_{i=1}^n \{ \widetilde{X}_i{\neq} X_i \text{ or } (\widetilde V_i{\neq} V_i)  \mathbf{1}(i>K)\}.
\end{align*} 
\emph{$L_2$ risk on the event $\tilde A$.} 
If the event $\tilde A$ occurs, then $\hat{f}_K(\omega)$ coincides with the untruncated estimator $\tilde{f}_K(\omega)$, which is defined by 
\begin{align*}
\tilde{f}_K(\omega)&=\frac{1}{2 \pi (n-K)} \sum_{j=K+1}^n  \left ( X_{j}^2+\sum_{1\leq |k|\leq K} a_kX_{j}(X_{j-|k|}+\xi_{j-|k|})\exp(-\texttt{i}\omega k) +\tilde{\xi}_{j}\right).
\end{align*}
In particular,
\begin{align*}\quad \quad \,\mathbb{E}[\{\hat{f}_K(\omega) - f(\omega)\}^2]&= \mathbb{E}[\mathbb{I}_{\tilde A}\{\hat{f}_K(\omega) - f(\omega)\}^2] + \mathbb{E}[\mathbb{I}_{\tilde A^c}\{\hat{f}_K(\omega) - f(\omega)\}^2] %\label{f0decompose}
\end{align*} and
\begin{align*}
\mathbb{E}[\mathbb{I}_{\tilde A}\{\hat{f}_K(\omega) - f(\omega)\}^2]&\leq \mathbb{E}[\{\tilde{f}_K(\omega) - f(\omega)\}^2]\\  %= \mathbb{E}[\{\tilde{f}_K(\omega) - \mathbb{E}[\tilde{f}_K(\omega)] + \{\mathbb{E}[\tilde{f}_K(\omega)] - f(\omega)\}^2] \\
&= \mathbb{E}[\{\tilde{f}_K(\omega) - \mathbb{E}[\tilde{f}_K(\omega)]\}^2]  + \{\mathbb{E}[\tilde{f}_K(\omega)] - f(\omega)\}^2.
\end{align*}
%\subsubsection{The Bias $ \{\mathbb{E}[\tilde{f}_K(\omega)] - f(\omega)\}^2$ }
We start with calculating the expectation of $\tilde{f}_K(\omega)$:
\begin{align*}
\mathbb{E}[ \tilde{f}_K(\omega)]& %\frac{1}{\pi(n-k)} \sum_{j=K+1}^n \mathbb{E} \left [ \sum_{|k|\leq K} X_j Z_{j-k} \exp(\texttt{i}\omega k) + \tilde{\xi}_j \right ] 
=  \frac{1}{2\pi(n-K)} \sum_{j=K+1}^n \left( \sum_{|k|\leq K}  a_k\mathbb{E} [ X_j X_{j-|k|} ] \exp(-\texttt{i}\omega k) + \mathbb{E}  [\tilde{\xi}_j  ] \right)\\
%=\mathbb{E}[ \bar{f}_K(\omega)]
&=\frac{1}{2\pi(n-K)} \sum_{j=K+1}^n \sum_{|k|\leq K}  a_k\sigma_k \exp(-\texttt{i}\omega k) = \frac{1}{2\pi} \sum_{|k|\leq K}  a_k\sigma_k \exp(-\texttt{i}\omega k).
\end{align*} 
As the last sum is the De la Valée-Poussin sum of $f$, Theorem 13.5 (see Chapter 3, \citealp{zygmund2002trigonometric}) implies
$$\left | \frac{1}{2\pi}\sum_{|k|\leq K}  a_k\sigma_k \exp(-\texttt{i}\omega k) -f(\omega)\right|\leq 4 \inf_{T_K \in \mathcal{T}_K} \|f-T_K\|_\infty ,$$
where $\mathcal{T}_K$ is the class of trigonometric polynomials of degree less or equal than $K$. 
By Jackson's inequality (see Theorem 13.6 in Chapter 3, \citealp{zygmund2002trigonometric}), we have
$$\inf_{T_K \in \mathcal{T}_K} \|f-T_K\|_\infty \leq C_s\omega(f^{s-\lfloor s\rfloor},2\pi/K)K^{-\lfloor s\rfloor},$$
where the constant $C_s>0$ depends only on  the smoothness $s$ of $f$, $\lfloor s\rfloor$ is the largest integer strictly less than $s$ and $\omega(f,\delta)=\max_{ |t_1-t_2|\leq \delta}|f(t_1)-f(t_2)|$ is the modulus of continuity. In particular,  $\omega(f^{s-\lfloor s\rfloor},2\pi/K)\leq L_0 K^{-(s-\lfloor s\rfloor)}$ and $\inf_{T_K \in \mathcal{T}_K} \|f-T_K\|_\infty\leq C^\prime_sK^{-s}$ for $f\in W^{s,\infty}(L_0,L)$. As $W^{s,2}(L)$ can be embedded in $W^{s-1/2,\infty}(\tilde{L}_0,\tilde{L})$ for some constants $\tilde{L}_0,\tilde{L}>0$ depending on $L$, it follows $\inf_{T_K \in \mathcal{T}_K} \|f-T_K\|_\infty\leq \tilde C_sK^{-(s-1/2)}$ for $f\in W^{s,2}(L)$. Thus, we obtain for the squared bias 
\begin{equation*}
\{\mathbb{E}[ \tilde{f}_K(\omega)]-f(\omega)\}^2\leq C \cdot
\begin{cases} 
K^{-2s}, & f\in W^{s,\infty}(L_0,L),\\
K^{-2s+1}, &f\in W^{s,2}(L),
\end{cases}
\end{equation*} where the constant $C>0$ depends on $s,\,L$, if $f\in W^{s,2}(L)$, and additionally on $L_0$, if $ f\in W^{s,\infty}(L_0,L).$

To analyze the variance, define the estimator without truncation or added Laplace noise by
$$
\bar{f}_K(\omega) = \frac{1}{2\pi (n-K)} \sum_{j=K+1}^n   \sum_{|k|\leq K} a_kX_{j}X_{j-|k|}\exp(-\texttt{i}\omega k).
$$
Then,
\begin{align*}
\tilde{f}_K(\omega)=&\frac{1}{2\pi (n-K)} \sum_{j=K+1}^n  \left ( X_{j}^2+\sum_{1\leq |k|\leq K} a_kX_{j}(X_{j-|k|}+\xi_{j-|k|})\exp(-\texttt{i} \omega k) +\tilde{\xi}_{j}
\right)\\
=&\frac{1}{2\pi (n-K)} \sum_{j=K+1}^n   \sum_{|k|\leq K} a_kX_{j}X_{j-|k|}\exp(-\texttt{i}\omega k)\\
& +\frac{1}{2\pi (n-K)} \sum_{j=K+1}^n   \sum_{1\leq |k|\leq K} a_kX_{j}\xi_{j-|k|}\exp(-\texttt{i} \omega k) +\frac{1}{2\pi (n-K)} \sum_{j=K+1}^n   \tilde{\xi}_{j}\\
=:& \bar{f}_K(\omega)+B+C.
\end{align*} 
Since all mixed terms have zero expectation and $\mathbb{E}[\tilde{f}_K(\omega)]=\mathbb{E}[\bar{f}_K(\omega)]$, we obtain
\begin{equation}
\mathbb{E}[\{\tilde{f}_K(\omega) - \mathbb{E}[\tilde{f}_K(\omega)]\}^2] =\mathrm{Var}[\bar{f}_K(\omega)]+\mathbb{E}[B^2]+\mathbb{E}[C^2]. \label{f0ABC}
\end{equation}
In particular, using $ \bar{\sigma}_k = \frac{1}{n-K} \sum_{j=K+1}^n   X_{j}X_{j-|k|}$, we have
\begin{align}
&(2\pi (n-K))^2   \mathrm{Var}[\bar{f}_K(\omega)] \nonumber\\
&= (n-K)^2  \sum_{l=-K}^K\sum_{k=-K}^K a_la_k\exp\{-\texttt{i}\omega (l-k)\}\text{Cov}(\bar{\sigma}_l, \bar{\sigma}_k)\nonumber\\
&= \sum_{l=-K}^K\sum_{k=-K}^K a_la_k\exp\{-\texttt{i}\omega (l-k)\}\sum_{r=K+1}^n \sum_{s=K+1}^n \text{Cov}(X_rX_{r-|l|},X_{s}X_{s-|k|}). \label{Varf0calc}
%&\sum_{j=-K}^K\sum_{k=-K}^K a_ja_k\exp\{\textt{i}\omega (j-k)\} \sum_{r=K+1}^n \sum_{s=K+1}^n 
% &\leq M_0 K/n.
\end{align}
Now,
\begin{align*}
&\text{Cov}(X_rX_{r-|l|},X_{s}X_{s-|k|})\\
&=\mathbb E[X_rX_{r-|l|}]\mathbb E[X_{s}X_{s-|k|}]+\mathbb E[X_rX_{s}] \mathbb E[X_{s-|k|}X_{r-|l|}]+\mathbb E[X_rX_{s-|k|}]\mathbb E[X_{s}X_{r-|l|}]-\sigma_l\sigma_k\\
&=\sigma_{|r-s|}\sigma_{|r-s-|k|+|l||}+\sigma_{|r-s+|k||} \sigma_{|r-s-|l||}.
\end{align*}
Substituting this into \eqref{Varf0calc}, we get 
\begin{align*}
&(2\pi (n-K))^2   \mathrm{Var}[\bar{f}_K(\omega)]\\
&= \sum_{l=-K}^K\sum_{k=-K}^K a_la_k\exp\{-\texttt{i}\omega (l-k)\}\sum_{h=-(n-K-1)}^{n-K-1}  (n{-}K{-}|h|)(\sigma_{|h|}\sigma_{|h-|k|+|l||}{+}\sigma_{|h+|k||} \sigma_{|h-|l||})\\
&=\sum_{l=-K}^K a_l\exp({-}\texttt{i}\omega l)\sum_{h=-(n-K-1)}^{n-K-1}  (n{-}K{-}|h|)\sum_{k=-K}^K a_k (\sigma_{|h|}\sigma_{|h-k+|l||}{+}\sigma_{|h+|k||} \sigma_{|h-|l||}) \exp(\texttt{i}\omega k)\\
&\leq \sum_{l=-K}^K a_l\sum_{h=-(n-K-1)}^{n-K-1}  (n-K-|h|)(|\sigma_{|h|}|+| \sigma_{|h-|l||}|)M. 
% \leq 2M^2(2K+1).
\end{align*}
Since $\sum_{h=-(n-K-1)}^{n-K-1}  \left (1-\frac{|h|}{n-K} \right)(|\sigma_{|h|}|+| \sigma_{|h-|l||}|) \leq 2M$, it follows $\mathrm{Var}[\bar{f}_K(\omega)]\leq CK(n-K)^{-1}.$
For the last term in Equation (\ref{f0ABC}) we have by independence of the $\xi_i$ that
$$
\mathbb{E}[C^2]=\frac{1}{4\pi^2(n-K)^2}\sum_{j=K+1}^n\mathbb{E}(\tilde\xi_j^2)=\frac{8\tilde \tau_n}{\alpha^2\pi^2(n-K)}.
$$
It remains to calculate $\mathbb{E}[B^2]$. For $j=K+1,...,n$, we have
\begin{align*}
&\sum_{1\leq|k|\leq K} a_kX_{j}\xi_{j-|k|}\exp(-\texttt{i}\omega k)
=  \sum_{k=1}^K a_kX_{j}\xi_{j-k}\exp(-\texttt{i}\omega k)+ \sum_{k=-K}^{-1} a_kX_{j}\xi_{j+k}\exp(-\texttt{i}\omega k)\\
&= \sum_{j-k=j-K}^{j-1} a_kX_{j}\xi_{j-k}\exp(-\texttt{i}\omega k)+ \sum_{j+k=j-K}^{j-1} a_kX_{j}\xi_{j+k}\exp(-\texttt{i}\omega k)\\
&= \sum_{s=j-K}^{j-1} a_{-s+j}X_{j}\xi_{s}\exp\{- \texttt{i} \omega (-s+j)\}+ \sum_{t=j-K}^{j-1} a_{t-j} X_{j}\xi_{t}\exp\{ -\texttt{i} \omega (t-j)\}.
%&=&\sum_{j=K+1}^nX_j\sum_{m=K+j}^{-K+j}\xi_m\exp\{-\texttt{i}\omega(m-j)\}=\sum_{j=K+1}^nX_j\exp(\texttt{i}\omega j)\sum_{m=K+j}^{-K+j}\xi_m\exp(-\texttt{i}\omega m).
\end{align*}
Thus,
\begin{align*}
&2\pi(n-K)B=\sum_{j=K+1}^n \sum_{1\leq|k|\leq K} a_kX_{j}\xi_{j-|k|}\exp(-\texttt{i}\omega k)\\
&= \sum_{j=K+1}^nX_{j}\exp\{ -\texttt{i} \omega j\}\sum_{s=j-K}^{j-1} a_{-s+j}\xi_{s}\exp\{ \texttt{i} \omega s\}\\
&\quad +  \sum_{j=K+1}^nX_{j}\exp\{ \texttt{i} \omega j\}\sum_{t=j-K}^{j-1} a_{t-j}\xi_{t}\exp\{ -\texttt{i} \omega t\} =: B_1 +B_2.
\end{align*}
Furthermore,
\begin{align*}
|B_1|^2&{=}    \sum_{j=K+1}^n\sum_{h=K+1}^nX_{j}X_{h}\exp\{{ -}\texttt{i} \omega (j{-}h)\}\sum_{s=j-K}^{j-1}\sum_{k=h-K}^{h-1} a_{-s+j}a_{-k+h}\xi_{s}\xi_{k}\exp\{ {-}\texttt{i} \omega (k{-}s)\}, \\
|B_2|^2&{=}  \sum_{j=K+1}^n\sum_{h=K+1}^nX_{j}X_{h}\exp\{ {-}\texttt{i} \omega (h{-}j)\}\sum_{t=j-K}^{j-1}\sum_{k=h-K}^{h-1} a_{t-j}a_{k-h} \xi_{t}\xi_{k}\exp\{ {-}\texttt{i} \omega (t{-}k)\}.
\end{align*}
Using $4\pi^2(n-K)^2 |B|^2\leq 2(|B_1|^2+|B_2|^2)$,  we obtain 
\begin{align*}
&2\pi^2(n-K)^2\mathbb{E}|B|^2\\
&=\sum_{j=K+1}^n\sum_{h=K+1}^n\mathbb{E}[X_jX_h]\exp\{{-}\mathtt{i}\omega(j{-}h)\}\sum_{s=j-K}^{j-1}\sum_{k=h-K}^{h-1}  a_{-s+j}a_{-k+h} \mathbb{E}[\xi_{s}\xi_{k}]\exp\{ {-}\texttt{i} \omega (k{-}s)\}\\
&\quad \,+\sum_{j=K+1}^n\sum_{h=K+1}^n\mathbb{E}[X_{j}X_{h}]\exp\{ {-}\texttt{i} \omega (h{-}j)\}\sum_{t=j-K}^{j-1}\sum_{k=h-K}^{h-1} a_{t-j}a_{k-h}  \mathbb{E}[\xi_{t}\xi_{k}]\exp\{ {-}\texttt{i} \omega (t{-}k)\}\\
&=\sum_{j=K+1}^n\sum_{h=K+1}^n\sigma_{|j-h|}\exp\{-\mathtt{i}\omega(j-h)\}\sum_{\substack{s=j-K,...,j-1\\k=h-K,...,h-1\\s=k}}a_{-s+j}a_{-k+h}\frac{32 \tau_n^{2}}{\alpha^2}\\
&\quad \,+\sum_{j=K+1}^n\sum_{h=K+1}^n\sigma_{|j-h|}\exp\{-\mathtt{i}\omega(j-h)\}\sum_{\substack{t=j-K,...,j-1\\k=h-K,...,h-1\\t=k}}a_{t-j}a_{k-h} \frac{32 \tau_n^{2}}{\alpha^2}\\
&\leq2\sum_{j=K+1}^n\sum_{h=K+1}^n\left |\sigma_{|j-h|} \right |\max\{0,K-|j-h|\}\frac{32 \tau_n^{2}}{\alpha^2}\\
&\leq\frac{64 \tau_n^{2}(K+1)}{\alpha^2}\sum_{j=K+1}^n\sum_{h=K+1}^n\left |\sigma_{|j-h|} \right| =\mathcal{O}\left\{32\tau_n^{2}K(n-K)\alpha^ {-2}\right\}%
%&=&.
\end{align*}
as $|a_{-s+j}a_{-k+h}|\leq1$ and $\sum_{\substack{s=j-K,...,j-1\\k=h-K,...,h-1\\s=k}}\frac{32 \tau_n^{2}}{\alpha^2}=\max\{0,K-|j-h|\}\frac{32 \tau_n^{2}}{\alpha^2}$. Furthermore, we used that $\sum_k|\sigma_k|<M$, see (\ref{constants}).
Summing up, we get
\begin{equation} \label{boundeventA}
\mathbb{E}[\{\tilde{f}_K(\omega) - f(\omega)]\}^2]= \mathcal{O}\left (\frac{K}{n} \right) +\mathcal{O}\left\{\frac{\max(\tau_n^{2}K,\tilde\tau_n^{2}) }{(n-K)\alpha^ {2}}\right\}  + \mathbb{E}[\mathbb{I}_{A^c}\{\hat{f}_K(\omega) - f(\omega)\}^2] .
\end{equation}
To conclude the proof, we show that $\mathbb{E}[\mathbb{I}_{\tilde A^c}\{\hat{f}_K(\omega) - f(\omega)\}^2]$ is not larger than the other terms.\\
%Next, we show that $\mathbb{E}[\mathbb{I}_{A^c}\{\hat{f}_K(\omega) - f(\omega)\}^2]=\mathcal{O}\{\tilde\tau_n^{2} (n-K)^{-1}\alpha^{-2}\} $.
\emph{$L_2$ risk on the event $\tilde A^c$.} Application of the Cauchy-Schwarz inequality gives
\begin{align*} %\label{f0CSbound}
\mathbb{E}[\mathbbm{I}_{\tilde A^c}\{\hat{f}_K(\omega) - f(\omega)\}^2] \leq \left(\mathbb E[\{\hat{f}_K(\omega) - f(\omega)\}^4] \cdot \Pr[ \tilde A^c ]\right)^{1/2} .
\end{align*}
Using $\|f\|_\infty<M$ and $|\widetilde{V}_j|\leq \tilde \tau_n$ yields for some constants $c_1,...,c_5>0$
\begin{align*}
\mathbb{E}[\{\hat{f}_K(\omega) - f(\omega)\}^4] &\leq c_1\mathbb{E}[\hat{f}_K(\omega)^4] +c_2 \\
&\leq \frac{c_3}{\{2\pi(n-K)\}^4} \left\{\mathbb{E}\left[ \left(\sum_{j=K+1}^n\widetilde V_j \right)^4\right ]+  \mathbb{E}\left[ \left(\sum_{j=K+1}^n \tilde{\xi}_{j} \right)^4\right] \right\} +c_2\\
&\leq  c_4 \tilde \tau_n^{4} +\frac{c_5\tilde \tau_n^{4} }{ \alpha^{4}} +c_2 =\mathcal{O}(\tilde \tau_n^{4} +\tilde\tau_n^{4}\alpha^{-4}).
\end{align*} 
Moreover,
\begin{align*}
\Pr[\tilde A^c] &\leq \sum_{i=1}^n  \Pr[\{ \widetilde{X}_i\neq X_i \text{ or }( \widetilde  V_i\not = V_i)\cdot \mathbf{1}(i>K)\}] \\ 
&\leq \sum_{i=1}^n \Pr[ \widetilde{X}_i\neq X_i ] + \Pr[\{ \widetilde{X}_i = X_i \text{ and }  (\widetilde V_i \not = V_i)\cdot I(i>K)\}].
\end{align*}
From (\ref{Gausstail}), we know that $\sum_{i=1}^n  \Pr[ \widetilde{X}_i\neq X_i ]=\mathcal{O}(n^{-3})$ for $\tau_n^2\geq 8\log^{1+\delta}(n)$. Next, we need to bound from above, for $i>K$ :
\begin{align}
&  \Pr[\{ \widetilde{X}_i = X_i \text{ and }\widetilde V_i \not = V_i\}] \nonumber\\ \nonumber
& \leq \Pr \left[\left |\widetilde{X}_{i}^2+\widetilde{X}_{i}\sum_{1\leq|k|\leq K}a_k(\widetilde{X}_{i-|k|} + \xi_{i-|k|} ) \exp(-\texttt{i}\omega k) \right |\geq \tilde \tau_n \right ]\\ \nonumber 
%& \leq \Pr \left[\left |\widetilde{X}_{i}\sum_{0\leq|k|\leq K}a_k(\widetilde{X}_{i-|k|} + \xi_{i-|k|} ) \exp(\texttt{i}\omega k) \right |\geq \tilde \tau_n \right ]\\ \nonumber 
&=\Pr \left[|\widetilde{X}_{i}|\cdot \left | \sum_{0\leq|k|\leq K}a_k(\widetilde{X}_{i-|k|} + \xi_{i-|k|} ) \exp(-\texttt{i}\omega k) \right|  \geq \tilde \tau_n \right ]\\ \nonumber 
&\lesssim\Pr \left[ \left | \sum_{0\leq|k|\leq K}  a_k\widetilde{X}_{i-|k|}\exp(-\texttt{i}\omega k) + \sum_{0\leq|k|\leq K}a_k\xi_{i-|k|} \exp(-\texttt{i}\omega k)\right|  \geq \frac{\tilde \tau_n}{\tau_n} \right ]\\
&\lesssim \Pr \left[ \left | \sum_{0 \leq|k|\leq K} a_k\widetilde{X}_{i-|k|} \cos (\omega k) \right | \geq \ \frac 12  \frac{\tilde \tau_n}{\tau_n}\right] + \Pr \left[\left | \sum_{0\leq|k|\leq K} a_k\xi_{i-|k|} \cos (\omega k)\right|  \geq \frac 12   \frac{\tilde \tau_n}{\tau_n} \right] \nonumber  \\ 
&\lesssim \Pr \left[ \left | \sum_{0 \leq k\leq K} a_k\widetilde{X}_{i-k} \cos (\omega k) \right | \geq  \frac{1}{4} \frac{\tilde \tau_n}{\tau_n}  \right] + \Pr \left[\left | \sum_{0\leq k\leq K}a_k \xi_{i-k}\cos (\omega k) \right|  \geq  \frac{1}{4} \ \frac{\tilde \tau_n}{\tau_n}\right] \nonumber \\
&\lesssim \Pr \left[ \left | \sum_{0 \leq k\leq K} {X}_{i-k}\right | \geq  \frac{1}{4} \frac{\tilde \tau_n}{\tau_n} \right] + \Pr \left[\left | \sum_{0\leq k\leq K} \xi_{i-k} \right|  \geq  \frac{1}{4}  \frac{\tilde \tau_n}{\tau_n}\right]. \label{boundthis}
\end{align}
In particular, the first sum consists of mean zero, Gaussian random variables and the second sum of i.i.d. mean zero Laplace random variables.
For the variance of the Gaussian sum we obtain
\begin{align*}
\mathrm{Var}\left( \sum_{0\leq k\leq K}{X}_{i-k}\right)&=\sum_{0\leq k\leq K} \sum_{0\leq l\leq K}\text{Cov} \left({X}_{i-k}, X_{i-l}\right)=\sum_{0\leq k\leq K} \sum_{0\leq l\leq K}\sigma_{|k-l|}\leq M(K+1). %\sigma_0(K+1).
\end{align*} 	
This implies
$$ \Pr \left[ \left | \sum_{0 \leq k\leq K} {X}_{i-k} \right | \geq  \frac{1}{4}  \frac{\tilde \tau_n}{\tau_n}  \right] \leq 2\exp\left (-\frac{\tilde \tau_n^{2}}{32M\tau_n^2(K+1)} \right ),$$
which is of the order $o(n^{-3})$ for $\tilde \tau_n^2 \tau_n^{-2}\geq 96(K+1)\log^{1+\delta}(n)$.
To bound the second term in (\ref{boundthis}) we apply the following version of the Bernstein inequality: Let $Y_1,...,Y_n$ independent zero-mean random variables. If for some $L>0$ and every integer ${\displaystyle k\geq 2}$,
${\displaystyle \mathbb {E} \left[\left|Y_{i}^{k}\right|\right]\leq {\frac {1}{2}}\mathbb {E} \left[Y_{i}^{2}\right]L^{k-2}k!},$
then
$$ {\displaystyle \Pr \left(\sum _{i=1}^{n}Y_{i}\geq 2t{\sqrt {\sum _{i=1}^{n} \mathbb {E} \left[Y_{i}^{2}\right]}}\right)<\exp(-t^{2}),\qquad {\text{for}}\quad 0\leq t\leq {\frac {1}{2L}}{\sqrt {\sum _{i=1}^{n} \mathbb {E} \left[Y_{i}^{2}\right]}}.}$$
In our case $\xi_0,...,\xi_K\overset{\mathrm{i.i.d.}}{\sim}\mathrm{Lap}(b)$ with $b=4\tau_n/\alpha$. Recall that $\mathbb{E}[ |\xi_i^k|] \leq \mathbb{E}[\xi_i^2]k!b^{k-2}/2 $ and thus $L=b$. Then, by symmetry
\begin{align*} 
\Pr \left[\left | \sum_{0\leq k\leq K}\xi_k\right | \geq 8\sqrt{2} t \sqrt{K+1}\tau_n/\alpha  \right] & \leq 2 \exp(-t^2),  \qquad {\text{for}}\quad 0\leq t\leq  \sqrt{(K+1)/2}.
\end{align*} 
If $\tilde \tau_n$ and $\tau_n$ are chosen such that
\begin{equation} \label{condtaun}
0\leq \frac{\alpha \tilde \tau_n }{32\sqrt{2}\sqrt{K+1}\tau_n^2}\leq \sqrt{(K+1)/2},
\end{equation} then
\begin{align}\label{sumLaptail}
\Pr \left[\left | \sum_{0\leq k\leq K}\xi_k\right | \geq  \frac{1}{4}  \frac{\tilde \tau_n}{\tau_n}\right]& \leq 2 \exp\left (-\frac{\tilde \tau_n^2 \alpha^2}{ 2048(K+1)\tau_n^4} \right ). %,  \qquad {\text{for}}\quad 0\leq z\leq  8 \sqrt{(K+1)}/\alpha.
\end{align}
Set $ \tau_n^2=8\log^{1+\delta}(n)$ for some $\delta>0$ and $\tilde{\tau}_n^2=1024\tau_n^{6}(K+1)$. In particular, $\tilde \tau_n^2 \tau_n^{-2}\geq 96(K+1)\log^{1+\delta}(n)$. If $K$ is chosen such that
\begin{equation} \label{condK}
\alpha\tau_n \leq \sqrt{(K+1)},
\end{equation} then condition (\ref{condtaun}) is satisfied for this choice of $\tau_n^2$ and $\tilde \tau_n^2$. Hence, the right-hand side of (\ref{sumLaptail}) is of order $\mathcal{O} \left \{  \exp\left (-4 \log(n)^{1+\delta}\alpha^2\right ) \right \}=\mathcal{O}(n^{-4})$, as for every  $\alpha,\delta>0$ there exists an $n_{\alpha,\delta}$ such that for all $n\geq n_{\alpha,\delta}$ it holds that $\log(n)^\delta\alpha^2 \geq 1.$ Thus, we have shown that for sufficiently large $n$,
\begin{equation} \label{boundeventAcompl}
\mathbb{E}[\mathbbm{I}_{\tilde A^c}\{\hat{f}_K(\omega) - f(\omega)\}^2] \leq \frac{c}{n}. %C\frac{\tilde \tau_n^2}{\alpha^2n}.
\end{equation}
Combining  (\ref{boundeventA}) and (\ref{boundeventAcompl}), setting  $ \tau_n^2=8\log^{1+\delta}(n)$ and $\tilde{\tau}_n^2=1024\tau_n^{6}(K+1)$, and choosing $K$ such that condition (\ref{condK}) is fulfilled gives
\begin{align*}\mathbb{E}|\hat{f}_K(\omega)-f(\omega)|^2&
&=	\mathcal{O}\left (\frac{K}{n} \right) +\mathcal{O}\left(\frac{K\tau_n^{6} }{(n-K)\alpha^ {2}}\right)   + \begin{cases}
\mathcal{O}(K^{-2s}), &f\in W^{s,\infty}(L_0,L),\\
\mathcal{O}(K^{-2s+1}), &f\in  W^{s,2}(L).
\end{cases}
%&=&\mathcal{O}\{\log(K)^2K^{-2s}\}+\mathcal{O}\left\{ \frac{\tau_n^6K}{\alpha^ {2}(n-K)}\right\}
\end{align*}
Optimizing $K^{-2s}+K \cdot( 1\lor \tau_n^{6}\alpha^{-2})/n$ with respect to $K$ gives $$K=\left( \frac{1}{n} \lor \frac{\tau_n^{6}}{n\alpha^2} \right)^{-\frac{1}{2s+1}}.$$
Analogously, optimizing $K^{-2s+1}+K\cdot (1 \lor  \tau_n^{6} \alpha^{-2})/n$ with respect to $K$ gives $$K=\left( \frac{1}{n} \lor \frac{\tau_n^{6}}{n\alpha^2} \right)^{-\frac{1}{2s}}.$$
As both choices of $K$ satisfy condition (\ref{condK}) for sufficiently large $n$, and  $f\in W^{s,\infty}(L_0,L)$ resp. $f\in W^{s,2}(L)$  was chosen arbitrary, we obtain 
\begin{align*}
\sup_{f\in W^{s,\infty}(L_0,L)}\mathbb{E}|\hat{f}_K(\omega)-f(\omega)|^2&\leq c_1 \left [\frac{1}{n} \lor\frac{\log^{3+3\delta}(n)}{n \alpha^2} \right]^{\frac{2s}{2s+1}}
\end{align*} 
and 
\begin{align*}
\sup_{f\in W^{s,2}(L)}\mathbb{E}|\hat{f}_K(\omega)-f(\omega)|^2&\leq c_2 \left [\frac{1}{n} \lor \frac{\log^{3+3\delta}(n)}{n \alpha^2} \right]^{\frac{2s-1}{2s}}.
\end{align*} \hfill \BlackBox

%%%%%%%%%%%%%%%%%%%%%%%%%%%%%%%%%%%%%%%%%%%%%%%%%%%%%%%%%%%%%%%%%%
%SECTION: Interactive privacy mechanism for $f$
%%%%%%%%%%%%%%%%%%%%%%%%%%%%%%%%%%%%%%%%%%%%%%%%%%%%%%%%%%%%%%%%%%

%%%%%%%%%%%%%%%%%%%%%%%%%%%
\subsection{Proof of Lemma~\ref{lemma:EaLDP_f}} 

\emph{Proof of $i)$}: Define $\check Z_i=0$ for $i=1,..,K$ and denote  $\tilde{x}$ for $x$ trimmed at $\pm \tau_n$. Then, for $i=1,...,n$, for all $z$ and $\check z$ and for all $x$ different from $x'$, it holds
\begin{align*}
&\frac{q^{(Z_i,\check Z_i)|X_i=x, (Z_{1},\check Z_{1}),\ldots , (Z_{i-1},\check Z_{i-1})}(z,\check z)}{q^{(Z_i,\check Z_i)|X_i=x', (Z_{1},\check Z_{1}),\ldots , (Z_{i-1},\check Z_{i-1})}(z,\check z)} \\
&= \frac{ \frac{\alpha}{8 \tau_n} \exp\left(-\frac{\alpha}{4 \tau_n } |z- \tilde x| \right)  q^{\check Z_i|X_i=x, Z_{1:(i-1)},\check Z_{1:(i-1)}}(\check z)}{ \frac{\alpha}{8 \tau_n} \exp\left(-\frac{\alpha}{4 \tau_n } |z- \tilde x^\prime| \right) q^{\check Z_i|X_i=x^\prime, Z_{1:(i-1)},\check Z_{1:(i-1)}}(\check z)}\\
&= \frac{ \exp\left(-\frac{\alpha}{4 \tau_n } |z- \tilde x| \right) q^{\check Z_i|X_i=x, W_i}(\check z)}{  \exp\left(-\frac{\alpha}{4 \tau_n } |z- \tilde x^\prime| \right) q^{\check Z_i|X_i=x^\prime, W_i}(\check z)}\leq e^{\alpha/2} \frac{  q^{\check Z_i|X_i=x, W_i}(\check z)}{ q^{\check Z_i|X_i=x^\prime, W_i}(\check z)},
%&\leq \exp\left(\frac{\alpha}{4 \tau_n } |\tilde x- \tilde x^\prime| \right ) e^{\alpha/2} \leq e^\alpha,
\end{align*} where we used $|\tilde x- \tilde x^\prime|\leq 2 \tau_n$. 
The proof that  $\check Z_i$ is an $\alpha/2$-LDP view of $  X_i, W_i$, i.e., that  for all  $\check z$ and all $x\neq x^\prime \in \mathbb{R}^d$ it holds
$$  \frac{  q^{\check Z_i|X_i=x, W_i}(\check z)}{ q^{\check Z_i|X_i=x^\prime, W_i}(\check z)}=\frac{  q^{\check Z_i|X_i=x, \widetilde W_i}(\check z)}{ q^{\check Z_i|X_i=x^\prime, \widetilde W_i}(\check z)}\leq e^{\alpha/2}, $$
follows from results of \cite{duchi2018minimax} and \cite{butucea2023phase}. \\

\emph{Proof of $ii)$}:
By construction, for $i=K+2,...,n$, conditioned on $X_i, Z_{1:(i-1)}$, the  $\check Z_{i}$ is independent from  $\check Z_{l}$ $,\, l=K+1,...,i-1.$  Thus,
\begin{align*}
&\mathbb E[\check Z_{i}\mid X_i, Z_{1:(i-1)}, \check Z_{K+1:(i-1)}] =\mathbb E[\check Z_{i}\mid X_i, Z_{1:(i-1)}].% = \tilde W_i\\
\end{align*}
Next, $$\mathbb E[\check Z_{i}\mid X_i, Z_{1:(i-1)}]=\mathbb E[\check Z_{i}\mid \widetilde W_{i}]= \mathbb E[ \mathbb E[\check Z_{i}\mid \widetilde Y_i]\mid \widetilde W_{i}].$$ By the same arguments as in the proof of Proposition 3.2. of  \cite{butucea2023phase}, one can show that almost surely
$$\mathbb E[\check Z_{i}\mid \widetilde Y_i]= \widetilde Y_i.$$
Furthermore, for every $j=0,..,K$,
$$\mathbb E[ \widetilde Y_{i,j}  \mid \widetilde W_{i,j}=w_{j}]= \tilde \tau_n \left(\frac{1}{2} + \frac{w_{j}}{2\tilde{\tau}_n} \right) - \tilde \tau_n \left(\frac{1}{2} - \frac{w_{j}}{2\tilde{\tau}_n} \right)=w_j.$$
Hence, the mechanism is conditionally unbiased. 
Next,
\begin{align*}
&\mathrm{Var}(  \check Z_{i,j}\mid X_i, Z_{1:(i-1)}, \check Z_{1:(i-1)} )\\
&= \mathbb{E} \left [ (\check Z_{i,j} - \mathbb{E}[ \check Z_{i,j}\mid  X_i, Z_{1:(i-1)}, \check Z_{1:(i-1)} ])^2 \mid  X_i, Z_{1:(i-1)}, \check Z_{1:(i-1)}\right]\\
&= \mathbb{E} \left [ (\check Z_{i,j} - \mathbb{E}[ \check Z_{i,j}\mid  X_i, Z_{1:(i-1)} ])^2 \mid  X_i, Z_{1:(i-1)}\right]{=}\mathrm{Var}(  \check Z_{i,j}\mid X_i, Z_{1:(i-1)} ) \leq \mathrm{Var}(  \check Z_{i,j})\leq  B^2,
\end{align*}
since $\check Z_{i,j}\in\{-B,B\}$ for all $i=K+1,...,n$ and $j=0,...,K$. Finally, it is shown by \citet[Lemma B.1.]{butucea2023phase}  that $C_K$ behaves asymptotically as $\sqrt{\pi/2}\sqrt{K}$, which implies $B=\mathcal{O}(\tilde \tau_n\sqrt{K}\alpha^{-1})$ for $K \to \infty$.
\hfill \BlackBox %\vspace{-.2cm}
\subsection{Proof of Theorem~\ref{theo:ratef}}
We rewrite the estimator $\check{f}_K$ as
\begin{align}
\check f_K(\omega) &= \frac 1{2\pi(n-K)} \sum_{i=K+1}^n  \sum_{0\leq |k|\leq K} \check Z_{i,|k|} \exp(-\texttt{i}k \omega)\nonumber \\
&= \frac 1{2\pi(n-K)} \sum_{0\leq |k|\leq K} \sum_{i=K+1}^n   \check Z_{i,|k|} \exp(-\texttt{i}k \omega) \nonumber\\
&= \frac 1{2\pi} \sum_{0\leq |k|\leq K} \check{\sigma}_k \exp(-\texttt{i}k \omega), \quad \omega\in[-\pi,\pi], \label{reprf}
\end{align}  where $\check Z_{i,|k|}$ is the $|k|$-th entry of $\check Z_{i}$ and $\check{\sigma}_j= \frac{1}{n-K} \sum_{i=K+1}^n   \check Z_{i,|j|}.$	\\
Let $f_K(\omega)=(2\pi)^{-1}\sum_{0\leq |k|\leq K} {\sigma}_k \exp(-\texttt{i}k \omega)$. Then, by Parseval's identity, (\ref{FSapprox}) and (\ref{FSapprox2}),
\begin{align*}
\mathbb{E}\|\check{f}_K-f\|_2^2&=\mathbb{E}\|\check{f}_K-f_K\|_2^2 +\|{f}_K-f\|_2^2\\
&=\sum_{0\leq|j|\leq K} \mathbb{E}[(\check{\sigma}_j-\sigma_j)^2] + \begin{cases}
\mathcal{O}(K^{-2s}), &f\in W^{s,\infty}(L_0,L),\\
\mathcal{O}(K^{-2s}), &f\in  W^{s,2}(L).
\end{cases}
\end{align*}
It is sufficient to consider  $j\in\{0,1,...,K\}$. 
In particular,
\begin{align*}\mathbb{E}[(\check{\sigma}_j - \sigma_j)^2] &= \mathbb{E}[(\check{\sigma}_j - \mathbb{E}[\check{\sigma}_j])^2]  +   (\mathbb{E}[\check{\sigma}_j ] - {\sigma}_j)^2
%&= \mathrm{Var}(\check{\sigma}_j ) +  \mathbb{E}[(\mathbb{E}[\check{\sigma}_j \mid X_{1:n}, Z_{1:n}] - {\sigma}_j)^2 ]  
= \mathrm{Var}(\check{\sigma}_j ) +   (\mathbb{E}[\check{\sigma}_j ] - {\sigma}_j)^2.
\end{align*}
By the law of total variance,
\begin{align*}
\mathrm{Var}(\check{\sigma}_j ) &= \mathbb{E}[ \mathrm{Var}(\check{\sigma}_j\mid \tilde W_{(K+1):n})]+ \mathrm{Var}(\mathbb{E}[\check{\sigma}_j\mid\tilde W_{(K+1):n}]).%\frac{1}{(n-K)^2} \sum_{i=K+1}^n  \sum_{k=K+1}^n  \text{Cov}(\check Z_{i,j},\check Z_{k,j}) \leq \frac{1}{(n-K)^2} \sum_{i=K+1}^n  \sum_{k=K+1}^n  \mathrm{Var}(\check Z_{i,j})\leq B^2
\end{align*} %as  $\check Z_{i,j} \in \{-B,B\}$.\\
Since conditioned on  $\widetilde W_{K+1},..,\tilde W_n$, the  $\check Z_{i,|j|},\, i=K+1,...,n,$ are independent, we have
\begin{align*}
\mathbb{E}[ \mathrm{Var}(\check{\sigma}_j\mid \tilde W_{(K+1):n})]&= \mathbb{E}\left [\frac{1}{(n-K)^2}\sum_{i=K+1}^n \mathrm{Var}(  \check Z_{i,j}\mid \tilde W_{(K+1):n}) \right]\\
&\leq \frac{1}{(n-K)^2}\sum_{i=K+1}^n \mathrm{Var}(  \check Z_{i,j}) \leq \frac{B^2}{n-K},
\end{align*} by Lemma~\ref{lemma:EaLDP_f}. %as  $\check Z_{i,j} \in \{-B,B\}$.
Again, with Lemma~\ref{lemma:EaLDP_f} it follows, that
$$
\mathbb{E}[\check Z_{i}\mid\widetilde W_{(K+1):n}]=\mathbb{E}[\check Z_{i}\mid\widetilde W_{(K+1):i}]=\mathbb{E}[\check Z_{i}\mid X_i, Z_{1:(i-1)}] = \widetilde W_i,
$$ 
and hence
\begin{align} \label{condEsigma}
\mathbb{E}[\check{\sigma}_j\mid \widetilde W_{(K+1):n}]&= \frac{1}{n-K}\sum_{i=K+1}^n  \mathbb{E}[ \check Z_{i,j} \mid  \widetilde W_{(K+1):n} ]= \frac{1}{n-K}\sum_{i=K+1}^n \widetilde W_{i,j}.
\end{align} 
Similar calculations as in the proof of  Theorem~\ref{theo:ratesigma}  yield that
$$\mathrm{Var}\left (\frac{1}{n-K}\sum_{i=K+1}^n \widetilde W_{i,j} \right)\leq \frac{2M_1}{(n-K)} + \frac{8M\tau_n^2}{\alpha^2(n-K)} +\frac{c\tilde{\tau}^2_n}{n-K},$$
if  $ \tau_n^2=8\log^{1+\delta}(n)$ and  $\widetilde \tau_n\geq 16\log^{1+\delta}(n) \tau_n^2$ for some $\delta>0$  is chosen.
Thus,
$$ \mathrm{Var}(\check{\sigma}_j )\leq \frac{2M_1}{(n-K)} + \frac{8M\tau_n^2}{\alpha^2(n-K)} +\frac{c\tilde{\tau}^2_n}{n^2} +\frac{B^2}{n-K}.$$
To compute the bias $(\mathbb{E}[\check{\sigma}_j ] - {\sigma}_j)$, we define, analogous to the proof of  Theorem~\ref{theo:ratesigma}, the event $A_j$ and its complement:
\begin{align*}
A_j&{=} \bigcap_{i{=}1}^n \{ \widetilde{X}_i{=} X_i \text{ and } (\widetilde W_{i,j} {=} W_{i,j})   {I}(i>K) \}, & A_j^c&{=}\bigcup_{i{=}1}^n \{ \widetilde{X}_i{\neq} X_i \text{ or } (\widetilde W_{i,j} {\neq} W_{i,j})  {I}(i>K)\}. 
\end{align*}
Using (\ref{condEsigma}), we get 
\begin{align*}
(\mathbb E [\check{\sigma}_j]-{\sigma}_j)^2&= \left(\mathbb{E} \left[ \frac{1}{n-K}\sum_{i=K+1}^n \widetilde W_{i,j}-{\sigma}_j\right] \right)^2\\
&\lesssim\mathbb{E} \left[ \mathbb{I}_{A_j}\left (\frac{1}{n-K}\sum_{i=K+1}^n \widetilde W_{i,j}-{\sigma}_j \right)^2\right] + \mathbb{E}\left[ \mathbb{I}_{A^c_j} \left (\frac{1}{n-K}\sum_{i=K+1}^n \widetilde W_{i,j}-{\sigma}_j \right)^2\right].
\end{align*}
In particular, if the event $A_j$ occurs, then  
$$\mathbb{E} \left[ \mathbb{I}_{A_j}\left (\frac{1}{n-K}\sum_{i=K+1}^n \widetilde W_{i,j}-{\sigma}_j \right)^2\right]=\mathbb E[\mathbb{I}_{A_j}(\tilde{\sigma}_j-\sigma_j)^2],$$ where $\tilde{\sigma}_j$ is the untruncated estimator defined by
\begin{align*}
\tilde{\sigma}_j&=\frac{1}{n-K} \sum_{i=K+1}^n  X_i \cdot \left ( X_{i-j}+\xi_{i-j}\right).
\end{align*}
Thus, 
%$$ (\mathbb{E}[\hat{\sigma}_j ] - {\sigma}_j)^2  = \mathbb{E}[(\check\sigma_j- {\sigma}_j)^2 ] = \mathbb{E}[\mathbb{I}_{A_j}(\check\sigma_j- {\sigma}_j)^2 ] + \mathbb{E}[[\mathbb{I}_{A_j^c}(\check\sigma_j- {\sigma}_j)^2 ] .$$ 
%$$\mathbb{E} [\check\sigma_j]=\mathbb{E} [\mathbb{I}_{A_j}\check\sigma_j]+\mathbb{E} [\mathbb{I}_{A^c_j}\check\sigma_j].$$
\begin{align*}
(\mathbb{E} [\check\sigma_j]-\sigma_j)^2%&\lesssim \mathbb{E}[\mathbb{I}_{A_j}(\mathbb{E} [\check\sigma_j]-\sigma_j)^2]+\mathbb{E}[\mathbb{I}_{A^c_j}(\mathbb{E} [\check\sigma_j]-\sigma_j)^2]\\
%&\leq \mathbb{E}[\mathbb{I}_{A_j}(\check\sigma_j-\sigma_j)^2]+\mathbb{E}[\mathbb{I}_{A^c_j}(\mathbb{E} [\check\sigma_j]-\sigma_j)^2]\\
&\lesssim\mathbb{E}[\mathbb{I}_{A_j}(\tilde\sigma_j-\sigma_j)^2]+ \left ( \Pr(A^c_j) \mathbb{E}\left[  \left (\frac{1}{n-K}\sum_{i=K+1}^n \widetilde W_{i,j}-{\sigma}_j \right )^4 \right]  \right )^{1/2}.
\end{align*}
Similar calculations as in the proof of Theorem~\ref{theo:ratesigma} yield 
\begin{align*}
\mathbb{E}[\mathbb{I}_{A_j}(\tilde\sigma_j- {\sigma}_j)^2 ] &\leq \frac{2M_1}{(n-K)} + \frac{8M\tau_n^2}{\alpha^2(n-K)}
\end{align*} and that $\Pr(A_j^c)=\mathcal{O}(n^{-3})$, if  $ \tau_n^2=8\log^{1+\delta}(n)$ and $\tilde \tau_n=16\log^{1+\delta}(n) \tau_n^2$ for some $\delta>0$  is chosen.
%As $(\mathbb{E} [\check\sigma_j]-\sigma_j)^2\leq 2B^2+M^2,$
Since $|\widetilde W_{i,j}|\leq \tilde{\tau}_n,$
it follows that
\begin{align*}
\mathbb{E}[(\check{\sigma}_j - \sigma_j)^2]&\leq \frac{2M_1}{n-K} + \frac{8M\tau_n^2}{\alpha^2(n-K)} +\frac{cB^2}{n-K} + o(n^{-1}). %.
\end{align*}
By Stirling's approximation, $B\lesssim \tilde \tau_n \sqrt{K+1}/\alpha$ and thus% see \cite{duchi2018minimax}, yields
\begin{align*}
\mathbb{E}\|\check{f}_K-f\|_2^2&\leq C \left\{ \frac{K}{n-K} +\frac{K\tau_n^2}{\alpha^2(n-K)} +\frac{\tilde \tau_n^2 K^2}{\alpha^2(n-K)} +\frac{1}{K^{2s}} \right \}.
\end{align*}
Optimizing $K^{-2s}+ K^2\tilde \tau_n^{2}/(\alpha^{2}n)+K/n$ with respect to $K$ gives $K\asymp\{\tilde \tau_n^{2}/(n\alpha^2)\}^{-1/(2s+2)}\land n^{-1/(2s+1)}$.
As $f\in W^{s,\infty}(L_0,L)$ resp. $f\in W^{s,2}(L)$  was chosen arbitrary, and $\alpha>0$, we obtain 
\begin{align*}
\sup_{f\in W^{s,\infty}(L_0,L)}\mathbb{E}\|\check{f}_K-f\|^2&\leq c_1  \left \{ \left [ \frac{\log^{4+4\delta}(n)}{n \alpha^2} \right]^{\frac{2s}{2s+2}} \lor \left( \frac{1}{n} \right)^{\frac{2s}{2s+1}}\right \},\\
\sup_{f\in W^{s,2}(L)}\mathbb{E}\|\check{f}_K-f\|^2 &\leq c_2  \left \{ \left[ \frac{\log^{4+4\delta}(n)}{n \alpha^2} \right]^{\frac{2s}{2s+2}} \lor \left( \frac{1}{n} \right)^{\frac{2s}{2s+1}}\right \}.% \lor \frac{\tau_n^4}{n}.
\end{align*} 
\hfill \BlackBox

%%%%%%%%%%%%%%%%%%%%%%%%%%%%%%%%%%%%
%COROLLARY
\subsection{Proof of Corollary~\ref{cor:ratef0}}
Let $\omega\in[-\pi,\pi]$. Then,
\begin{align*}
\mathbb{E}|\check{f}_K(\omega)-f(\omega)|^2&=\{\mathbb{E}[\check{f}_K(\omega)]-f(\omega)\}^2 +\mathbb{E}(\check{f}_K(\omega)-\mathbb{E}[\check{f}_K(\omega)] )^2.
\end{align*}
Using (\ref{reprf}) and (\ref{condEsigma}) yields
\begin{align*}
\{\mathbb{E}[\check{f}_K(\omega)]-f(\omega)\}^2 &= \left \{ \frac 1{2\pi} \sum_{0\leq |j|\leq K} \mathbb{E}[\check{\sigma}_j]\exp(-\texttt{i}j \omega) -  f(\omega) \right  \}^2\\
&= \left \{ \frac 1{2\pi} \sum_{0\leq |j|\leq K} \mathbb{E}[\mathbb{E}[\check{\sigma}_j\mid  \widetilde W_{(K+1):n}]]\exp(-\texttt{i}j \omega) -  f(\omega) \right  \}^2\\
&= \left \{  \mathbb{E} \left [\frac{1}{2\pi(n-K)} \sum_{0\leq |j|\leq K}\sum_{i=K+1}^n \widetilde W_{i,|j|}\exp(-\texttt{i}j \omega) \right ] -  f(\omega) \right  \}^2.
\end{align*}
In particular, if the event $\cup_{j=0}^K A_j$ without any truncation occurs, then the estimator \begin{equation} \label{theest} \frac{1}{2\pi(n-K)} \sum_{0\leq |j|\leq K}\sum_{i=K+1}^n \widetilde W_{i,|j|}\exp(-\texttt{i}j \omega) \end{equation} equals the estimator $\hat{f}^{\text{NI}}_K(\omega)$ in Section~\ref{sec:1}, but with less Laplace noise. %such that the terms with $\alpha^4$ do not appear in the calculation of the variance below. 
The probability of the event that truncation occurs can be computed similarly as in the proof of Theorem~\ref{theo:ratesigma}.
Combining these results, we obtain the following bound for the bias: if  $ \tau_n^2=8\log^{1+\delta}(n)$ and $\tilde \tau_n=16\log^{1+\delta}(n) \tau_n^2$ for some $\delta>0$, then
\begin{align*}
\{\mathbb{E}[\check{f}_K(\omega)]-f(\omega)\}^2 &\lesssim \frac{K}{n} + 
\begin{cases} 
\log(K)K^{-2s}, & f\in W^{s,\infty}(L_0,L),\\
K^{-2s+1}, &f\in W^{s,2}(L),
\end{cases}
\end{align*} where the constant $C>0$ depends on $s,\,L$, if $f\in W^{s,2}(L)$, and additionally on $L_0$, if $ f\in W^{s,\infty}(L_0,L).$
Next, by the law of total variance,
\begin{align*}
\mathrm{Var} \left( \check{f}_K(\omega) \right ) = & \,\mathrm{Var} \left (  \frac 1{2\pi} \sum_{0\leq |k|\leq K} \check{\sigma}_k \exp(-\texttt{i}k \omega) \right )\\
=& \, \mathbb{E} \left [  \mathrm{Var} \left (  \frac 1{2\pi} \sum_{0\leq |k|\leq K} \check{\sigma}_k \exp(-\texttt{i}k \omega) \mid\widetilde W_{(K+1):n} \right)\right ]\\
&\,+  \mathrm{Var} \left (    \mathbb{E} \left [ \frac 1{2\pi} \sum_{0\leq |k|\leq K} \check{\sigma}_k \exp(-\texttt{i}k \omega) \mid\widetilde W_{(K+1):n} \right ] \right).
\end{align*}
As $\mathbb E[ \check \sigma_j\mid  \widetilde W_{(K+1):n}]= \frac{1}{n-K}\sum_{i=K+1}^n \widetilde W_{i,|j|}$, again, the second variance term  above can be bounded by similar calculations as for $\hat{f}^{\text{NI}}_K(\omega)$ in the proof of Proposition~\ref{prop:ratekroll}. In particular,  the terms with $\alpha^4$ do not appear, since (\ref{theest}) includes less added Laplace noise. The probability of the event that truncation occurs can be computed similarly as in the proof of Theorem~\ref{theo:ratesigma}. Thus,
$$\mathrm{Var} \left (    \mathbb{E} \left [ \frac 1{2\pi} \sum_{0\leq |k|\leq K} \check{\sigma}_k \exp(-\texttt{i}k \omega) \mid\widetilde W_{(K+1):n} \right ] \right) \lesssim \frac{K}{n} + \frac{ \tau_n^2 K}{n\alpha^2} .$$
Since conditioned on $W_{(K+1):n},$ the $\check Z_{i,|k|}$ for $i=K+1,...,n$ and $k=0,...,K$ are independent, it follows 
\begin{align*}
& \mathrm{Var} \left (  \frac 1{2\pi} \sum_{0\leq |k|\leq K} \check{\sigma}_k \exp(-\texttt{i}k \omega) \mid\widetilde W_{(K+1):n} \right)\\
&= \frac 1{2\pi} \sum_{0\leq |k|\leq K}  \frac{1}{(n-K)^2} \sum_{i=K+1}^n  \mathrm{Var} \left (   \check Z_{i,|k|} \mid\widetilde W_{(K+1):n} \right)  \lesssim \frac{KB^2}{2\pi(n-K)},
\end{align*}
where we used that $ \check Z_{i,|k|}\in\{-B,B\}$.
%\textcolor{red}{KK: The problem is this term $\mathrm{Var} \left (   \check Z_{i,|k|} \mid\widetilde W_{(K+1):n} \right)$ }.
Combining the previous bounds results in the following bound for the variance
\begin{align*}
\mathrm{Var} \left( \check{f}_K(\omega) \right )=  \mathcal{O}\left (\frac{K}{n} \right) +\mathcal{O}\left(\frac{\tau_n^{2}K }{n\alpha^ {2}}\right) + \mathcal{O}\left( \frac{KB^2}{2\pi(n-K)} \right).
\end{align*}
Using $B\lesssim \tilde \tau_n \sqrt{K+1}/\alpha$  by Lemma~\ref{lemma:EaLDP_f}, we obtain for $ \tau_n^2=8\log^{1+\delta}(n)$ and $\tilde \tau_n=16\log^{1+\delta}(n) \tau_n^2$ for some $\delta>0$, that
\begin{align*}
\mathbb{E}|\check{f}_K(\omega)-f(\omega)|^2= &  \mathcal{O}\left (\frac{K}{n} \right) +\mathcal{O}\left(\frac{K\tau_n^{2} }{n\alpha^ {2}}\right) + \mathcal{O}\left(\frac{K^2\tilde \tau_n^2}{2\pi(n-K)\alpha^2} \right)\\
&+\begin{cases}
\mathcal{O}(\log(K)K^{-2s}), &f\in W^{s,\infty}(L_0,L),\\
\mathcal{O}(K^{-2s+1}), &f\in  W^{s,2}(L).
\end{cases}
\end{align*}
Optimizing with respect to $K$ completes the proof. \hfill \BlackBox
%%%%%%%%%%%%%%%%%%%%%%%%%%%%%%%%%%%%%%%%%%%%%%%%%%%%%%%%%%%%%%%%%%
%SECTION: Theoretical guarantees
%%%%%%%%%%%%%%%%%%%%%%%%%%%%%%%%%%%%%%%%%%%%%%%%%%%%%%%%%%%%%%%%%%

\subsection{Proof of Lemma~\ref{FisherInfo}}
Now, we prove Equation~\eqref{term1}, that is, for all $i=1,\ldots,n$, 
\begin{align*}
I(Z_i|Z_{1:(i-1)})  \lesssim (e^\alpha -1)^2 I(X_i|Z_{1:(i-1)}) ,
\end{align*}
where, by convention, there is no conditioning when $i=1$. Note that
\begin{align} \label{denom}
\mathcal{L}_\theta^{Z_i|Z_{1:(i-1)}}(z_i) = \int q^{Z_i|X_i, Z_{1:(i-1)}} (z_i)\mathcal{L}_\theta^{X_i|Z_{1:(i-1)}}\, \mathrm{d}x_i \geq \inf_x q^{Z_i|X_i=x, Z_{1:(i-1)}} (z_i).
\end{align}
Moreover, $\int \dot{\mathcal{L}}_\theta^{X_i|Z_{1:(i-1)}}(z) \, \mathrm{d}z = 0$ implies that 
\begin{align*}
\int \left[\dot{\mathcal{L}}_\theta^{X_i|Z_{1:(i-1)}}(x_i)\right]_+ \, \mathrm{d}x_i &=  \int \left[\dot{\mathcal{L}}_\theta^{X_i|Z_{1:(i-1)}}(x_i)\right]_- \, \mathrm{d}x_i \\
&= \frac 12 \int \left|\dot{\mathcal{L}}_\theta^{X_i|Z_{1:(i-1)}}(x_i)\right| \, \mathrm{d}x_i \leq \frac12 I^{1/2}(X_i|Z_{1:(i-1)})  .  
\end{align*}
Thus, see also Duchi et al.,
\begin{align}\label{num}
\dot{\mathcal{L}}_\theta^{Z_i|Z_{1:(i-1)}}(z_i) & = 
\int q^{Z_i|X_i, Z_{1:(i-1)}} (z_i) \dot{\mathcal{L}}_\theta^{X_i|Z_{1:(i-1)}}(x_i) \, \mathrm{d}x_i \nonumber \\
& = 
\int q^{Z_i|X_i, Z_{1:(i-1)}} (z_i) \left( \left[\dot{\mathcal{L}}_\theta^{X_i|Z_{1:(i-1)}}(x_i)\right]_+
- \left[\dot{\mathcal{L}}_\theta^{X_i|Z_{1:(i-1)}}(x_i)\right]_- \right) \, \mathrm{d}x_i \nonumber 
\\
& \leq | \sup_x q^{Z_i|X_i=x, Z_{1:(i-1)}} (z_i) - \inf_x q^{Z_i|X_i=x, Z_{1:(i-1)}} (z_i)| \frac 12 \int \left|\dot{\mathcal{L}}_\theta^{X_i|Z_{1:(i-1)}}(x_i)\right| \, \mathrm{d}x_i \nonumber \\
&\leq (e^\alpha - 1) \cdot (2 \wedge e^\alpha) \cdot  \inf_x q^{Z_i|X_i=x, Z_{1:(i-1)}} (z_i) \cdot \frac 12 \int \left|\dot{\mathcal{L}}_\theta^{X_i|Z_{1:(i-1)}}(x_i)\right| \, \mathrm{d}x_i .
\end{align}
We conclude that \eqref{term1} holds by bounding $(2 \wedge e^\alpha)\cdot \frac 12 \leq 1$ and by using \eqref{denom} and \eqref{num} that
\begin{align*}
I(Z_i|Z_{1:(i-1)}) 
& \leq (e^\alpha - 1)^2  \mathbb E_{Z_{1:(i-1)}} \left[ \int \inf_x q^{Z_i|X_i=x, Z_{1:(i-1)}} (z_i) \, \mathrm{d}z_i \left( \int \left|\dot{\mathcal{L}}_\theta^{X_i|Z_{1:(i-1)}}(x_i)\right|  \mathrm{d}x_i\right)^2
\right]\\
& \leq (e^\alpha - 1)^2  I(X_i|Z_{1:(i-1)}).
\end{align*}

For $i\geq 2$, we prove \eqref{term2}, that is, for an arbitrary $\epsilon>0$,
\begin{align*} 
I(X_i|Z_{1:(i-1)})
&\leq (1+ \frac 1\epsilon) (e^\alpha - 1)^2 I(X_i,X_{i-1}|Z_{1:(i-2)}) + (1+\epsilon) e^{2\alpha} I(X_i|Z_{1:(i-2)}) .
\end{align*}
Thus,
\begin{align*}
I(X_i|Z_{1:(i-1)}) & = I(X_i,Z_{1:(i-1)}) {- }I(Z_{1:(i-1)}) = I(X_i,Z_{i-1}|Z_{1:(i-2)}) {+} I(Z_{1:(i-2)}) {- }I(Z_{1:(i-1)})\\
&\leq \mathbb E_{Z_{1:(i-2)}}  \left[ \int \frac{(\dot{\mathcal{L}}_\theta^{X_i,Z_{i-1}|Z_{1:(i-2)}} (x_i,z_{i-1}))^2}{\mathcal{L}_\theta^{X_i,Z_{i-1}|Z_{1:(i-2)}}(x_i, z_{i-1}) } \, \mathrm{d} x_i \, \mathrm{d}z_{i-1} \right] - I(Z_{1:(i-1)}|Z_{1:(i-2)}).
\end{align*}
We use now that, conditionally on $X_{i-1}$, the variables $X_i$ and $Z_{i-1}$ are independent to write that
\begin{align*} % \label{denom2}
\mathcal{L}_\theta^{X_i,Z_{i-1}|Z_{1:(i-2)}}(x_i, z_{i-1}) 
&= \int q^{Z_{i-1}| X_{i-1}, Z_{1:(i-2)}}(z_{i-1}) \cdot \mathcal{L}_\theta^{X_i,X_{i-1}|Z_{1:(i-2)}}(x_i, x_{i-1})  \, \mathrm{d}x_{i-1} \nonumber\\
& \geq \inf_x q^{Z_{i-1}| X_{i-1}=x, Z_{1:(i-2)}}(z_{i-1}) \cdot \mathcal{L}_\theta^{X_i|Z_{1:(i-2)}}(x_i).
\end{align*}
Moreover, 
\begin{align*}
& \dot{\mathcal{L}}_\theta^{X_i,Z_{i-1}|Z_{1:(i-2)}} (x_i,z_{i-1})    \\&= \int q^{Z_{i-1}| X_{i-1}, Z_{1:(i-2)}}(z_{i-1}) \left( \left[\dot{\mathcal{L}}_\theta^{X_i,X_{i-1}|Z_{1:(i-2)}}(x_i, x_{i-1})  \right]_+  {-} \left[\dot{\mathcal{L}}_\theta^{X_i,X_{i{-}1}|Z_{1:(i-2)}}(x_i, x_{i-1})  \right]_-  \right) \mathrm{d}x_{i-1} \\
& \leq  \sup_x q^{Z_{i-1}| X_{i-1}=x, Z_{1:(i-2)}}(z_{i-1}) \cdot 
\int \left[\dot{\mathcal{L}}_\theta^{X_i,X_{i-1}|Z_{1:(i-2)}}(x_i, x_{i-1})  \right]_+  \mathrm{d}x_{i-1} \\
&- \inf_x q^{Z_{i-1}| X_{i-1}=x, Z_{1:(i-2)}}(z_{i-1})\cdot  \int \left[\dot{\mathcal{L}}_\theta^{X_i,X_{i-1}|Z_{1:(i-2)}}(x_i, x_{i-1})  \right]_-  \mathrm{d}x_{i-1}  . 
\end{align*}
Since $\int \dot{\mathcal{L}}_\theta^{X_i,X_{i-1}|Z_{1:(i-2)}}(x_i, x_{i-1})  \, \mathrm{d}x_{i-1}  = \dot{\mathcal{L}}_\theta^{X_i|Z_{1:(i-2)}}(x_i )  $, we get
\begin{align*}
\int \left[\dot{\mathcal{L}}_\theta^{X_i,X_{i-1}|Z_{1:(i-2)}}(x_i, x_{i-1})  \right]_+  \mathrm{d}x_{i-1} 
&{=} \frac 12 \int \left|\dot{\mathcal{L}}_\theta^{X_i,X_{i-1}|Z_{1:(i-2)}}(x_i, x_{i-1})  \right|  \mathrm{d}x_{i-1}  {+} \frac 12 \dot{\mathcal{L}}_\theta^{X_i|Z_{1:(i-2)}}(x_i ) \\
\int \left[\dot{\mathcal{L}}_\theta^{X_i,X_{i-1}|Z_{1:(i-2)}}(x_i, x_{i-1})  \right]_-  \mathrm{d}x_{i-1} 
&{=} \frac 12 \int \left|\dot{\mathcal{L}}_\theta^{X_i,X_{i-1}|Z_{1:(i-2)}}(x_i, x_{i-1})  \right|  \mathrm{d}x_{i-1} { - }\frac 12 \dot{\mathcal{L}}_\theta^{X_i|Z_{1:(i-2)}}(x_i ).
\end{align*}
Therefore,
\begin{align}
& \dot{\mathcal{L}}_\theta^{X_i,Z_{i-1}|Z_{1:(i-2)}} (x_i,z_{i-1})  \nonumber  \\
& \leq (e^\alpha-1) \cdot (2 \wedge e^\alpha) \cdot \inf_x q^{Z_{i-1}| X_{i-1}=x, Z_{1:(i-2)}}(z_{i-1})\cdot \frac 12  \int \left|\dot{\mathcal{L}}_\theta^{X_i,X_{i-1}|Z_{1:(i-2)}}(x_i, x_{i-1})  \right| \, \mathrm{d}x_{i-1} \nonumber  \\
&+ \frac 12 \left(\sup_x q^{Z_{i-1}| X_{i-1}=x, Z_{1:(i-2)}}(z_{i-1}) + \inf_x q^{Z_{i-1}| X_{i-1}=x, Z_{1:(i-2)}}(z_{i-1}) \right) \left| \dot{\mathcal{L}}_\theta^{X_i|Z_{1:(i-2)}}(x_i ) \right|  \nonumber \\
& \leq (e^\alpha-1) \cdot \inf_x q^{Z_{i-1}| X_{i-1}=x, Z_{1:(i-2)}}(z_{i-1})\cdot  \int \left|\dot{\mathcal{L}}_\theta^{X_i,X_{i-1}|Z_{1:(i-2)}}(x_i, x_{i-1})  \right| \, \mathrm{d}x_{i-1} \nonumber  \\
&+ \frac 12 (e^\alpha +1) \inf_x q^{Z_{i-1}| X_{i-1}=x, Z_{1:(i-2)}}(z_{i-1}) \left| \dot{\mathcal{L}}_\theta^{X_i|Z_{1:(i-2)}}(x_i )\right| . \label{term31}
\end{align}
Thus, using that $(a+b)^2 \leq (1+\frac 1\epsilon) a^2 + (1+\epsilon) b^2$ for any $\epsilon >0$, we get that
\begin{align*}
& \mathbb E_{Z_{1:(i-2)}}  \left[ \int \frac{(\dot{\mathcal{L}}_\theta^{X_i,Z_{i-1}|Z_{1:(i-2)}} (x_i,z_{i-1}))^2}{\mathcal{L}_\theta^{X_i,Z_{i-1}|Z_{1:(i-2)}}(x_i, z_{i-1}) } \, \mathrm{d}x_i \, \mathrm{d}z_{i-1} \right] \\
& \leq \mathbb E_{Z_{1:(i-2)}}  \left[ 
\int \inf_x q^{Z_{i-1}| X_{i-1}=x, Z_{1:(i-2)}}(z_{i-1}) \, \mathrm{d} z_{i-1} \right.\\
& \left.\left( (1+\frac 1\epsilon) (e^\alpha - 1)^2 \int \frac{(\int \left|\dot{\mathcal{L}}_\theta^{X_i,X_{i-1}|Z_{1:(i-2)}}(x_i, x_{i-1})  \right| dx_{i-1})^2}{\mathcal{L}_\theta^{X_i|Z_{1:(i-2)}}(x_i)}  \, \mathrm{d} x_{i}  \right. \right.\\
& \left.\left. + (1+\varepsilon) \frac 14 (e^\alpha+1)^2 \int \frac { \left(\dot{\mathcal{L}}_\theta^{X_i|Z_{1:(i-2)}}(x_i )\right)^2}{{\mathcal{L}_\theta^{X_i|Z_{1:(i-2)}}(x_i)} } \, \mathrm{d} x_i
\right) \right].
\end{align*}
Note first that $\mathbb E_{Z_{1:(i-2)}}  \left[ \int \frac { \left(\dot{\mathcal{L}}_\theta^{X_i|Z_{1:(i-2)}}(x_i )\right)^2}{{\mathcal{L}_\theta^{X_i|Z_{1:(i-2)}}(x_i)} } \, \mathrm{d} x_i \right]$ is the Fisher information $I(X_i|Z_{1:(i-2)})$ and then 
\begin{align*}
&\mathbb E_{Z_{1:(i-2)}}  \left[ 
\int \frac{(\int \left|\dot{\mathcal{L}}_\theta^{X_i,X_{i-1}|Z_{1:(i-2)}}(x_i, x_{i-1})  \right| \, \mathrm{d} x_{i-1})^2}{\mathcal{L}_\theta^{X_i|Z_{1:(i-2)}}(x_i)} \, \mathrm{d} x_{i} \right]\\
&\leq \mathbb E_{Z_{1:(i-2)}}  \left[ \int \int \frac{\left|\dot{\mathcal{L}}_\theta^{X_i,X_{i-1}|Z_{1:(i-2)}}(x_i, x_{i-1})  \right|^2}{\mathcal{L}_\theta^{X_i,X_{i-1}|Z_{1:(i-2)}}(x_i)} \, \mathrm{d} x_{i-1} \, \mathrm{d}x_{i}
\right]
\end{align*}
which is the Fisher information $I(X_i,X_{i-1}| Z_{1:(i-2)})$ and this finishes the proof of \eqref{term2}.

For $i\geq 2$, we prove \eqref{term3}, that is
\begin{align*} 
I(X_i|Z_{1:(i-1)})  \leq (1 + \frac 94 (e^\alpha-1)^2) \left( I(X_i,X_{i-1}|Z_{1:(i-2)}) 
+  I(X_i|Z_{1:(i-2)})
\right).
\end{align*}    
Let us go back to \eqref{term31}. We decompose $\frac 12(e^\alpha +1) = \frac 12 (e^\alpha -1) +1$ and use that:
$$
\left| \dot{\mathcal{L}}_\theta^{X_i|Z_{1:(i-2)}}(x_i )  \right| 
= \left| \int \dot{\mathcal{L}}_\theta^{X_i, X_{i-1}|Z_{1:(i-2)}}(x_i ,x_{i-1} ) \, \mathrm{d} x_{i-1} \right| 
\leq \int \left| \dot{\mathcal{L}}_\theta^{X_i, X_{i-1}|Z_{1:(i-2)}}(x_i , x_{i-1})  \right| \, \mathrm{d} x_{i-1}.
$$
Therefore,
\begin{align*}
& \dot{\mathcal{L}}_\theta^{X_i,Z_{i-1}|Z_{1:(i-2)}} (x_i,z_{i-1})  \nonumber  \\
% & \leq (e^\alpha-1) \cdot (2 \wedge e^\alpha) \cdot \inf_x q^{Z_{i-1}| X_{i-1}=x, Z_{1:(i-2)}}(z_{i-1})\cdot \frac 12  \int \left|\dot{\mathcal{L}}_\theta^{X_i,X_{i-1}|Z_{1:(i-2)}}(x_i, x_{i-1})  \right| dx_{i-1} \nonumber  \\
%&+ \frac 12 \left(\sup_x q^{Z_{i-1}| X_{i-1}=x, Z_{1:(i-2)}}(z_{i-1}) + \inf_x q^{Z_{i-1}| X_{i-1}=x, Z_{1:(i-2)}}(z_{i-1}) \right) \left| \dot{\mathcal{L}}_\theta^{X_i|Z_{1:(i-2)}}(x_i ) \right|  \nonumber \\
& \leq \frac 32 (e^\alpha-1) \cdot \inf_x q^{Z_{i-1}| X_{i-1}=x, Z_{1:(i-2)}}(z_{i-1})\cdot  \int \left|\dot{\mathcal{L}}_\theta^{X_i,X_{i-1}|Z_{1:(i-2)}}(x_i, x_{i-1})  \right| \, \mathrm{d} x_{i-1} \nonumber  \\
&+ \inf_x q^{Z_{i-1}| X_{i-1}=x, Z_{1:(i-2)}}(z_{i-1}) \left| \dot{\mathcal{L}}_\theta^{X_i|Z_{1:(i-2)}}(x_i )\right| . 
\end{align*}
Similarly to the previous proof we obtain for an arbitrary $\epsilon >0$ that
\begin{align*}
I(X_i|Z_{1:(i-1)})  \leq (1 + \frac 1\epsilon) \frac 94 (e^\alpha-1)^2  I(X_i,X_{i-1}|Z_{1:(i-2)}) 
+ (1+\epsilon) I(X_i|Z_{1:(i-2)})
\end{align*}
and conclude by choosing $\epsilon = \frac 94 (e^\alpha-1)^2$. \hfill \BlackBox

\subsection{Proof of Theorem~\ref{thm:LBcoeff}}
Let us consider a centered Gaussian time series. Define the parametric family of covariance matrices $\Sigma(\theta)$, where $\theta$ belongs to $[-T, T]$ for some given $T>0$, such that, for some constant $\sigma>0$
$$
\sigma_0 = 1; \, \sigma_j = 0, \text{ for } j\not \in \{0,K\} \quad \text{ and } \quad \sigma_K(\theta) = \theta  \cdot \sigma \cos(K w_0), \, \, K\geq 1 .
$$
We assume that we find $w_0$ in $(-\pi,\pi)$ such that $\cos(K w_0) \not = 0$. These stationary Gaussian processes are associated to the parametric family of spectral densities
\begin{align*}
f_\theta(w) &= \frac{\sigma_0}{2\pi} +  \frac{\sigma_K( \theta )}{\pi} \cdot \cos(K w) , \\
&= \frac{\sigma_0}{2\pi} +  \theta \frac{ \sigma \cos(K w_0) }{\pi}\cdot \cos(K w) \\
&= \frac{\sigma_0}{2\pi} +  \theta \frac{\sigma}{2 \pi} \cdot \left[\cos(K(w-w_0))+\cos(K(w+w_0)) \right].
\end{align*}
Indeed, $f_\theta$ is still a symmetric function and bounded from below by some positive constant
$$
\frac{\sigma_0 }{2\pi} -  T \frac \sigma \pi \geq c >0,
$$
for any $\theta$ in the set $[-T,T]$ and properly chosen $T$.\\
Let us see that estimating $\sigma_K(\theta)$ can be reduced to estimating $f_\theta(w_0)$ in our model, since 
$$
\sigma_K(\theta) = \left( f_\theta(w_0) - \frac{\sigma_0}{2 \pi} \right) \cdot \frac{\pi}{\cos(Kw_0)}
$$ and $\sigma_0$ is known. Thus, we bound the minimax risk from below as follows:
\begin{align*}
R:=\inf_{\hat \sigma_K} \sup_{Q \in \mathcal{Q}^{\mathrm{NI}}_\alpha} \sup_{\sigma:\|\sigma\|_2^2\leq M_1} \mathbb E_\theta \left[ |\hat \sigma_K - \sigma_K|^2\right] 
& \geq \inf_{\hat \sigma_K} \sup_{Q \in \mathcal{Q}^{\mathrm{NI}}_\alpha} \sup_{\theta \in [-T,T]} \mathbb E_\theta \left[ |\hat \sigma_K - \sigma_K(\theta)|^2\right]\\
&\gtrsim \inf_{\hat f_n} \sup_{Q \in \mathcal{Q}^{\mathrm{NI}}_\alpha} \sup_{\theta \in [-T,T]}  \mathbb E_\theta \left[|\hat f_n(w_0) - f_\theta (w_0)|^2\right]
\end{align*}
by first reducing the risk over all possible covariance matrices $\Sigma$ (uniformly bounded in spectral norm) to the parametric family $\Sigma(\theta)$. Then we changed the problem of recovering $\sigma_K(\theta)$ into that of recovering $f_\theta$. Next, we reduce the problem to a Bayesian problem:
\begin{align*}
R    &\gtrsim \inf_{\hat f_n} \sup_{Q \in \mathcal{Q}^{\mathrm{NI}}_\alpha} \int_{[-T,T]} \mathbb E_\theta \left[|\hat f_n(w_0) - f_\theta (w_0)|^2\right] \, \mathrm{d} \lambda_{T} (\theta),
\end{align*}
where $\lambda$ is an a priori measure on [-1,1] with finite Fisher information  
$$
I(\lambda) := \int_{[-1,1]} (\lambda'(u))^2/\lambda(u) \, \mathrm{d}u
$$ and $\lambda_{T}$ is $\lambda$ re-scaled to the interval $[-T,T]$. Thus, $I(\lambda_{T}) = I(\lambda)/T^2$.\\
We apply the van Trees inequality to get
\begin{align*}
\int_{[-T,T]} \mathbb E_\theta \left[|\hat f_n(w_0) - f_\theta (w_0)|^2\right] \, \mathrm{d}\lambda_{T} (\theta)
& \geq \frac{(\int_{[-T,T]} \frac{\partial}{\partial \theta} f_\theta(w_0) \, \mathrm{d}\lambda_{T} (\theta))^2}{\int_{[-T,T]}I(Z_{1:n}) \, \mathrm{d}\lambda_{T} (\theta) + I(\lambda_{T})} \\
&\geq \frac{(\int_{[-T,T]} \sigma \cos^2(K w_0)/{\pi} \, \mathrm{d}\lambda_{T} (\theta))^2}{\int_{[-T,T]}I(Z_{1:n}) \, \mathrm{d}\lambda_{T} (\theta) + I(\lambda)/T^2} \\
& \geq \frac{ \sigma^2 \cos^4(Kw_0) /\pi^2}{ \int_{[-T,T]}I(Z_{1:n}) \, \mathrm{d}\lambda_{T} (\theta) + I(\lambda)/T^2} . 
\end{align*}
Let $(Z_1,...,Z_n)$ be the time series obtained by a non-interactive $\alpha-$differentially private mechanism, $Q$ in $\mathcal{Q}^{\mathrm{NI}}_\alpha$.
The corresponding likelihoods of the vectors $X_{1:n}$ and of $Z_{1:n}$ are denoted by $\mathcal{L}_\theta^{X_{1:n}}$ and $\mathcal{L}_\theta^{Z_{1:n}}$, respectively.

%\bigskip

We want to give an upper bound for $I(Z_{1:n})$ in order to finish the proof. Let us note that by our choice of the covariance matrix $\Sigma(\theta)$ the likelihood $\mathcal{L}_\theta^{X_{1:n}}$ can be decomposed as a product of $K$ likelihoods :
$$
\mathcal{L}_\theta^{X_{n},X_{n-K},\ldots , X_{n-N_0 K}} \cdot \mathcal{L}_\theta^{X_{n-1},X_{n-1-K},\ldots , X_{n-1-N_1 K}} \cdot \ldots \cdot \mathcal{L}_\theta^{X_{n-K+1}, X_{n-2K+1},\ldots, X_{n-N_{K-1}K+1}},
$$
where $N_i = \lfloor (n-i-K)/K\rfloor$ for $i$ from 0 to $K-1$. Each vector is 1-dependent. The non-interactive privacy mechanism preserves this structure, so the same holds for $\mathcal{L}_\theta^{Z_{1:n}}$. Since the Fisher information of jointly independent vectors is additive, we write :
\begin{equation}
\label{term4.0}
I(Z_{1:n}) = \sum_{i=0}^{K-1} I(Z_{n-i}, Z_{n-i-K},\ldots, Z_{n-i-N_iK}).    
\end{equation}
We shall first give the proof for the case $K=1$ when the whole vector $X_{1:n}$ is 1-dependent. Next, we will adapt the calculations to the case $K>1$ using the decomposition above. 

Let us suppose $K=1$. First we apply \eqref{term1} and next \eqref{term2} in Lemma~\ref{FisherInfo} to get:
\begin{align}
I(Z_{1:n}) & = \sum_{i=1}^n I(Z_i | Z_{1:(i-1)}) \nonumber \\
& \lesssim (e^\alpha - 1)^2 \sum_{i=1}^n I(X_i | Z_{1:(i-1)}) \label{term4} \\
& \lesssim (e^\alpha - 1)^2 \sum_{i=1}^n (A \cdot I(X_i, X_{i-1} | Z_{1:(i-2)}) + B \cdot I(X_i | Z_{1:(i-2)}) ) , \nonumber
\end{align}
with $A=(1+ \frac 1\epsilon)(e^\alpha-1)^2$ and $B = (1+\epsilon) e^{2\alpha}$ for an arbitrary $\epsilon>0$ such that $A <1$. Since $X_i$ is independent of $X_{1:(i-2)}$, then we can easily show that $X_i$ is independent of $Z_{1:(i-2)}$. Thus $I(X_i | Z_{1:(i-2)}) ) = I(X_i)$. By iterating and repeating this reasoning we get that:
\begin{equation}\label{term5}
I(X_i | Z_{1:(i-1)}) \leq A^{i-1} I(X_{1:i}) + A^{i-2}B I(X_{2:i}) + \ldots + ABI(X_{i-1},X_{i}) + B I(X_i).    
\end{equation}
Note that $X_i$ has a likelihood free of the parameter $\theta$ so it does not bring any information on $\theta$, that is $I(X_i)=0$. 

Let us prove now that $I(X_{i-\ell},\ldots,X_i) \leq C \ell \sigma^2 $, for any $\ell$ from 1 to $i-1$, where $C= (4\pi^2 c^2)^{-1}$ depends on $c$ introduced in the definition of the spectral density $f_\theta$. The vector $X_{(i-\ell):i}$ has a centered Gaussian distribution with covariance matrix $\Sigma_1^{\ell+1}(\theta)$ given by:
$$
\Sigma_1^{\ell+1}(\theta) = \sigma_0 \cdot \mathbb{I}_{\ell+1} + \sigma_1(\theta) \cdot M, \quad [M]_{i,j}  = \mathbf{1}(|i-j|=1),
$$
where $\mathbb{I}_{\ell+1}$ denotes the identity matrix of size $\ell+1$.  Let us denote by $u_0,\ldots,u_{\ell}$ the eigenvalues of the matrix $M$. Note that $M$ and $ \Sigma_1^{\ell+1}(\theta) $ share the same eigenvectors.
We can write:
\begin{align*}
I(X_{i-\ell},\ldots,X_i) &= \frac 12 tr\left[ \left( \Sigma_1^{\ell+1}(\theta)^{-1} \frac{\partial}{\partial \theta} \Sigma_1^{\ell+1}(\theta) \right)^2\right]\\
&= \frac 12 \left\| \Sigma_1^{\ell+1}(\theta)^{-1} \cdot \sigma \cos(w_0) M \right\|_F^2\\
&= \frac{1}{2 } \sum_{k=0}^\ell \frac{\sigma^2 \cos^2(w_0) u_k^2}{(\sigma_0+\sigma_1(\theta) \cdot u_k)^2}\\
&\leq \frac{\sigma^2 \cos^2(w_0)  \|M\|_F^2}{2 (\sigma_0 - T \sigma \|M\|_\infty)^2}.
\end{align*}
It is easy to calculate that $\|M\|_F^2 = 2 \ell$ and that its operator norm is $\|M\|_\infty = 2 \cos(\pi/(\ell+2)) \leq 2$. By our construction of the spectral density we get that:
$$
I(X_{i-\ell},\ldots,X_i) \leq \frac{\sigma^2 \cdot 2 \ell}{2 (2\pi c)^2} \leq \frac{\sigma^2}{4\pi^2 c^2} \ell = C \ell \sigma^2.
$$
Use this in \eqref{term5} to get:
$$
I(X_i | Z_{1:(i-1)}) \leq AB \sum_{\ell = 1}^{i-1} C \ell \sigma^2 A^{\ell-1} \lesssim (e^\alpha -1)^2 \frac{BC \sigma^2}{(1-A)^2}.
$$
Finally, we use this in \eqref{term4} to get:
$$
I(Z_{1:n}) \lesssim n (e^\alpha - 1)^4  \sigma^2,
$$
up to constants depending on fixed positive constants $\epsilon$ and $c$. This finishes the proof for the case $K=1$.

Now, for $K>1$ and each $i$ in $0,\ldots,K-1$, let us bound similarly each term in \eqref{term4.0}, of the form  $I(Z_{n-i},Z_{n-i-K},...,Z_{n-i-N_i K})$. We get
\begin{eqnarray*}
&& I(Z_{n-i},Z_{n-i-K},...,Z_{n-i-N_i K})\\ 
&=& 
I(Z_{n-i}|Z_{n-i-K},...,Z_{n-i-N_i K})
{+}I(Z_{n-i-K}|Z_{n-i-2K},...,Z_{n-i-N_i K})    {+}\ldots {+} I(Z_{n-i-N_i K})\\
&\lesssim &(e^\alpha - 1)^2 (I(X_{n-i}|Z_{n-i-K},...,Z_{n-i-N_i K})
+I(X_{n-i-K}|Z_{n-i-2K},...,Z_{n-i-N_i K})\\
&  &  +\ldots + I(X_{n-i-N_i K})).
\end{eqnarray*}
We iterate the procedure by applying \eqref{term2} in Lemma~\ref{FisherInfo} in order to get geometric series analogous to \eqref{term5}. Recall that the subvector $X_{n-i},X_{n-i-K},...,X_{n-i-\ell K}$ has a centered Gaussian distribution with covariance matrix 
$$
\Sigma_K^{\ell +1} (\theta ) = \sigma_0 \cdot \mathbb{I}_{\ell+1} + \sigma_K(\theta) \cdot M, \quad \ell = 1,...,N_i.
$$
Following the same steps as for $K=1$ we get:
\begin{align*}
I(X_{n-i-\ell K},\ldots,X_{n-i}) &{=} \frac 12 tr\left[ \left( \Sigma_K^{\ell+1}(\theta)^{-1} \frac{\partial}{\partial \theta} \Sigma_K^{\ell+1}(\theta) \right)^2\right]{= }\frac 12 \left\| \Sigma_K^{\ell+1}(\theta)^{-1} \cdot \sigma \cos(K w_0) M \right\|_F^2\\
&= \frac{1}{2 } \sum_{k=0}^\ell \frac{\sigma^2 \cos^2(K w_0) u_k^2}{(\sigma_0+\sigma_K(\theta) \cdot u_k)^2}\\
&\leq \frac{\sigma^2 \cos^2(K w_0)  \|M\|_F^2}{2 (\sigma_0 - T \sigma \|M\|_\infty)^2} \leq \frac{\sigma^2 \ell}{(\sigma_0 - 2 T \sigma)^2} \leq \frac{\sigma^2 \ell}{4 \pi^2 c^2},
\end{align*}
where $c>0$ is chosen in the definition of the spectral density $f_\theta$. Further we get that:
$$
I(X_{n-i}|Z_{n-i-K},...,Z_{n-i-N_i K}) \lesssim  (e^\alpha - 1)^2 \frac{BC \sigma^2}{(1-A)^2}
$$
and, finally, since $N_i \leq n/K$:
\begin{eqnarray*}
I(Z_{1:n}) & \lesssim & \sum_{i=0}^{K-1} \frac nK (e^\alpha - 1)^4 \sigma^2 = n (e^\alpha - 1)^4 \sigma^2,
\end{eqnarray*}
up to constants depending on fixed positive constants $\epsilon$ and $c$.

\bigskip

%%%%%%%%%%%%%%%%%%%%%%%
Finally, 
$$
I(Z_{1:n}) \lesssim (e^\alpha-1)^2 \sum_{i=1}^n I(X_i|Z_{1:(i-1)}) \lesssim n (e^\alpha-1)^4 \sigma^2.
$$
We conclude by recalling that $w_0$ was chosen such that  $\cos(K w_0) \not= 0$ and by plugging this into the Van Trees inequality above, to get that, for large enough $n$,
\begin{align*}
\inf_{\hat \sigma_K} \sup_{Q \in \mathcal{Q}^{\mathrm{NI}}_\alpha} \sup_{\sigma:\|\sigma\|_2^2\leq M_1} \mathbb E_\theta \left[ |\hat \sigma_K - \sigma_K|^2\right]
& \gtrsim \frac{\sigma^2 \cos^4(Kw_0) /\pi^2}{ n (e^\alpha-1)^4 \sigma^2 + I(\lambda)/T^2} \asymp \frac{1}{n \alpha^4} . 
\end{align*}
\hfill \BlackBox

%%%%%%%%%%%%%%%%%%%%%%%%%%%%%%%%%%%%%%%%%%%%%%%%%%%%%%%%%%%%%%%%%
%APPENDIX: Simulation
%%%%%%%%%%%%%%%%%%%%%%%%%%%%%%%%%%%%%%%%%%%%%%%%%%%%%%%%%%%%%%%%%%
%\clearpage
\section{Additional Simulation Results}

\FloatBarrier
%The simulation results for Examples~(2) and~(3) are presented in
%Figures~\ref{EX2_loglog_plot_MSE_varyalpha} and
%\ref{EX3_loglog_plot_MSE_varyalpha}.
%\subsection{Examples (2)--(4): Sample size $n=1000$}
\begin{figure}[H]
\begin{subfigure}{0.48\textwidth} 
\includegraphics[width=0.85\textwidth,  height=0.21\textheight]{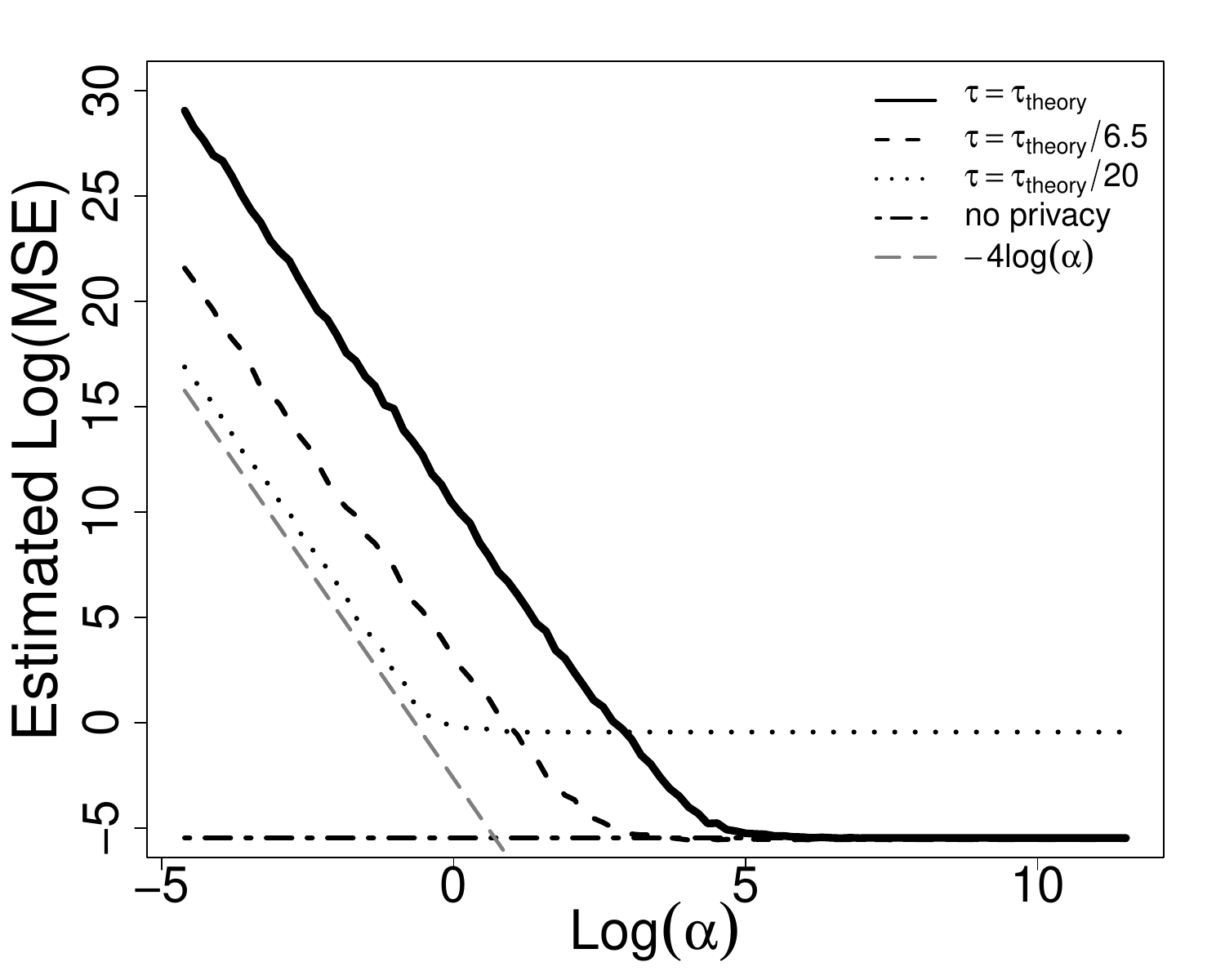}
\caption{Estimation of $\sigma_0{=}1.44$ from NI-privatized data.}
\label{EX2_loglog_plot_MSE_varyalpha_parta}
\end{subfigure}
\begin{subfigure}{0.48\textwidth}
\includegraphics[width=0.85\textwidth, height=0.21\textheight]{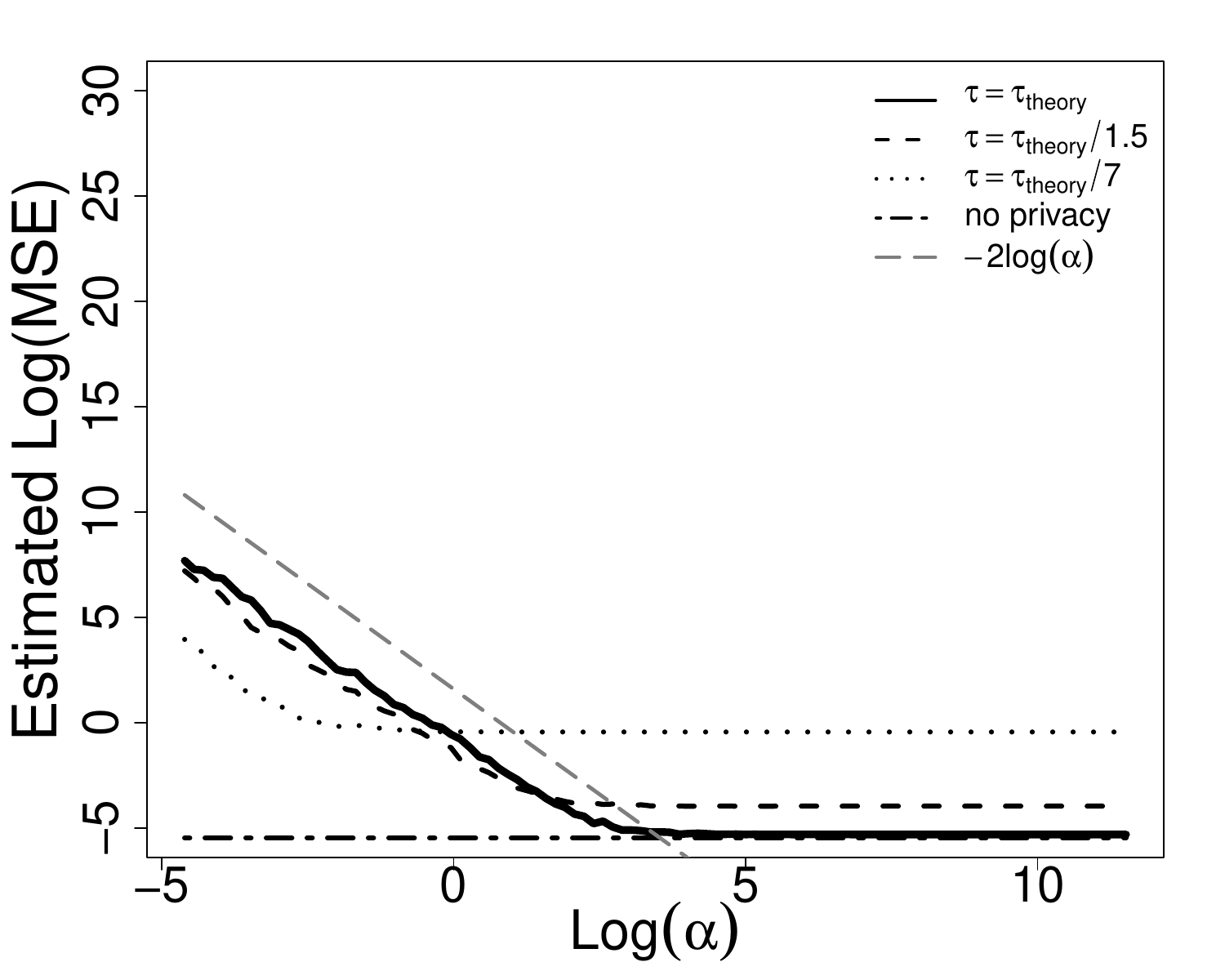}
\caption{Estimation of $\sigma_0{=}1.44$ from SI-privatized data.}
\end{subfigure}
\begin{subfigure}{0.48\textwidth}
\includegraphics[width=0.85\textwidth, height=0.21\textheight]{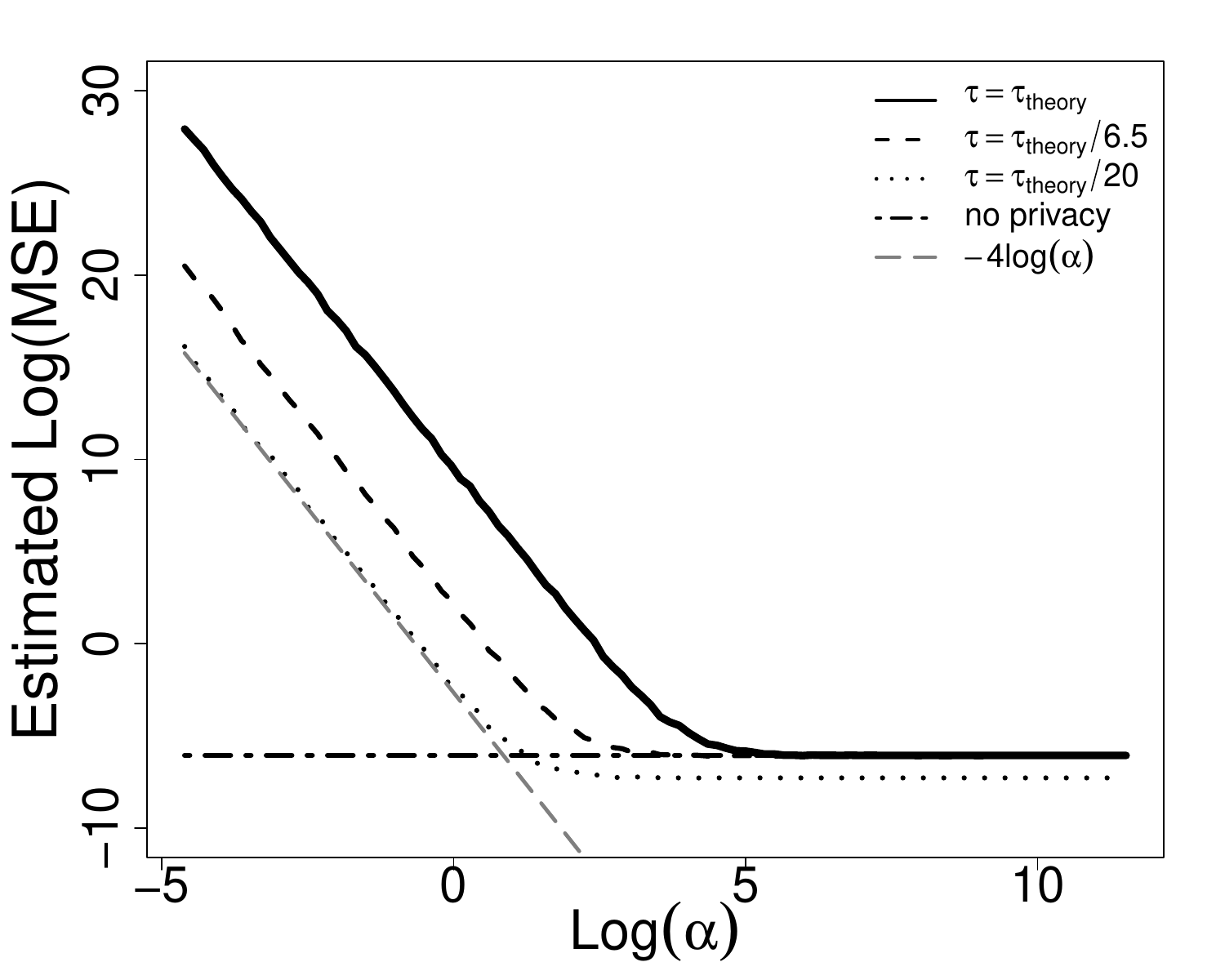} 
\caption{Estimation of $\sigma_2{=}0.02$ from NI-privatized data.}
\end{subfigure}
\begin{subfigure}{0.48\textwidth}
\includegraphics[width=0.85\textwidth, height=0.21\textheight]{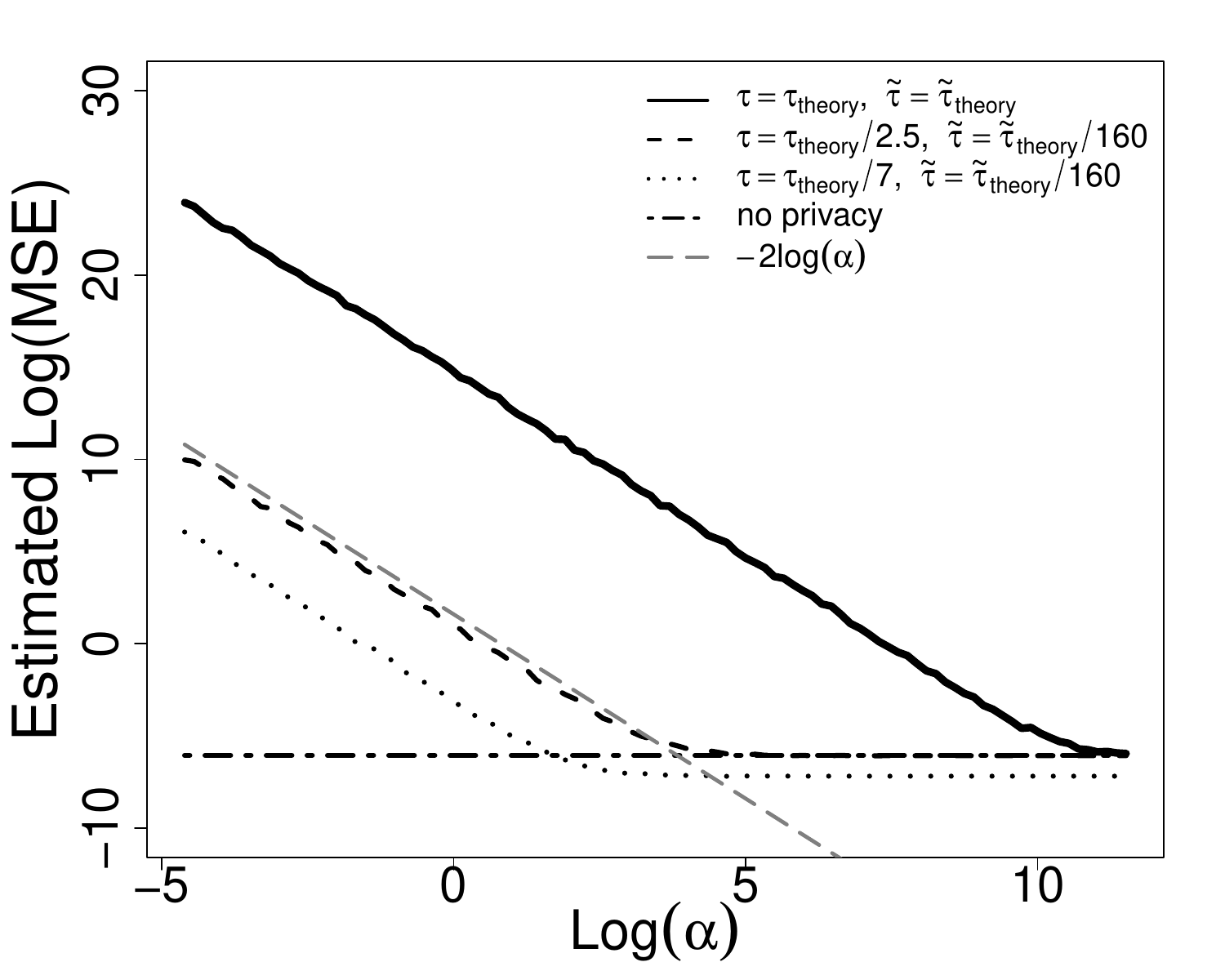} 
\caption{Estimation of $\sigma_2{=}0.02$ from SI-privatized data.}
\end{subfigure}
\begin{subfigure}{0.48\textwidth}
\includegraphics[width=0.85\textwidth, height=0.21\textheight]{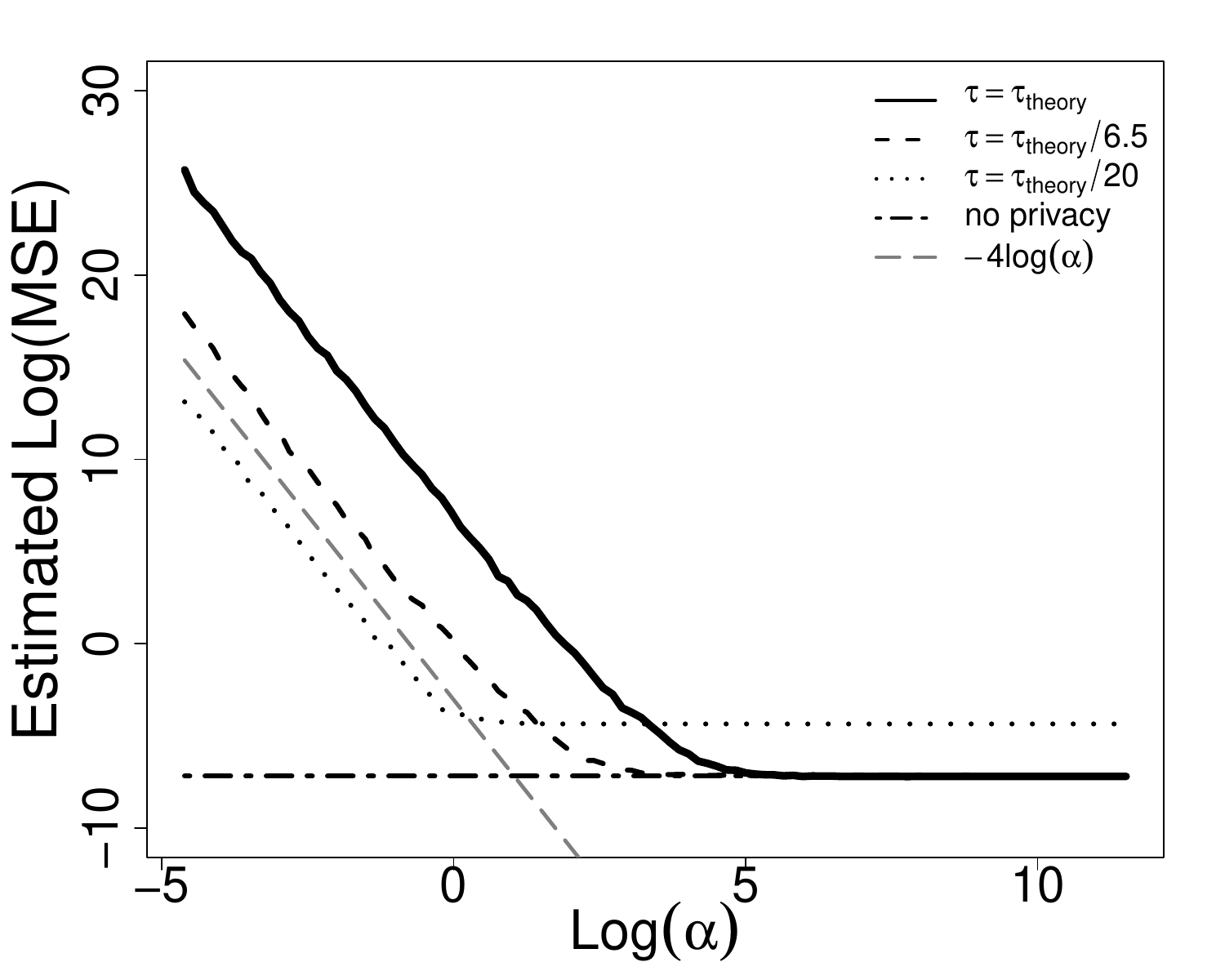} 
\caption{Estimation of $f(\pi/5){=} 0.19$ from NI-privatized data.}
\end{subfigure} \hspace*{.31cm}	
\begin{subfigure}{0.48\textwidth}
\includegraphics[width=0.85\textwidth, height=0.21\textheight]{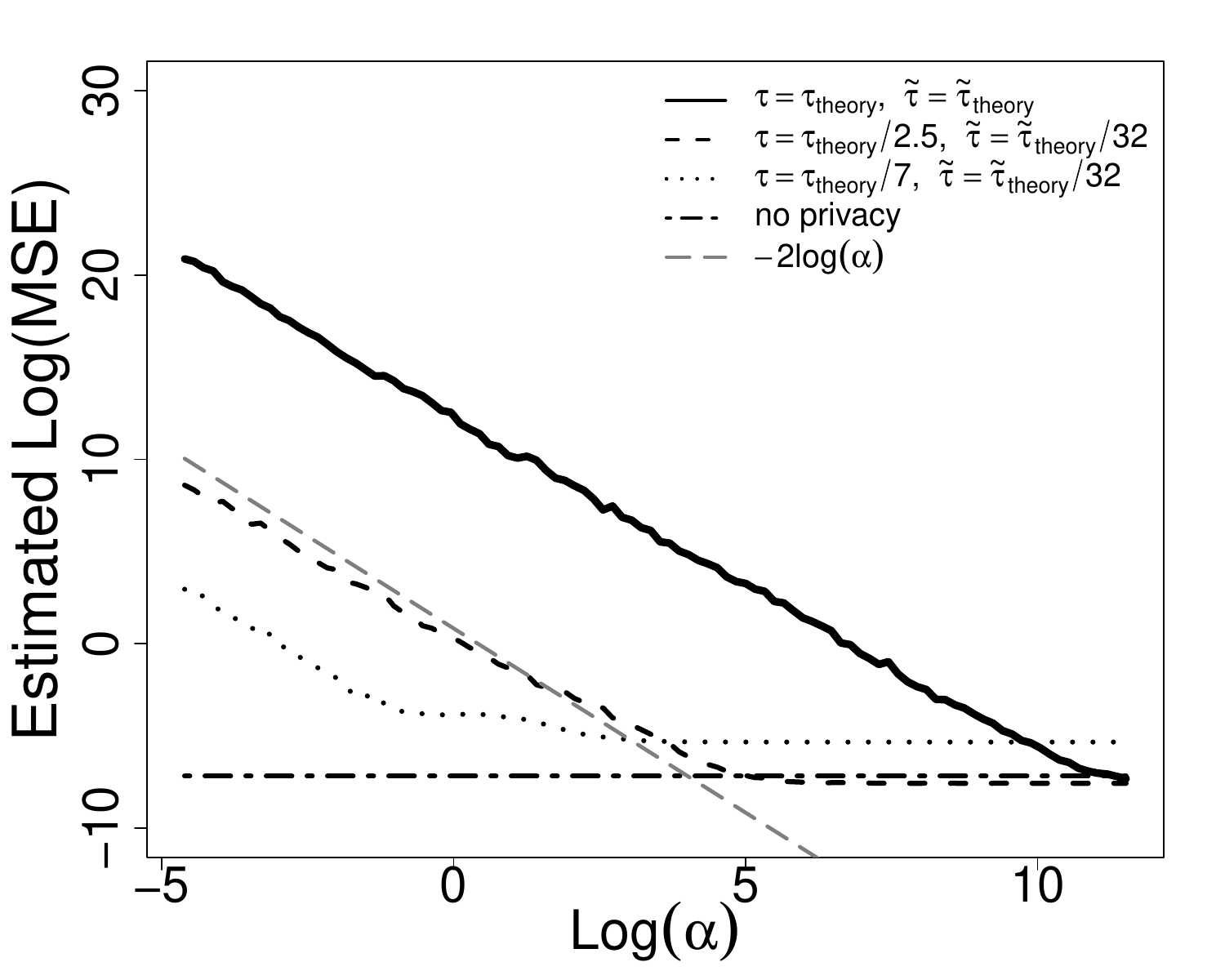} 
\caption{Estimation of $f(\pi/5){=} 0.19$  from SI-privatized data.}
\end{subfigure}
\caption{
Example~(2): Log-log plot of the estimated mean squared error (MSE) versus the privacy level $\alpha$ for estimation of  variance  (top), covariance coefficient (center) and  pointwise spectral density (bottom) using NI‑privatized data (left) and  SI‑privatized data (right).
Gray dashed lines indicate the theoretical $\alpha$‑scaling of the MSE for NI‑ and SI‑based estimators, including $\log(n)$-factors.
}
\label{EX2_loglog_plot_MSE_varyalpha}
\end{figure}
\begin{figure}[H]
\begin{subfigure}{0.48\textwidth} 
\includegraphics[width=0.85\textwidth, height=0.21\textheight]{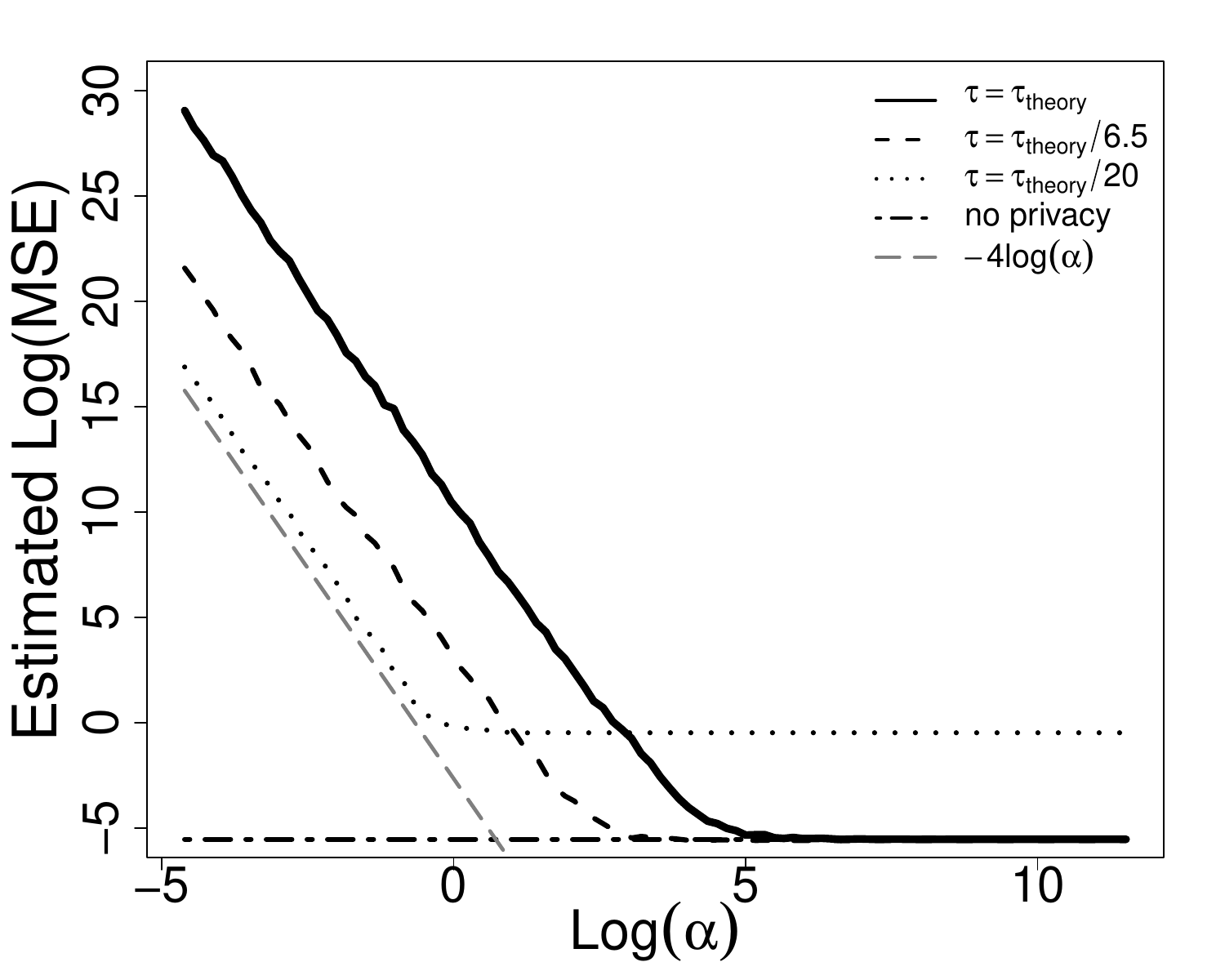}
\caption{Estimation of $\sigma_0{=}1.44$ from NI-privatized data.}
\label{EX3_loglog_plot_MSE_varyalpha_parta}
\end{subfigure}
\begin{subfigure}{0.48\textwidth}
\includegraphics[width=0.85\textwidth, height=0.21\textheight]{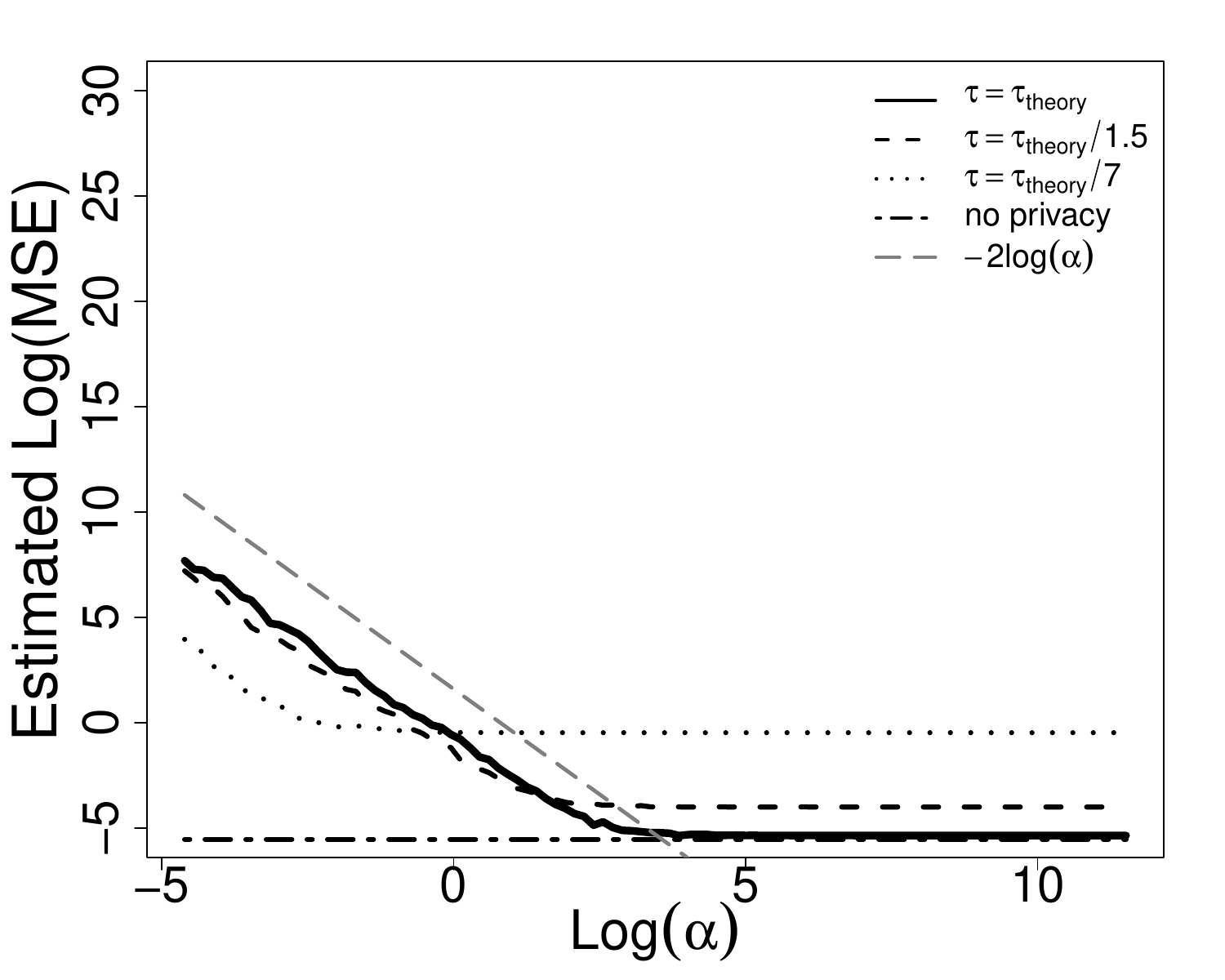}
\caption{Estimation of $\sigma_0{=}1.44$ from SI-privatized data.}
\end{subfigure}
\begin{subfigure}{0.48\textwidth}
\includegraphics[width=0.85\textwidth, height=0.21\textheight]{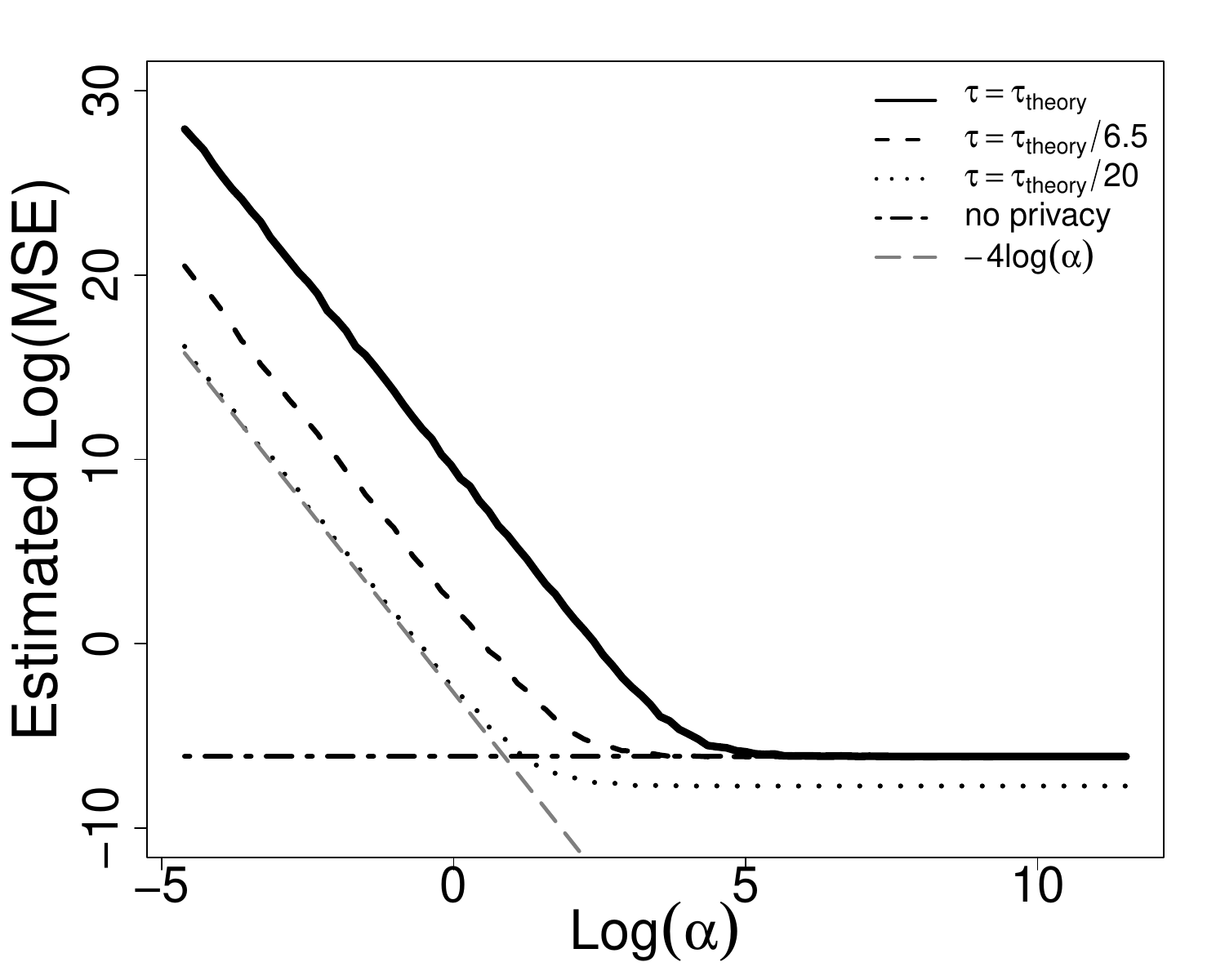} 
\caption{Estimation of $\sigma_2{=}0.005$ from NI-privatized data.}
\end{subfigure}
\begin{subfigure}{0.48\textwidth}
\includegraphics[width=0.85\textwidth, height=0.21\textheight]{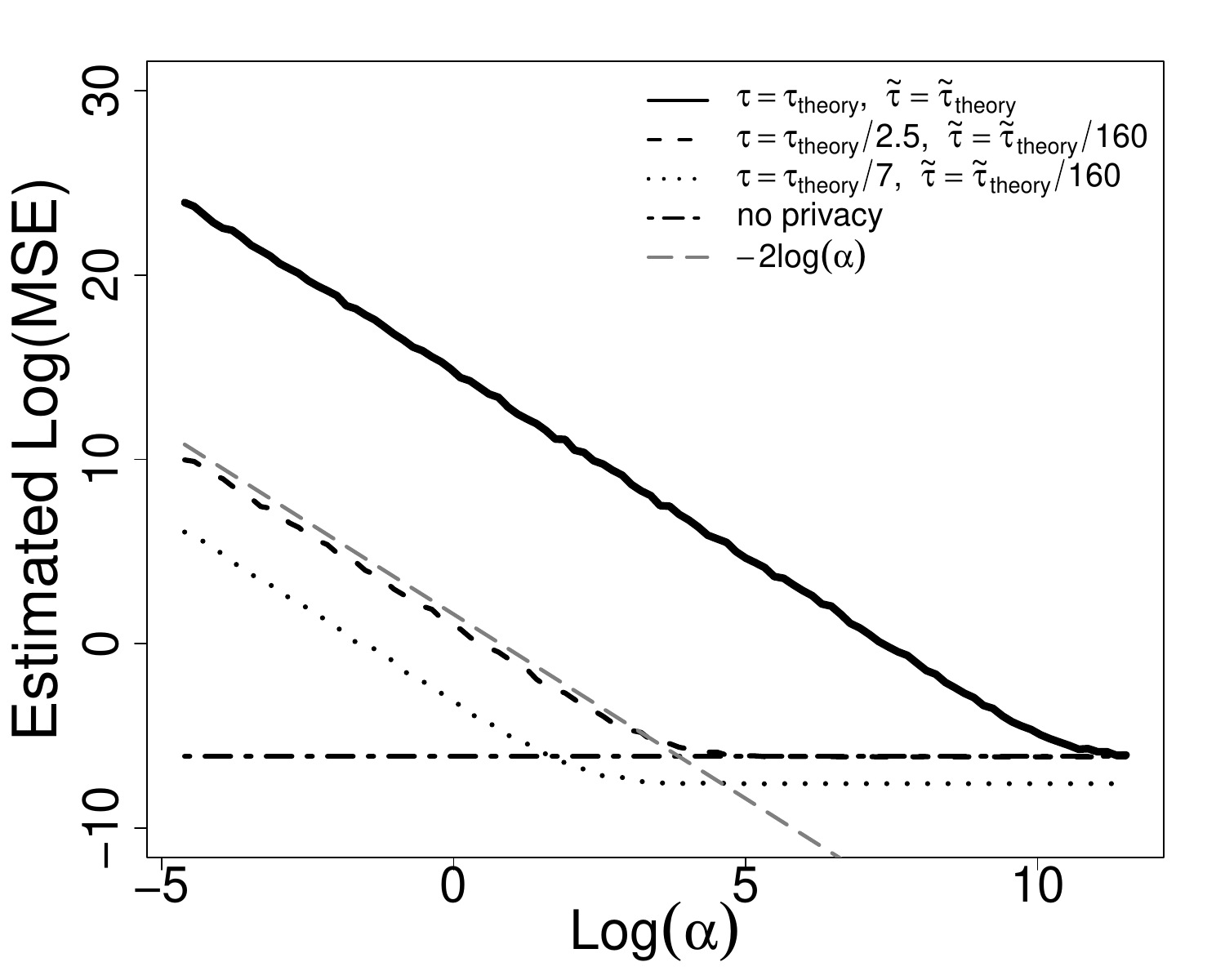} 
\caption{Estimation of $\sigma_2{=}0.005$ from SI-privatized data.}
\end{subfigure}
\begin{subfigure}{0.48\textwidth}
\includegraphics[width=0.85\textwidth, height=0.21\textheight]{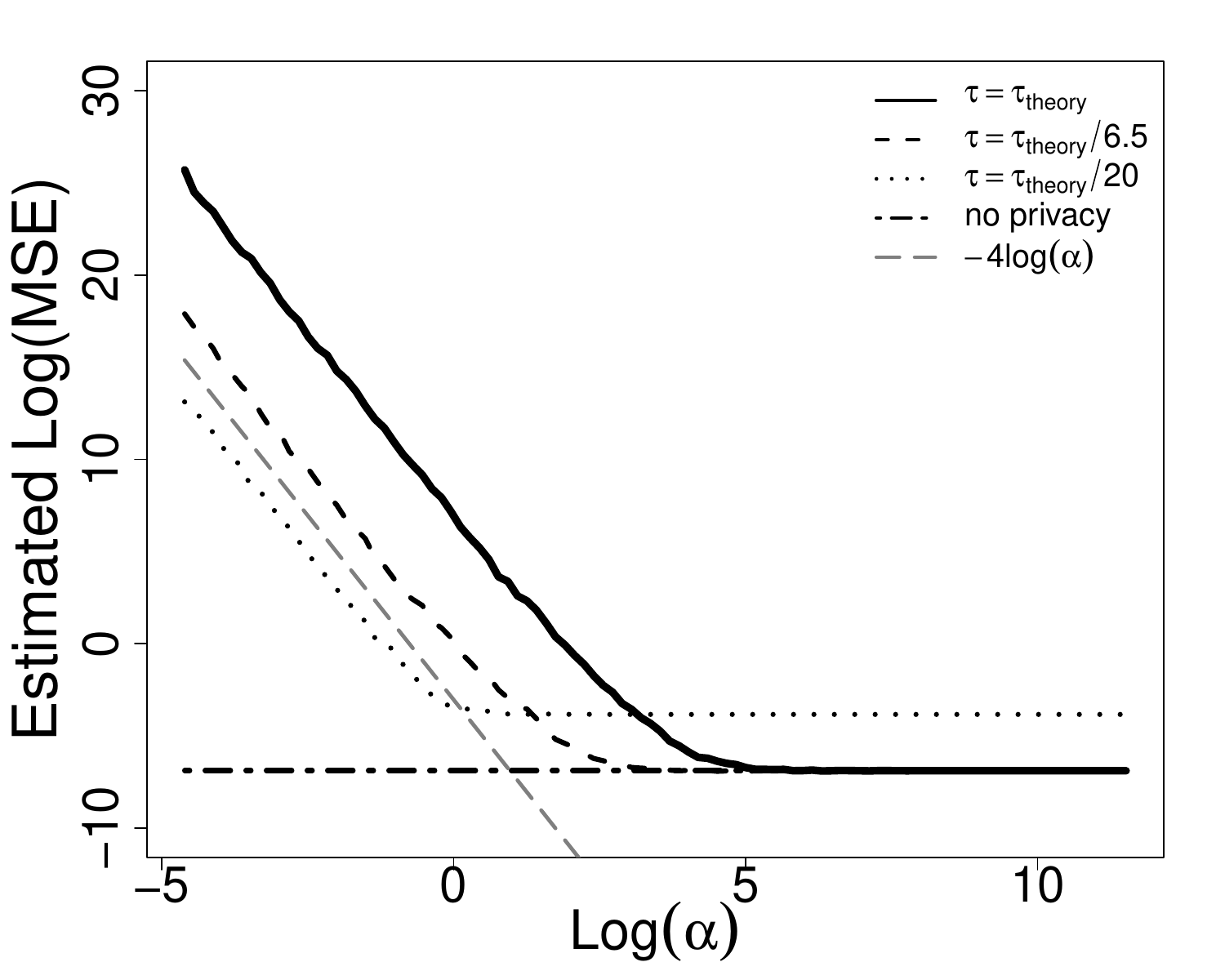} 
\caption{Estimation of $f(\pi/5){=} 0.24$ from NI-privatized data.}
\end{subfigure} \hspace*{.31cm}	
\begin{subfigure}{0.48\textwidth}
\includegraphics[width=0.85\textwidth, height=0.21\textheight]{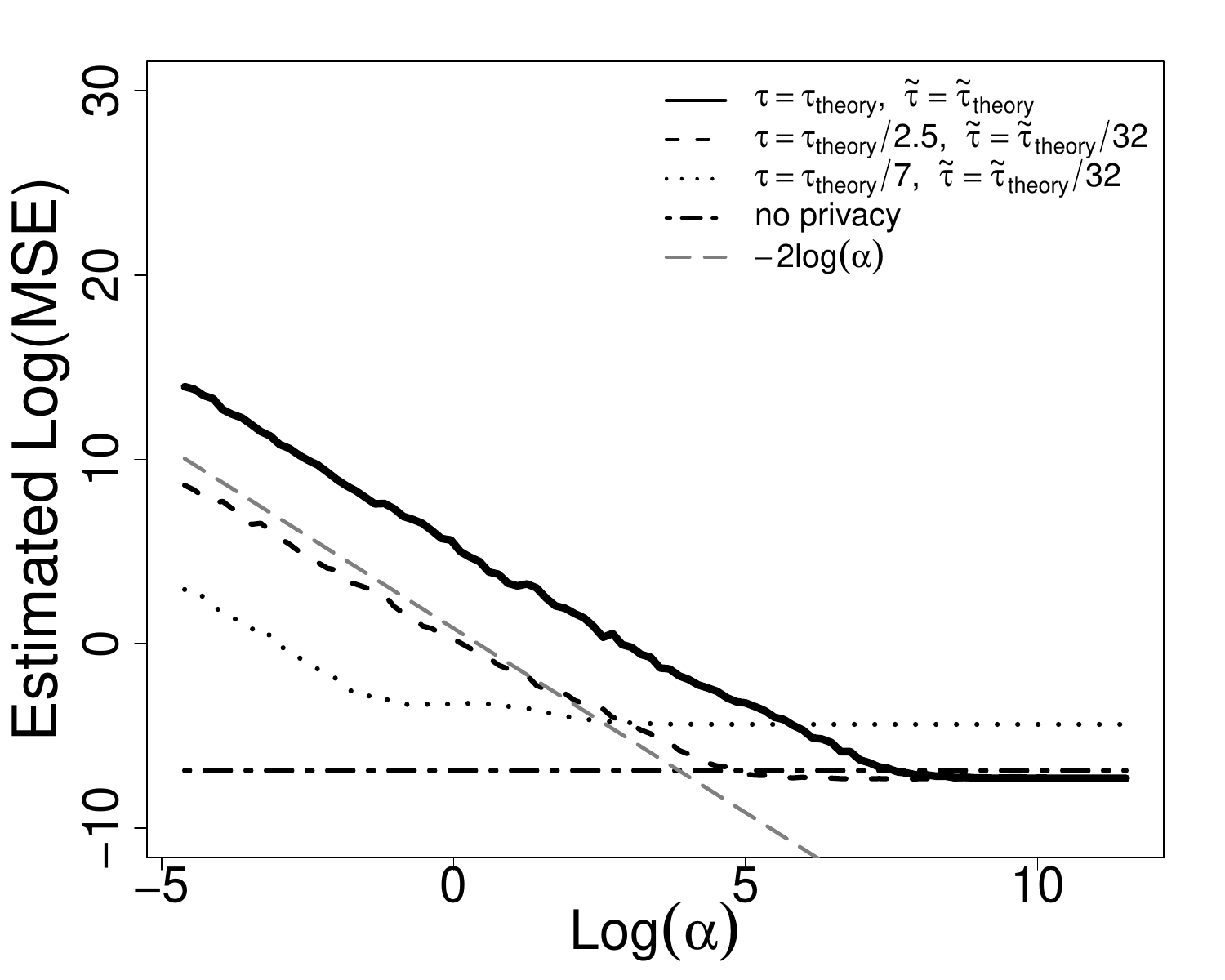} 
\caption{Estimation of $f(\pi/5){=} 0.24$  from SI-privatized data.}
\end{subfigure}
\caption{
Example~(3): Log-log plot of the estimated mean squared error (MSE) versus the privacy level $\alpha$ for estimation of  variance  (top), covariance coefficient (center) and  pointwise spectral density (bottom) using NI‑privatized data (left) and  SI‑privatized data (right).
Gray dashed lines indicate the theoretical $\alpha$‑scaling of the MSE for NI‑ and SI‑based estimators, including $\log(n)$-factors.
}
\label{EX3_loglog_plot_MSE_varyalpha}
\end{figure}
\FloatBarrier
We report the estimated mean squared error (MSE) for private estimation of $\sigma_0$, $\sigma_2$, and $f(\pi/5)$ from centered Gaussian stationary processes with a $0.8$‑Hölder continuous spectral density (Example~2) and with polynomially decaying covariances $\sigma_k = 1.44(1+|k|)^{-5.1}$ (Example~3), using both sequentially interactive (SI) and non‑interactive (NI) mechanisms. Tuning parameters and simulation settings are chosen as described in Section~\ref{sec:simulation}. The corresponding results are shown in Figure~\ref{EX2_loglog_plot_MSE_varyalpha} for Example~(2) and  Figure~\ref{EX3_loglog_plot_MSE_varyalpha} for Example~(3).

\bibliographystyle{apalike}
\bibliography{literature}

\end{document}